\documentclass[10pt]{article}
\usepackage{amssymb}

\newcommand{\bin}{\mathop{\mbox{bin}}\nolimits}

\title{The Busy Beaver Competition: a historical survey}
\author{Pascal MICHEL\thanks{Corresponding address: 59 rue du Cardinal Lemoine,
75005 Paris, France.}\\
{\small \'Equipe de Logique Math\'ematique, Institut de Math\'ematiques
de Jussieu--Paris Rive Gauche,}\\
{\small  UMR 7586, B\^atiment Sophie Germain, case 7012, 75205 Paris
Cedex 13, France}\\
{\small and Universit\'e de Cergy-Pontoise, ESPE, F-95000 Cergy-Pontoise, France}\\
{\small michel@math.univ-paris-diderot.fr}}

\date{Version 7\\ November 30, 2022}

\begin{document}
\maketitle

\begin{abstract}
Tibor Rado defined the Busy Beaver Competition in 1962. He used Turing machines
to give explicit definitions for some functions that are not computable
and grow faster than any computable function. He put forward the problem
of computing the values of these functions on numbers 1, 2, 3, $\ldots$.
More and more powerful computers have made possible the computation
of lower bounds for these values. In 1988, Brady extended the definitions
to functions on two variables.

We give a historical survey of these works. The successive record holders
in the Busy Beaver Competition are displayed, with their discoverers,
the date they were found, and, for some of them, an analysis of their
behavior.

We also survey the relations between busy beaver functions,
the variants of their definitions, and the links with
logical unprovability.

\bigskip


\bigskip

\end{abstract}

\section{Introduction}

\subsection{Noncomputable functions}

In 1936, Turing succeeded in making formal the intuitive notion of a
function computable by a finite, mechanical, procedure. He defined what
is now called a {\em Turing machine} and stated that a function on
integers is intuitively computable if and only if it is computable by a
Turing machine. Other authors, such as Church, Kleene, Post, and, later,
Markov, defined other models of computation that turn out to compute the
same functions as Turing machines do. See Soare (1996, 2007, 2009) for more
details about the history of the {\em Church-Turing Thesis}, as is now
named the capture of the intuitive notion of computability by the formal
notion of Turing machine.

Given a model of computation, a {\em noncomputable} function can
easily be defined by {\em diagonalization}. The list of all computable
functions is written, and then a function is defined such that it is
distinct from each function in the list. Then this function is
noncomputable. Such a definition by diagonalization leaves too much room
in the choice of the list and in the choice of the values of the final
function. What is needed is a function whose definition is simple,
natural and without ambiguousness.

In 1962, Rado succeeded in providing a natural definition for
noncomputable functions on the integers. He defined a {\em Busy Beaver}
game, leading to two functions $\mathnormal{\Sigma}$ and $S$ which are still the best
examples of noncomputable functions that one can give nowadays.
The values $\mathnormal{\Sigma}(n)$ and $S(n)$ are defined by considering the finite
set of carefully defined Turing machines with two symbols and $n$
states, and picking among these machines those with some maximal
behavior. It makes sense to compute the values $\mathnormal{\Sigma}(n)$, $S(n)$ of
these functions on small integers $n = 1,2,\ldots$. We have $\mathnormal{\Sigma}(1) =
S(1) = 1$, trivially. Lin and Rado (1965) gave proofs for the values
$\mathnormal{\Sigma}(2)$, $S(2)$, $\mathnormal{\Sigma}(3)$ and $S(3)$, and Brady (1983) did for $\mathnormal{\Sigma}(4)$
and $S(4)$. Only lower bounds had been provided for $\mathnormal{\Sigma}(5)$ and $S(5)$, by
the works of Green, Lynn, Schult, Uhing and eventually Marxen and
Buntrock. The lower bounds for $\mathnormal{\Sigma}(6)$ and $S(6)$ are still an
ongoing quest.

The initial Busy Beaver game, as defined by Rado, used Turing machines
with two symbols. Brady (1988) generalized the problem to Turing
machines with $k$ symbols, $k \ge 3$. He defined a function $S(n,k)$
such that $S(n,2)$ is Rado's $S(n)$, and gave some lower bounds.
Michel (2004) resumed the computation of lower bounds for $S(n,k)$ and
another function $\mathnormal{\Sigma}(n,k)$, and the search is going on, with the
works of Brady, Souris, Lafitte and Papazian, T. and S. Ligocki.

Since 2004, results are sent by email to Marxen and to Michel, who record
them on their websites. This paper aims to give a published version of
these records.

\subsection{Big numbers}

Consider Rado's functions $S$ and $\mathnormal{\Sigma}$.
Not only they are not computable, but they grow faster than any
computable function. That is, for any computable function $f$, there
exists an integer $N$ such that, for all $n > N$, we have $S(n) > f(n)$.
This property can be used to write big numbers. For example, if $S^k(n)$ denotes
$S(S(\ldots S(n) \ldots ))$, iterated $k$ times, then $S^{9^9}(9)$ is a very big number,
bigger than any number that was written with six symbols before the definition
of the $S$ function.

Bigger numbers can be obtained by defining functions growing much faster
than Rado's busy beaver functions. A natural idea to get such functions
is to define Turing machines of order $k$ as follows.
{\em Turing machines of order 1} are usual Turing machines without oracle,
and, for $k \ge 2$, {\em Turing machines of order $k$} are Turing machines
with oracle, where the oracle is the halting problem for Turing machines
of order $k-1$. Then the {\em $k$-th busy beaver function} $B_k(n)$
is the maximum number of steps taken by a Turing machine of order $k$
with $n$ states and two symbols that stops when it is launched on a blank tape.
So $B_1(n) = S(n)$, and $B_k(n)$ grows faster than any function computable
by a Turing machine of order $k$.

Unfortunately, there is no canonical way to define a Turing machine with oracle,
so Scott Aaronson, in his paper {\em Who can name the bigger number?}
(see the website), asked for naturally defined functions growing as fast as the
$k$-th busy beaver functions for $k \ge 2$.
Such functions were found by Nabutovsky and Weinberger (2007). By using
homology of groups, they defined a function growing as fast as
the third busy beaver function, and another one growing as fast as
the fifth busy beaver function. Michel (2010) went on studying these functions.
Scott Aaronson (2020) defined the {\em beeping busy beaver function},
which grows as fast as the second busy beaver function.
See the definition and values in Section 7.5.

\subsection{Contents}

The paper is structured as follows.
\begin{enumerate}
\item Introduction.
\item Preliminaries.
\item Historical overview.
\item Historical survey 
(lower bounds for $S(n,k)$ and $\mathnormal{\Sigma}(n,k)$, and
tables of the Turing machines that achieve these lower bounds).
\item Behaviors of busy beavers. We also display
the relations between these behaviors and open problems in mathematics called
Collatz-like problems and we resume some machines with non-Collatz-like
behaviors. We also present pairs of machines that have the same behaviors,
but not the same numbers of states and symbols.
\item Properties of the busy beaver functions and relations between
$S(n)$ and $\mathnormal{\Sigma}(n)$.
\item Variants of busy beavers:

 - Busy beavers defined by 4-tuples.

 - Busy beavers whose head can stand still.

 - Busy beavers on a one-way infinite tape.
 
 - Two-dimensional busy beavers.

 - Beeping busy beavers
\item The methods.
\item Busy beavers and unprovability.
\end{enumerate}

\section{Preliminaries}

There are many possible definitions for a Turing machine. We will follow
the conventions chosen by Rado (1962) in his definition of functions $\mathnormal{\Sigma}$
and $S$. A Turing machine has a tape, made of cells, infinite to the
left and to the right. On each cell a {\em symbol} is written. There is
a finite set $S = \{0,1,\ldots\}$ of symbols. The symbol 0 is the
{\em blank symbol}. A Turing machine has a tape head, which reads and
writes symbols on the tape, and can move in both {\em directions} left or
right, denoted by $L$ and $R$. A Turing machine has a finite set of
{\em states} $Q = \{A,B,\ldots\}$, plus a special state $H$, the
{\em halting state}. A Turing machine has a {\em next move function}
$$\delta : Q \times S \longrightarrow (S \times \{L,R\} \times Q) \cup
\{(1,R,H)\}.$$
If we have $\delta(q,a) = (b,d,p)$, then it means that, when the Turing
machine is in state $q$ and reads symbol $a$ on the tape, then it writes
symbol $b$ instead of $a$ on the cell currently read, it moves one cell in
the direction $d \in \{L,R\}$, and it changes the state from $q$ to $p$.
Each application of next move function $\delta$ is a {\em step} of the
computation. If $\delta(q,a) = (1,R,H)$, then, when the machine is in
state $q$ reading symbol $a$, it writes a 1, moves right, enters state
$H$, and stops. We follow Rado (1962) in not allowing the center
direction, that is in compelling the tape head to move left or right at
each step. Like Rado, we keep the halting state $H$ out of the set of
states. We differ from Rado in not allowing transitions
$\delta(q,a) = (b,d,H)$ with $b \ne 1$, $d \ne R$.

Note that such a machine is a universal model of computation. That is, any
computable function on integers can be computed by a Turing machine as
defined above. Initially, a finite string of symbols is written on the
tape. It is called the input, and can be a code for an integer. All other
cells contain the blank symbol. The tape head reads the leftmost symbol
of the input and the state is the initial state $A$. Then the
computation is launched according to the next move function. If it
stops, by entering the halting state $H$, then the string of symbols
written on the tape is the output, which can be a code for an integer. So a
Turing machine defines a partial function on integers. Reciprocally, any
computable partial function on integers can be computed by a Turing
machine as defined above.

In order to define functions $\mathnormal{\Sigma}$ and $S$, Rado (1962) considers
Turing machines with $n$ states and two symbols 0 and 1. His definitions
can be easily extended to Turing machines with $n$ states and $k$
symbols, $k \ge 3$, as Brady (1988) does. We consider the set $TM(n,k)$
of Turing machines with $n$ states and $k$ symbols. With our
definitions, it is a finite set with $(2kn+1)^{kn}$ members. We launch
each of these $(2kn+1)^{kn}$ Turing machines on a blank tape, that is a
tape with the blank symbol 0 in each cell. Some of these machines never
stop. The other ones, that eventually stop, are called {\em busy beavers},
and they are competing in two competitions, for the maximum number of
steps and for the maximum number of non-blank symbols left on the tape.
Let $s(M)$ be the number of computation steps taken by the busy beaver
$M$ to stop. Let $\sigma(M)$ be the number of non-blank symbols left on
the tape by the busy beaver $M$ when it stops. Then the busy beaver
functions are
$$S(n,k) =\max\{s(M):M \mbox{ is a busy beaver with $n$ states and $k$ symbols}\},$$
$$\mathnormal{\Sigma}(n,k) =\max\{\sigma(M):M \mbox{ is a busy beaver
with $n$ states and $k$ symbols}\}.$$
For $k = 2$, we find Rado's functions $S(n) = S(n,2)$ and
$\mathnormal{\Sigma}(n) = \mathnormal{\Sigma}(n,2)$.

Note that a permutation of the states, symbols or directions does not
change the behavior of a Turing machine. The choice between machines that
differ only by such permutations is settled by the following normalizing rule:
when a Turing machine is launched on a blank tape, it enters states in the
order $A$, $B$, $C,\ldots$, it writes symbols in the order 1, 2$\ldots$,
and it first moves right. So, normally, the first transition is
$\delta(A,0) = (1,R,B)$ or $\delta(A,0) = (0,R,B)$.

\bigskip

\noindent{\bf Note about terminology and notations}

Many names are used by authors: busy beaver game, busy beaver contest,
busy beaver problem, busy beaver competition.
All of them were first used early: game by Rado (1962), contest and
problem by Rado (1963), competition by Green (1964).

What is exactly a busy beaver is rarely specified. Let us give some
exceptions. For some authors, such as Green (1964) and Oberschelp et al.\
(1988), a busy beaver is any Turing machine that participates to the
busy beaver competition and halts. For others, such as Dewdney (1984)
and Ben-Amram and Petersen (2002), a busy beaver is a winner of this
competition. Rado (1962) called the winner a champion, and this term
has been used sometimes afterwards.

The number of ones left on the tape by the Turing machine $M$ when
it stops is often called the \emph{score} and denoted by
$\sigma(M)$, since Rado (1962). It is called the \emph{productivity}
by Boolos and Jeffrey (1974), a term used again by Hertel (2009) and
Harland (2013,2016). Harland uses the term \emph{activity}
for the number of moves of a Turing machine.

Almost all authors use the notations $\mathnormal{\Sigma(n)}$
and $S(n)$ for the busy beaver functions. Notable exceptions are:
ones ($n$) and time($n$) by Ben-Amram et al.\ (1996) and
Ben-Amram and Petersen (2002); $bb(n)$ and $ff(n)$
by Harland (2013,2016,2022).

\section{Historical overview}

The search for champions in the busy beaver competition can be roughly
divided into the following stages. Note that, from the beginnings,
computers have been tools to find good competitors, so better results follow
more powerful computers.

\bigskip

\begin{table}
\begin{tabular}{|c|c|l|}
\hline
1963 & Rado, Lin & \qquad $S(2,2)$ = 6, $\mathnormal{\Sigma}(2,2)$ = 4\\
     &           & \qquad $S(3,2)$ = 21, $\mathnormal{\Sigma}(3,2)$ = 6\\
\hline
1964 & Brady & (4,2)-TM: $s$ = 107, $\sigma$ = 13\\
\hline
1964 & Green & (5,2)-TM: $\sigma$ = 17\\
     &       & (6,2)-TM: $\sigma$ = 35\\
     &       & (7,2)-TM: $\sigma$ = 22,961\\
\hline
1972 & Lynn & (5,2)-TM: $s$ = 435, $\sigma$ = 22\\
     &      & (6,2)-TM: $s$ = 522, $\sigma$ = 42\\
\hline
1973 & Weimann & (5,2)-TM: $s$ = 556, $\sigma$ = 40\\
\hline
1974 & Lynn & (5,2)-TM: $s$ = 7,707, $\sigma$ = 112\\
\hline
1974 & Brady & \qquad $S(4,2)$ = 107, $\mathnormal{\Sigma}$(4,2) = 13\\
\hline
1983 & Brady & (6,2)-TM: $s$ = 13,488, $\sigma$ = 117\\
\hline
 1982 & Schult & (5,2)-TM: $s$ = 134,467, $\sigma$ = 501\\
      &        & (6,2)-TM: $s$ = 4,208,824, $\sigma$ = 2,075\\
\hline
December 1984 & Uhing & (5,2)-TM: $s$ = 2,133,492, $\sigma$ = 1,915\\
\hline
February 1986 & Uhing & (5,2)-TM: $s$ = 2,358,064\\
\hline
1988  & Brady & (2,3)-TM: $s$ = 38, $\sigma$ = 9\\
      &       & (2,4)-TM: $s$ = 7,195, $\sigma$ = 90\\
\hline
February 1990 & Marxen, Buntrock & {\bf (5,2)-TM}: $s$ = 47,176,870, $\sigma$ = 4,098\\
              &                  & (6,2)-TM: $s$ = 13,122,572,797, $\sigma$ = 136,612\\
\hline
September 1997 & Marxen, Buntrock & (6,2)-TM: $s$ =
 8,690,333,381,690,951, $\sigma$ = 95,524,079\\
\hline
August 2000 & Marxen, Buntrock & (6,2)-TM: $s > 5.3 \times 10^{42}$, $\sigma > 2.5 \times 10^{21}$\\
\hline
October 2000 & Marxen, Buntrock & (6,2)-TM: $s > 6.1 \times 10^{925}$, $\sigma > 6.4 \times 10^{462}$\\
\hline
March 2001 & Marxen, Buntrock & (6,2)-TM: $s > 3.0 \times 10^{1730}$, $\sigma > 1.2 \times 10^{865}$\\
\hline
\end{tabular}
\caption{Busy Beaver Competition from 1963 to 2001. In the last column, an ($n$,$k$)-Turing machine
is a Turing machine with $n$ states and $k$ symbols. Number $s$ is the number of steps,
and number $\sigma$ is the number of non-blank symbols left by the Turing machine when it stops.
When ($n$,$k$)-TM is in bold type, the Turing machine is the current record holder.
When values of $S(n,k)$ and $\mathnormal{\Sigma}(n,k)$ are indicated, the line refers to the proof
that the functions have these values.}
\end{table}

\begin{table}
\begin{tabular}{|c|c|l|}
\hline
October 2004  & Michel & (3,3)-TM: $s$ = 40,737, $\sigma$ = 208\\
\hline
November 2004 & Brady & (3,3)-TM: $s$ = 29,403,894, $\sigma$ = 5,600\\
\hline
December 2004 & Brady & (3,3)-TM: $s$ = 92,649,163, $\sigma$ = 13,949\\
\hline
February 2005 & T. and S. Ligocki & {\bf (2,4)-TM}: $s$ = 3,932,964, $\sigma$ = 2,050\\
              &                   & (2,5)-TM: $s$ = 16,268,767, $\sigma$ = 4,099\\
              &                   & (2,6)-TM: $s$ = 98,364,599, $\sigma$ = 10,574\\

\hline
April 2005 & T. and S. Ligocki & (4,3)-TM: $s$ = 250,096,776, $\sigma$ = 15,008\\
 & & (3,4)-TM: $s$ = 262,759,288, $\sigma$ = 17,323\\
 & & (2,5)-TM: $s$ = 148,304,214, $\sigma$ = 11,120\\
 & & (2,6)-TM: $s$ = 493,600,387, $\sigma$ = 15,828\\
\hline
July 2005 & Souris & (3,3)-TM: $s$ = 544,884,219, $\sigma$ = 36,089\\
\hline
August 2005 & Lafitte, Papazian & (3,3)-TM: $s$ = 4,939,345,068, $\sigma$ = 107,900\\
 & & (2,5)-TM: $s$ = 8,619,024,596, $\sigma$ = 90,604\\
\hline
September 2005 & Lafitte, Papazian & (3,3)-TM: $s$ = 987,522,842,126, $\sigma$ = 1,525,688\\
 & & (2,5)-TM: $\sigma$ = 97,104\\
\hline
October 2005 & Lafitte, Papazian & (2,5)-TM: $s$ = 233,431,192,481, $\sigma$ = 458,357\\
 & & (2,5)-TM: $s$ = 912,594,733,606, $\sigma$ = 1,957,771\\
\hline
December 2005 & Lafitte, Papazian & (2,5)-TM: $s$ = 924,180,005,181\\
\hline
April 2006 & Lafitte, Papazian & (3,3)-TM: $s$ = 4,144,465,135,614, $\sigma$ = 2,950,149\\
\hline
May 2006 & Lafitte, Papazian & (2,5)-TM: $s$ = 3,793,261,759,791, $\sigma$ = 2,576,467\\
\hline
June 2006 & Lafitte, Papazian & (2,5)-TM: $s$ = 14,103,258,269,249, $\sigma$ = 4,848,239\\
\hline
July 2006 & Lafitte, Papazian & (2,5)-TM: $s$ = 26,375,397,569,930\\
\hline
August 2006 & T. and S. Ligocki & (3,3)-TM: $s$ = 4,345,166,620,336,565,  $\sigma$ = 95,524,079\\
 & & (2,5)-TM: $s > 7.0 \times 10^{21}$, $\sigma$ = 172,312,766,455\\
\hline
\end{tabular}
\caption{Busy Beaver Competition from 2004 to 2006}
\end{table}

\begin{table}
\begin{tabular}{|c|c|l|}
\hline
June 2007 & Lafitte, Papazian
   & \qquad $S(2,3)$ = 38, $\mathnormal{\Sigma}(2,3)$ = 9\\
\hline
September 2007 & T. and S. Ligocki & (3,4)-TM: $s > 5.7 \times 10^{52}$, $\sigma > 2.4 \times 10^{26}$\\
 & & (2,6)-TM: $s > 2.3 \times 10^{54}$, $\sigma > 1.9 \times 10^{27}$\\
\hline
October 2007 & T. and S. Ligocki
   & (4,3)-TM: $s > 1.5 \times 10^{1426}$, $\sigma > 1.1 \times 10^{713}$\\
 & & (3,4)-TM: $s > 4.3 \times 10^{281}$, $\sigma > 6.0 \times 10^{140}$\\
 & & (3,4)-TM: $s > 7.6 \times 10^{868}$, $\sigma > 4.6 \times 10^{434}$\\
 & & (3,4)-TM: $s > 3.1 \times 10^{1256}$, $\sigma > 2.1 \times 10^{628}$\\
 & & (2,5)-TM: $s > 5.2 \times 10^{61}$, $\sigma > 9.3 \times 10^{30}$\\
 & & (2,5)-TM: $s > 1.6 \times 10^{211}$, $\sigma > 5.2 \times 10^{105}$\\
\hline
November 2007 & T. and S. Ligocki
   & (6,2)-TM: $s > 8.9 \times 10^{1762}$, $\sigma > 2.5 \times 10^{881}$\\
 & & {\bf (3,3)-TM}: $s$ = 119,112,334,170,342,540, $\sigma$ = 374,676,383\\
 & & (4,3)-TM: $s > 7.7 \times 10^{1618}$, $\sigma > 1.6 \times 10^{809}$\\
 & & (4,3)-TM: $s > 3.7 \times 10^{1973}$, $\sigma > 8.0 \times 10^{986}$\\
 & & (4,3)-TM: $s > 3.9 \times 10^{7721}$, $\sigma > 4.0 \times 10^{3860}$\\
 & & (4,3)-TM: $s > 3.9 \times 10^{9122}$, $\sigma > 2.5 \times 10^{4561}$\\
 & & (3,4)-TM: $s > 8.4 \times 10^{2601}$, $\sigma > 1.7 \times 10^{1301}$\\
 & & (3,4)-TM: $s > 3.4 \times 10^{4710}$, $\sigma > 1.4 \times 10^{2355}$\\
 & & (3,4)-TM: $s > 5.9 \times 10^{4744}$, $\sigma > 2.2 \times 10^{2372}$\\
 & & {\bf (2,5)-TM}: $s > 1.9 \times 10^{704}$, $\sigma > 1.7 \times 10^{352}$\\
 & & (2,6)-TM: $s > 4.9 \times 10^{1643}$, $\sigma > 8.6 \times 10^{821}$\\
 & & (2,6)-TM: $s > 2.5 \times 10^{9863}$, $\sigma > 6.9 \times 10^{4931}$\\
\hline
December 2007 & T. and S. Ligocki
   & (6,2)-TM: $s > 2.5 \times 10^{2879}$, $\sigma > 4.6 \times 10^{1439}$\\
 & & (4,3)-TM: $s > 7.9 \times 10^{9863}$, $\sigma > 8.9 \times 10^{4931}$\\
 & & (4,3)-TM: $s > 5.3 \times 10^{12068}$, $\sigma > 4.2 \times 10^{6034}$\\
 & & {\bf (3,4)-TM}: $s > 5.2 \times 10^{13036}$, $\sigma > 3.7 \times 10^{6518}$\\
\hline
January 2008 & T. and S. Ligocki
   & {\bf (4,3)-TM}: $s > 1.0 \times 10^{14072}$, $\sigma > 1.3 \times 10^{7036}$\\
 & & {\bf (2,6)-TM}: $s > 2.4 \times 10^{9866}$, $\sigma > 1.9 \times 10^{4933}$\\
\hline
May 2010 & Kropitz
   & (6,2)-TM: $s > 3.8 \times 10^{21132}$, $\sigma > 3.1 \times 10^{10566}$\\
\hline
June 2010 & Kropitz
   & (6,2)-TM: $s > 7.4 \times 10^{36534}$, $\sigma > 3.4 \times 10^{18267}$\\
\hline
March 2014 & ``Wythagoras''
   & (7,2)-TM: $s > \sigma > 10^{10^{10^{10^{18,705,352}}}}$\\
\hline
May 2022 & S. Ligocki
   & (6,2)-TM: $s > 9.6 \times 10^{78913}$, $\sigma > 6.0 \times 10^{39456}$\\
\hline
May 2022 & Kropitz
   & (6,2)-TM: $s > 5.4 \times 10^{197282}$, $\sigma > 2.0 \times 10^{98641}$\\
 & & (6,2)-TM: $s > 8.2 \times 10^{1,292,913,985}$, $\sigma > 1.7 \times 10^{646,456,993}$\\
\hline
May 2022 & S. Ligocki
   & (6,2)-TM: $s > \sigma > 10^{10^{10^{10^{20823}}}}$\\
\hline
May 2022 & Kropitz
  & {\bf (6,2)-TM}: $s > \sigma > 10^{\wedge \wedge}15$\\
\hline
\end{tabular}
\caption{Busy Beaver Competition since 2007}
\end{table}

\noindent{\bf First stage: Following the definitions}.
The definitions of the busy beaver functions $\mathnormal{\Sigma}(n)$ and $S(n)$ by
Rado (1962) were quickly followed by conjectures and proofs for
$n = 2,3$, by Rado and Lin. Brady (1964) gave a conjecture for $n = 4$,
and Green (1964) gave lower bounds for many values of $n$. Lynn (1972)
improved these lower bounds for $n = 5,6$. Brady proved his conjecture
for $n = 4$ in 1974, and published the result in 1983. Details on this
first stage can be found in the articles of Lynn (1972) and Brady
(1983, 1988).

\bigskip

\noindent{\bf Second stage: Following the Dortmund contest}.
More results for $n = 5,6$ followed the contest that was organized at
Dortmund in 1983, and was wun by Schult. Uhing improved twice the result
in 1984 and in 1986. Marxen and Buntrock began a search for competitors
for $n = 5,6$ in 1989. They quickly found a conjectural winner for
$n = 5$, and went on finding many good machines for $n = 6$, up to 2001.
Michel (1993) studied the behaviors of many competitors for $n = 5$,
proving that they depend on well known open problems in number theory.
Details on this second stage can be found in the articles of Dewdney
(1984ab,1985ab), Brady (1988), and Marxen and Buntrock (1990).
From 1997, results began to be put on the web, either on Google groups,
or on personal websites.

\bigskip

\noindent{\bf Third stage: Machines with more than two symbols}.
As soon as 1988, Brady extended the busy beaver competition to machines
with more than two symbols and gave some lower bounds. Michel (2004)
resumed the search, and his lower bounds were quickly overtaken by those
from Brady. Between 2005 and 2008, more than forty new machines, each one
breaking a record, were found by two teams: the French one made of
Gr\'egory Lafitte and Christophe Papazian, and the father-and-son
collaboration of Terry and Shawn Ligocki.
Four new machines for the classical busy beaver competition of
machines with 6 states and 2 symbols were also found, by the Ligockis
and by Pavel Kropitz.
In May 2022, Shawn Ligocki and Pavel Kropitz found many new machines
with 6 states and 2 symbols.

With the coming of the web
age, researchers have faced two problems: how to announce results, and how to
store them. In 1997, Heiner Marxen chose to post them on Google groups,
but it seems that the oldest reports are no longer available. From 2004,
most results have been announced by sending them by email to several
people (for example, the new machines with 6 states and 2 symbols found
by Terry and Shawn Ligocki in November and December 2007 were sent by
email to six persons: Allen H. Brady, Gr\'egory Lafitte, Heiner Marxen,
Pascal Michel, Christophe Papazian and Myron P. Souris).
Storing results have been made on web pages (see websites list after the
references). Brady has stored results on machines with 3 states and 3
symbols on his own website. Both Marxen and Michel have kept account of
all results on their websites. Moreover, Marxen has held simulations,
with four variants, of each discovered machine. Michel has held
theoretical analyses of many machines.
In February 2022, Shawn Logocki created a Google Group email list:

\verb+https://groups.google.com/g/busy-beaver-discuss+

\noindent on which many results have been posted.

\section{Historical survey}
\subsection{Turing machines with 2 states and 2 symbols}

\begin{itemize}

\item
 Rado (1963) claimed that $\mathnormal{\Sigma}(2,2) = 4$, but that $S(2,2)$
 was yet unknown.

\item
 The value $S(2,2) = 6$ was probably set by Lin in 1963.
 See 

\verb+http://turbotm.de/~heiner/BB/simTM22_bb.html+

 for a study of the winner by H. Marxen.
\end{itemize}

\begin{center}
\begin{tabular}{|c|c|c|c|}
\hline
1963 & Rado, Lin & $S(2,2) = 6$ & $\mathnormal{\Sigma}(2,2) = 4$\\
\hline
\end{tabular}
\end{center}

\bigskip

 The winner and some other good machines:

\smallskip

\begin{tabular}{cccccc}
  A0  & A1  & B0  & B1 &$s(M)$&$\sigma(M)$\\
  1RB & 1LB & 1LA & 1RH  &  6  &     4    \\
  1RB & 1RH & 1LB & 1LA  &  6  &     3    \\
  1RB & 0LB & 1LA & 1RH  &  6  &     3    
\end{tabular}

\subsection{Turing machines with 3 states and 2 symbols}

\begin{itemize}

\item
 Soon after the definition of the functions $S$ and $\mathnormal{\Sigma}$,
 by Rado (1962), it was conjectured that $S(3,2) = 21$,
 and $\mathnormal{\Sigma}(3,2) = 6$.

\item
 Lin (1963) proved this conjecture and this proof was
 eventually published by Lin and Rado (1965).
 See studies by Heiner Marxen of the winners for
 the $S$ function in

\verb+http://turbotm.de/~heiner/BB/simTM32_bbS.html+

 and for the $\mathnormal{\Sigma}$ function in

\verb+http://turbotm.de/~heiner/BB/simTM32_bbO.html+

\end{itemize}

\begin{center}
\begin{tabular}{|c|c|c|c|}
\hline
1963 & Rado, Lin & $S(3,2) = 21$ & $\mathnormal{\Sigma}(3,2) = 6$\\
\hline
\end{tabular}
\end{center}

\bigskip

 The winners and some other good machines:

\smallskip

\begin{tabular}{cccccccc}
 A0  &  A1 &  B0 &  B1 &  C0 &  C1 &$s(M)$&$\sigma(M)$\\
 1RB & 1RH & 1LB & 0RC & 1LC & 1LA  &  21 &     5   \\
 1RB & 1RH & 0LC & 0RC & 1LC & 1LA  &  20 &     5   \\
 1RB & 1LA & 0RC & 1RH & 1LC & 0LA  &  20 &     5   \\
 0RB & 1RH & 0LC & 1RA & 1RB & 1LC  &  17 &     4   \\
 0RB & 1LC & 1LA & 1RB & 1LB & 1RH  &  16 &     5   \\
 1RB & 1RH & 0RC & 1RB & 1LC & 1LA  &  14 &     6   \\
 1RB & 1RC & 1LC & 1RH & 1RA & 0LB  &  13 &     6   \\
 1RB & 1LC & 1LA & 1RB & 1LB & 1RH  &  13 &     6   \\
 0RB & 1LC & 1RC & 1RB & 1LA & 1RH  &  13 &     5   \\
 1RB & 1RA & 1LC & 1RH & 1RA & 1LB  &  12 &     6   \\
 1RB & 1LC & 1RC & 1RH & 1LA & 0LB  &  11 &     6   
\end{tabular}

\subsection{Turing machines with 4 states and 2 symbols}

\begin{itemize}

 \item
 Brady (1964,1965,1966) found a machine $M$ such
 that $s(M) = 107$ and $\sigma(M) = 13$.
 See study by H. Marxen in

\verb+http://turbotm.de/~heiner/BB/simTM42_bb.html+

 Brady conjectured that
 $S(4,2) = 107$ and $\mathnormal{\Sigma}(4,2) = 13$.

 \item
 Brady (1974,1975) proved this conjecture, and the proof
 was eventually published in Brady (1983).

 \item
 Independently, Machlin and Stout (1990) published
 another proof of the same result, first reported by
 Kopp (1981) (Kopp is the maiden name of Machlin).

\item
  Independently, Weimann, Casper and Fenzl (1973)
  claimed that they proved this conjecture.
 
\end{itemize} 

\begin{center}
\begin{tabular}{|c|c|c|c|}
\hline
1964 & Brady & $s$ = 107 & $\sigma$ = 13\\
\hline
1974 & Brady & $S(4,2) = 107$ & $\mathnormal{\Sigma}(4,2) = 13$\\
\hline
\end{tabular}
\end{center}

\bigskip

 The winner and some other good machines:

\smallskip

\begin{tabular}{cccccccccc} 
 A0  & A1  & B0  & B1  & C0  & C1  & D0  & D1  & $s(M)$&$\sigma(M)$\\
 1RB & 1LB & 1LA & 0LC & 1RH & 1LD & 1RD & 0RA  &  107 &    13    \\ 
 1RB & 1LD & 1LC & 0RB & 1RA & 1LA & 1RH & 0LC  &   97 &     9    \\ 
 1RB & 0RC & 1LA & 1RA & 1RH & 1RD & 1LD & 0LB  &   96 &    13    \\ 
 1RB & 1LB & 0LC & 0RD & 1RH & 1LA & 1RA & 0LA  &   96 &     6    \\ 
 1RB & 1LD & 0LC & 0RC & 1LC & 1LA & 1RH & 0LA  &   84 &    11    \\ 
 1RB & 1RH & 1LC & 0RD & 1LA & 1LB & 0LC & 1RD  &   83 &     8    \\ 
 1RB & 0RD & 1LC & 0LA & 1RA & 1LB & 1RH & 0RC  &   78 &    12     
\end{tabular}

\subsection{Turing machines with 5 states and 2 symbols} 

\begin{itemize}

 \item
 Green (1964) found a machine $M$ with $\sigma(M) = 17$.

 \item
 Lynn (1972) found machines $M$ and $N$ with $s(M) = 435$
 and $\sigma(N) = 22$.

\item
  Weimann (1973) found a machine $M$ with $s(M) = 556$
  and $\sigma(M) = 40$. 

 \item
 Lynn, cited by Brady (1983), found in 1974 machines $M$
 and $N$ with $s(M) = 7,707$ and $\sigma(N) = 112$.

 \item
   Uwe Schult, cited by Ludewig et al.\ (1983) and by Dewdney (1984a),
   found, in August 1982, a machine $M$ with $s(M) = 134,467$
   and $\sigma(M) = 501$.
   This machine was analyzed independently by Ludewig (in
   Ludewig et al. (1983)), by Robinson
   (cited by Dewdney (1984b)), and by Michel (1993).

 \item
 George Uhing, cited by Dewdney (1985a,b), found, in December 1984,
 a machine $M$ with $s(M) = 2,133,492$ and $\sigma(M) = 1,915$.
 This machine was analyzed by Michel (1993).

 \item
 George Uhing, cited by Brady (1988), found, in February 1986,
 a machine $M$ with $s(M) = 2,358,064$ (and $\sigma(M) = 1,471$).
 This machine was analyzed by Michel (1993). Machine 7 in
 Marxen bb-list, in

\verb+http://turbotm.de/~heiner/BB/bb-list+

 can be obtained from Uhing's one,
 as given by Brady (1988), by the permutation of states
 (A D B E). See study by H. Marxen in

\verb+http://turbotm.de/~heiner/BB/simmbL5_7.html+

 \item
 Heiner Marxen and J\"{u}rgen Buntrock found, in August 1989,
 a machine $M$ with $s(M) = 11,798,826$ and $\sigma(M) = 4,098$.
 This machine was cited by Marxen and Buntrock (1990),
 and by Machlin and Stout (1990),
 and was analyzed by Michel (1993).
 See study by H. Marxen in

\verb+http://turbotm.de/~heiner/BB/simmbL5_2.html+

 \item
 Heiner Marxen and J\"{u}rgen Buntrock found, in September 1989,
 a machine $M$ with $s(M) = 23,554,764$ (and $\sigma(M) = 4,097$).
 This machine was cited by Machlin and Stout (1990),
 and was analyzed by Michel (1993).
 See study by H. Marxen in

\verb+http://turbotm.de/~heiner/BB/simmbL5_3.html+

 and analysis by P. Michel in Section \ref{sec:tm52b}.

 \item
 Heiner Marxen and J\"{u}rgen Buntrock found, in September 1989,
 a machine $M$ with  $s(M) = 47,176,870$ and $\sigma(M) = 4,098$.
 This machine was cited by Marxen and Buntrock (1990),
 and was analyzed by Buro (1990) and by Michel (1993).
 See study by H. Marxen in

\verb+http://turbotm.de/~heiner/BB/simmbL5_1.html+

analysis by Buro in (p.\ 64-67)

\verb+https://skatgame.net/mburo/ps/diploma.pdf+

 and analysis by P. Michel in Section \ref{sec:tm52a}.
 It is the current record holder.

 \item
 Marxen gives a list of machines $M$ with high values of $s(M)$ and $\sigma(M)$ in

\verb+http://turbotm.de/~heiner/BB/bb-list+

 \item
 The study of Turing machines with 5 states and 2 symbols is still going on.
 Marxen and Buntrock (1990),
 Skelet, and Hertel (2009) created programs to detect never halting
 machines, and manually proved that some machines, undetected by their programs,
 never halt. In each case, about a hundred holdouts were resisting computer and
 manual analyses. See Skelet's study in

\verb+http://skelet.ludost.net/bb/index.html+

The number of holdouts is gradually shrinking, due to the work of many people.
See the 42 holdouts of Skelet in

\verb+http://skelet.ludost.net/bb/nreg.html+

and the study of 14 of them in

\verb+http://googology.wikia.com/wiki/Forum:Sigma_project+

\item
  Daniel Briggs did some work on this question: see

\verb+https://web.archive.org/web/20121026023118/http://web.mit.edu/~dbriggs/www+

See his analysis of some machines in

\verb+https://github.com/danbriggs/Turing/blob/master/paper/HNRs.pdf+

 He wrote, in October 2021, that only 10 machines are truly difficult.

 \item
 Norbert B\'atfai, allowing transitions where
 the head can stand still, found, in August 2009, a machine M with $s(M) = 70,740,810$
 and $\sigma(M) =  4098$. Note that this machine does not follow the current rules of
 the busy beaver competition. See B\'atfai's study in

\verb+http://arxiv.org/abs/0908.4013+

\end{itemize}

\begin{center}
\begin{tabular}{|c|c|c|c|}
\hline
1964 & Green & &  $\sigma$ = 17\\
\hline
1972 & Lynn & $s$ = 435 & $\sigma$ = 22\\
\hline
1973 & Weimann & $s$ = 556 & $\sigma$ = 40\\
\hline
1974 & Lynn & $s$ = 7,707 & $\sigma$ = 112\\
\hline
August 1982 & Schult & $s$ = 134,467 & $\sigma$ = 501\\
\hline
December 1984 & Uhing & $s$ = 2,133,492 & $\sigma$ = 1,915\\
\hline
February 1986 & Uhing & $s$ = 2,358,064 & \\
\hline
February 1990 & Marxen, Buntrock & $s$ = 47,176,870 & $\sigma$ = 4,098\\
\hline
\end{tabular}
\end{center}

\bigskip

 The record holder and some other good machines:
 
\smallskip

\begin{small}
\begin{tabular}{cccccccccccc}
A0  & A1  & B0  & B1  & C0  & C1  & D0  & D1  & E0  & E1  &  $s(M)$  &$\sigma(M)$\\
1RB & 1LC & 1RC & 1RB & 1RD & 0LE & 1LA & 1LD & 1RH & 0LA  & 47,176,870 & 4098\\  
1RB & 0LD & 1LC & 1RD & 1LA & 1LC & 1RH & 1RE & 1RA & 0RB  & 23,554,764 & 4097\\
1RB & 1RA & 1LC & 1LB & 1RA & 0LD & 0RB & 1LE & 1RH & 0RB  & 11,821,234 & 4097\\
1RB & 1RA & 1LC & 1LB & 1RA & 0LD & 1RC & 1LE & 1RH & 0RB  & 11,821,220 & 4097\\
1RB & 1RA & 0LC & 0RC & 1RH & 1RD & 1LE & 0LA & 1LA & 1LE  & 11,821,190 & 4096\\
1RB & 1RA & 1LC & 0RD & 1LA & 1LC & 1RH & 1RE & 1LC & 0LA  & 11,815,076 & 4096\\
1RB & 1RA & 1LC & 1LB & 1RA & 0LD & 0RB & 1LE & 1RH & 1LC  & 11,811,040 & 4097\\
1RB & 1RA & 1LC & 1LB & 0RC & 1LD & 1RA & 0LE & 1RH & 1LC  & 11,811,040 & 4097\\
1RB & 1RA & 1LC & 1LB & 1RA & 0LD & 1RC & 1LE & 1RH & 1LC  & 11,811,026 & 4097\\
1RB & 1RA & 0LC & 0RC & 1RH & 1RD & 1LE & 1RB & 1LA & 1LE  & 11,811,010 & 4096\\
1RB & 1RA & 1LC & 1LB & 1RA & 1LD & 0RE & 0LE & 1RH & 1LC  & 11,804,940 & 4097\\
1RB & 1RA & 1LC & 1LB & 1RA & 1LD & 1RA & 0LE & 1RH & 1LC  & 11,804,926 & 4097\\
1RB & 1RA & 1LC & 0RD & 1LA & 1LC & 1RH & 1RE & 0LE & 1RB  & 11,804,910 & 4096\\  
1RB & 1RA & 1LC & 0RD & 1LA & 1LC & 1RH & 1RE & 1LC & 1RB  & 11,804,896 & 4096\\  
1RB & 1RA & 1LC & 1LB & 1RA & 1LD & 1RA & 1LE & 1RH & 0LC  & 11,798,826 & 4098\\  
1RB & 1RA & 1LC & 1RD & 1LA & 1LC & 1RH & 0RE & 1LC & 1RB  & 11,798,796 & 4097\\  
1RB & 1RA & 1LC & 1RD & 1LA & 1LC & 1RH & 1RE & 0LE & 0RB  & 11,792,724 & 4097\\
1RB & 1RA & 1LC & 1RD & 1LA & 1LC & 1RH & 1RE & 1LA & 0RB  & 11,792,696 & 4097\\
1RB & 1RA & 1LC & 1RD & 1LA & 1LC & 1RH & 1RE & 1RA & 0RB  & 11,792,682 & 4097\\
0RB & 0LC & 1RC & 1RD & 1LA & 0LE & 1RE & 1RH & 1LA & 1RA  &  2,358,065 & 1471\\
1RB & 1RH & 1LC & 1RC & 0RE & 0LD & 1LC & 0LB & 1RD & 1RA  &  2,358,064 & 1471\\  
1RB & 1LC & 0LA & 0LD & 1LA & 1RH & 1LB & 1RE & 0RD & 0RB  &  2,133,492 & 1915\\  
1RB & 0LC & 1RC & 1RD & 1LA & 0RB & 0RE & 1RH & 1LC & 1RA  &    134,467 &  501 
\end{tabular}
\end{small}

\smallskip

(All these machines can be found in Buro (1990), pp.\ 69-70.
The machines $M$ with $\sigma(M) > 1471$ were discovered by Marxen and
Buntrock. The machine with the transition $(A,0)\to (0,R,B)$ was discovered
by Buro, the next two ones were by Uhing, and the last one was by Schult.
 Heiner Marxen says there are no other $\sigma$ values within the $\sigma$
 range above).

\subsection{Turing machines with 6 states and 2 symbols}

\begin{itemize}

 \item
 Green (1964) found a machine $M$ with $\sigma(M) = 35$.

 \item
 Lynn (1972) found a machine $M$ with
 $s(M) = 522$ and $\sigma(M) = 42$.

 \item
 Brady (1983) found machines $M$ and $N$ with $s(M) = 13,488$
 and $\sigma(N) = 117$.

 \item
   Uwe Schult, cited by Ludewig et al.\ (1983) and by Dewdney (1984a),
   found, in December 1982, a machine $M$  with $s(M) = 4,208,824$
   and $\sigma(M) = 2,075$.

 \item
 Heiner Marxen and J\"{u}rgen Buntrock found, in January 1990,
 a machine $M$ with $s(M) = 13,122,572,797$ and $\sigma(M) = 136,612$.
 This machine was cited by Marxen and Buntrock (1990).
 See  study by H. Marxen in

 \verb+http://turbotm.de/~heiner/BB/simmbL6_1.html+

 \item
 Heiner Marxen and J\"{u}rgen Buntrock found, in January 1990,
 a machine $M$ with $s(M) = 8,690,333,381,690,951$
 and $\sigma(M) = 95,524,079$.
 This machine was posted on the web (Google groups) on September 3, 1997.
 See machine 2 in Marxen's bb-list in

\verb+http://turbotm.de/~heiner/BB/bb-list+

 See study by H. Marxen in

\verb+http://turbotm.de/~heiner/BB/simmbL6_2.html+

 See analysis by R. Munafo in his website:

 \verb+http://mrob.com/pub/math/ln-notes1-4.html#mb-bb-1+

and in Section \ref{sec:tm62d}.

 \item
 Heiner Marxen and J\"{u}rgen Buntrock found, in July 2000,
 a machine $M$ with $s(M) > 5.3 \times 10^{42}$
 and $\sigma(M) > 2.5 \times 10^{21}$.
 This machine was posted on the web (Google groups) on August 5, 2000.
 See machine 3 in Marxen's bb-list:

\verb+http://turbotm.de/~heiner/BB/bb-list+

and machine k in Marxen's bb-6list:

\verb+http://turbotm.de/~heiner/BB/bb-6list+

 See  study by H. Marxen in:

\verb+http://turbotm.de/~heiner/BB/simmbL6_3.html+

 \item
 Heiner Marxen and J\"{u}rgen Buntrock found, in August 2000,
 a machine $M$ with $s(M) > 6.1 \times 10^{119}$
 and $\sigma(M) > 1.4 \times 10^{60}$.
 This machine was posted on the web (Google groups) on October 23, 2000.
 See machine o in Marxen's bb-6list in:

\verb+http://turbotm.de/~heiner/BB/bb-6list+

 See study by H. Marxen in:

\verb+http://turbotm.de/~heiner/BB/simmbL6_o.html+

 See analysis by P. Michel in Section \ref{sec:tm62c}.

 \item
 Heiner Marxen and J\"{u}rgen Buntrock found, in August 2000,
 a machine $M$ with $s(M) > 6.1 \times 10^{925}$
 and $\sigma(M) > 6.4 \times 10^{462}$.
 This machine was posted on the web (Google groups) on October 23, 2000.
 See machine q in Marxen's bb-6list in

\verb+http://turbotm.de/~heiner/BB/bb-6list+
 
 See study by Marxen in

\verb+http://turbotm.de/~heiner/BB/simmbL6_q.html+

 See analyses by R. Munafo, the short one in

\verb+http://mrob.com/pub/math/ln-notes1-5.html#mb6q+

or the long one in

\verb+http://mrob.com/pub/math/ln-mb6q.html+

and see analysis by P. Michel in Section \ref{sec:tm62b}.

 \item
 Heiner Marxen and J\"{u}rgen Buntrock found, in February 2001,
 a machine $M$ with $s(M) > 3.0 \times 10^{1730}$
 and $\sigma(M) > 1.2 \times 10^{865}$.
 This machine was posted on the web (Google groups) on March 5, 2001.
 See machine r in Marxen's bb-6list in

\verb+http://turbotm.de/~heiner/BB/bb-6list+

 See study by Marxen in

\verb+http://turbotm.de/~heiner/BB/simmbL6_r.html+

 See analysis by P. Michel in Section \ref{sec:tm62a}.

 \item
 Marxen gives a list of machines $M$ with high values
of $s(M)$ and $\sigma(M)$ in

\verb+http://turbotm.de/~heiner/BB/bb-6list+
 
 \item
 Terry and Shawn Ligocki found, in November 2007,
 a machine $M$ with $s(M) > 8.9 \times 10^{1762}$
 and $\sigma(M) > 2.5 \times 10^{881}$.
 See study by H. Marxen in

\verb+http://turbotm.de/~heiner/BB/simLig62_a.html+

 See analysis by P. Michel in Section \ref{sec:tm62f}.

 \item
 Terry and Shawn Ligocki found, in December 2007,
 a machine $M$ with $s(M) > 2.5 \times 10^{2879}$
 and $\sigma(M) > 4.6 \times 10^{1439}$.
 See study by H. Marxen in

\verb+http://turbotm.de/~heiner/BB/simLig62_b.html+

 See analysis by P. Michel in Section \ref{sec:tm62e}.

 \item
 Pavel Kropitz found, in May 2010,
 a machine $M$ with $s(M) > 3.8 \times 10^{21132}$
 and $\sigma(M) > 3.1 \times 10^{10566}$.
 See study by H. Marxen in

\verb+http://turbotm.de/~heiner/BB/simKro62_a.html+

 See analysis in Section \ref{sec:tm62g}.

 \item
 Pavel Kropitz found, in June 2010,
 a machine $M$ with $s(M) > 7.4 \times 10^{36534}$
 and $\sigma(M) > 3.5 \times 10^{18267}$.
 See study by H. Marxen in

\verb+http://turbotm.de/~heiner/BB/simKro62_b.html+

 See analysis by P. Michel in Section \ref{sec:tm62h}.

 \item
 Shawn Ligocki found, in May 2022,
 a machine $M$ with $s(M) > 9.6 \times 10^{78913}$
 and $\sigma(M) > 6.0 \times 10^{39456}$.
 See study by S. Ligocki in

\verb+https://www.sligocki.com/2022/05/13/bb-6-2-e78913.html+

  \item
 Pavel Kropitz found, in May 2022,
 a machine $M$ with $s(M) > 5.4 \times 10^{197282}$
 and $\sigma(M) > 2.0 \times 10^{98641}$.

  \item
 Pavel Kropitz found, in May 2022,
 a machine $M$ with $s(M) > 8.2 \times 10^{1,292,913,985}$
 and $\sigma(M) > 1.7 \times 10^{646,456,993}$.

  \item
 Shawn Ligocki found, in May 2022,
 a machine $M$ with $s(M) > \sigma(M) > 10^{10^{10^{10^{20823}}}}$.

  \item
 Pavel Kropitz found, in May 2022,
 a machine $M$ with $s(M) > \sigma(M) > 10^{\wedge \wedge}15$,
 that is a tower of powers of 10 with 15 occurences of 10.
 It is the current record holder.

 See analysis in Section \ref{sec:tm62i}.
 
  \item
  See more analyses of the machines found in May 2022 in
  
\verb+https://groups.google.com/g/busy-beaver-discuss+
 
\end{itemize}

\begin{center}
\begin{tabular}{|c|c|c|c|}
\hline
1964 & Green & & $\sigma$ = 35\\
\hline
1972 & Lynn & $s$ = 522 & $\sigma$ = 42\\
\hline
1983 & Brady & $s$ = 13,488 & $\sigma$ = 117\\
\hline
December 1982 & Schult & $s$ = 4,208,824 & $\sigma$ = 2,075\\
\hline
February 1990 & Marxen, Buntrock & $s$ = 13,122,572,797 & $\sigma$ =
 136,612\\
\hline
September 1997 & Marxen, Buntrock & $s$ = 8,690,333,381,690,951 &
 $\sigma$ = 95,524,079\\
\hline
August 2000 & Marxen, Buntrock & $s > 5.3 \times 10^{42}$ & $\sigma >
 2.5 \times 10^{21}$\\
\hline
October 2000 & Marxen, Buntrock & $s > 6.1 \times 10^{925}$ & $\sigma >
 6.4 \times 10^{462}$\\
\hline
March 2001 & Marxen, Buntrock & $s > 3.0 \times 10^{1730}$ & $\sigma >
 1.2 \times 10^{865}$\\
\hline
November 2007 & T. and S. Ligocki & $s > 8.9 \times 10^{1762}$ &
 $\sigma > 2.5 \times 10^{881}$\\
\hline
December 2007 & T. and S. Ligocki & $s > 2.5 \times 10^{2879}$ &
 $\sigma > 4.6 \times 10^{1439}$\\
\hline
May 2010 & Kropitz & $s > 3.8 \times 10^{21132}$ & 
 $\sigma > 3.1 \times 10^{10566}$\\
\hline
June 2010 & Kropitz & $s > 7.4 \times 10^{36534}$ & 
 $\sigma > 3.5 \times 10^{18267}$\\
\hline
May 2022 & S.\ Ligocki & $s > 9.6 \times 10^{78913}$ & 
 $\sigma > 6.0 \times 10^{39456}$\\
\hline
May 2022 & Kropitz & $s > 5.4 \times 10^{197282}$ & 
 $\sigma > 2.0 \times 10^{98641}$\\
 & & $s > 8.2 \times 10^{1,292,913,985}$ & 
 $\sigma > 1.7 \times 10^{646,456,993}$\\
\hline
May 2022 & S.\ Ligocki & $s > 10^{10^{10^{10^{20823}}}}$ & 
 $\sigma > 10^{10^{10^{10^{20823}}}}$\\
\hline
May 2022 & Kropitz & $s > 10^{\wedge \wedge}15$ & 
 $\sigma > 10^{\wedge \wedge}15$\\
\hline
\end{tabular}
\end{center}

\bigskip

 The record holder and some other good machines:

\smallskip

\begin{tiny}
\begin{tabular}{cccccccccccccc}
A0  & A1  & B0  & B1  & C0  & C1  & D0  & D1  & E0  & E1  & F0  & F1
&$s(M) >$& $\sigma(M)>$ \\
1RB & 0LD & 1RC & 0RF & 1LC & 1LA & 0LE & 1RH & 1LF & 0RB & 0RC & 0RE
& $10^{\wedge \wedge}15$ & $10^{\wedge \wedge}15$\\
1RB & 0LA & 1LC & 1LF & 0LD & 0LC & 0LE & 0LB & 1RE & 0RA & 1RH & 1LD
& $10^{10^{10^{10^{20823}}}}$ & $10^{10^{10^{10^{20823}}}}$\\
1RB & 1RH & 0LC & 0LD & 1LD & 1LC & 1RE & 1LB & 1RF & 1RD & 0LD & 0RA
& $8.2 \times 10^{1,292,913,985}$ & $1.7 \times 10^{646,456,993}$\\
1RB & 1RH & 1RC & 1RA & 1RD & 0RB & 1LE & 0RC & 0LF & 0LD & 0LB & 1LA
& $5.4 \times 10^{197282}$ & $2.0 \times 10^{98641}$\\
1RB & 1RC & 1LC & 0RF & 1RA & 0LD & 0LC & 0LE & 1LD & 0RA & 1RE & 1RH
 & $9.6 \times 10^{78913}$ & $6.0 \times 10^{39456}$\\
1RB & 1LE & 1RC & 1RF & 1LD & 0RB & 1RE & 0LC & 1LA & 0RD & 1RH & 1RC
 & $7.4 \times 10^{36534}$ & $3.5 \times 10^{18267}$\\
1RB & 0LD & 1RC & 0RF & 1LC & 1LA & 0LE & 1RH & 1LA & 0RB & 0RC & 0RE
 & $3.8 \times 10^{21132}$ & $3.1 \times 10^{10566}$\\
1RB & 0LE & 1LC & 0RA & 1LD & 0RC & 1LE & 0LF & 1LA & 1LC & 1LE & 1RH
 & $2.5 \times 10^{2879}$ & $4.6 \times 10^{1439}$\\ 
1RB & 0RF & 0LB & 1LC & 1LD & 0RC & 1LE & 1RH & 1LF & 0LD & 1RA & 0LE
 & $8.9 \times 10^{1762}$ & $2.5 \times 10^{881}$\\
1RB & 0LF & 0RC & 0RD & 1LD & 1RE & 0LE & 0LD & 0RA & 1RC & 1LA & 1RH
 & $3.0 \times 10^{1730}$ & $1.2 \times 10^{865}$\\ 
1RB & 0LB & 0RC & 1LB & 1RD & 0LA & 1LE & 1LF & 1LA & 0LD & 1RH & 1LE
 & $6.1 \times 10^{925}$ & $6.4 \times 10^{462}$\\ 
1RB & 0LC & 1LA & 1RC & 1RA & 0LD & 1LE & 1LC & 1RF & 1RH & 1RA & 1RE
 & $6.1 \times 10^{119}$ & $1.4 \times 10^{60}$\\ 
1RB & 0LB & 1LC & 0RE & 1RE & 0LD & 1LA & 1LA & 0RA & 0RF & 1RE & 1RH
 & $5.5 \times 10^{99}$ & $6.9 \times 10^{49}$\\ 
1RB & 0LC & 1LA & 1LD & 1RD & 0RC & 0LB & 0RE & 1RC & 1LF & 1LE & 1RH
 & $3.2 \times 10^{98}$ & $1.1 \times 10^{49}$\\ 
1RB & 0LC & 1LA & 1RD & 1RA & 0LE & 1RA & 0RB & 1LF & 1LC & 1RD & 1RH
 & $2.0 \times 10^{95}$ & $6.7 \times 10^{47}$\\ 
1RB & 0LC & 1LA & 1RD & 0LB & 0LE & 1RA & 0RB & 1LF & 1LC & 1RD & 1RH
 & $2.0 \times 10^{95}$ & $6.7 \times 10^{47}$\\ 
1RB & 0RC & 0LA & 0RD & 1RD & 1RH & 1LE & 0LD & 1RF & 1LB & 1RA & 1RE
 & $5.3 \times 10^{42}$ & $2.5 \times 10^{21}$
\end{tabular}
\end{tiny}

\subsection{Turing machines with 7 states and 2 symbols}

\begin{itemize}

 \item Green (1964) found a machine $M$ with $\sigma(M) = 22,961$.

 \item This machine was superseded by the machine with 6 states and 2 symbols
   found in January 1990 by Heiner Marxen and J\"{u}rgen Buntrock.

 \item ``Wythagoras'' found, in March 2014, a machine $M$ with
   $s(M) > \sigma(M) > 10^{10^{10^{10^{18,705,352}}}}$.
   This machine comes from the (6,2)-TM found by Pavel Kropitz in June 2010, as follows:
   A seventh state G is added, with the transition (G,0) $\to$ (1,L,E).
   This state G becomes the initial state. Then the machine is normalized
   by swapping Left and Right and by the circular permutation of states (A C F E B D G). 
   See

   \verb+http://googology.wikia.com/wiki/User_blog:Wythagoras/A_good_bound_for_S(7)%3F+

 \item This machine was superseded by the machine with 6 states and 2 symbols found in May 2022 by Pavel Kropitz.

\end{itemize}
   
\begin{center}
\begin{tabular}{|c|c|c|c|}
\hline
1964 & Green & \multicolumn{2}{c|}{$\sigma$ = 22,961}\\
\hline
1990 & Marxen, Buntrock & \multicolumn{2}{c|}{superseded by a (6,2)-TM}\\ 
\hline
March 2014 & ``Wythagoras'' & \multicolumn{2}{c|}{$s > \sigma > 10^{10^{10^{10^{18,705,352}}}}$}\\
\hline
May 2022 & Kropitz & \multicolumn{2}{c|}{superseded by a (6,2)-TM}\\ 
\hline
\end{tabular}
\end{center}

The ``Wythagoras'' machine:

\smallskip

\begin{tabular}{cccccccccccccccc}
A0  & A1  & B0  & B1  & C0  & C1  & D0  & D1  & E0  & E1  & F0  & F1  & G0  & G1\\
1RB &     & 1RC & 0LG & 1LD & 1RB & 1LF & 1LE & 1RH & 1LF & 1RG & 0LD & 1LB & 0RF
\end{tabular}

\bigskip

\subsection{Turing machines with 2 states and 3 symbols}

\begin{itemize}

 \item
 Brady (1988) found a machine $M$ with $s(M)$ = 38
 See study by H. Marxen in

\verb+http://turbotm.de/~heiner/BB/simTM23_cb.html+

 \item
 This machine was found independently by Michel (2004),
 who gave $\sigma(M) = 9$ and conjectured that 
 $S(2,3) = 38$ and $\mathnormal{\Sigma}(2,3) = 9$.

 \item
 Lafitte and Papazian (2007) proved this conjecture.
 T.\ and S.\ Ligocki (unpublished) proved this conjecture, independently.

\end{itemize}

\begin{center}
\begin{tabular}{|c|c|c|c|}
\hline
1988 & Brady & $s$ = 38 & $\sigma$ = 9\\
\hline
2007 & Lafitte, Papazian & $S(2,3) = 38$ & $\mathnormal{\Sigma}(2,3) = 9$\\ 
\hline
\end{tabular}
\end{center}

\bigskip

The winner and some other good machines:

\smallskip

\begin{tabular}{cccccccc}
A0  & A1  & A2  & B0  & B1  & B2  &$s(M)$&$\sigma(M)$\\ 
1RB & 2LB & 1RH & 2LA & 2RB & 1LB &  38  &     9    \\
1RB & 0LB & 1RH & 2LA & 1RB & 1RA &  29  &     8    \\
0RB & 2LB & 1RH & 1LA & 1RB & 1RA &  27  &     6    \\ 
1RB & 2LA & 1RH & 1LB & 1LA & 0RA &  26  &     6    \\ 
1RB & 1LA & 1LB & 0LA & 2RA & 1RH &  26  &     6    \\ 
1RB & 1LB & 1RH & 2LA & 2RB & 1LB &  24  &     7    
\end{tabular}

\subsection{Turing machines with 3 states and 3 symbols}

\begin{itemize}
  \item
  Chester Y.\ Lee (1963) gave two machines: a machine $M$ found by
  David Jefferson, with $s(M) = 44$ and $\sigma(M) = 12$,
  and a machine $N$ found by R.\ Blodgett, with $s(M) = 57$ and
  $\sigma(M) = 9$. Note that the definition of $\mathnormal{\Sigma}(3,3)$
  is different from the usual definition (number of 1s instead of
  number of non-blank symbols). These machines are cited by
  Korthage (1966) (p.\ 114).

 \item
 Michel (2004) found machines $M$ and $N$ with
 $s(M) = 40,737$ and $\sigma(N) = 208$.

 \item
 Brady found, in November 2004, a machine $M$ with
 $s(M) = 29,403,894$ and $\sigma(M) = 5600$.
 See  \verb+http://www.cse.unr.edu/~al/BusyBeaver.html+
 
 See study by H. Marxen in

\verb+http://turbotm.de/~heiner/BB/simAB3Y_b.html+

 \item
 Brady found, in December 2004, a machine $M$ with
 $s(M) = 92,649,163$ and $\sigma(M) = 13,949$.
 See  \verb+http://www.cse.unr.edu/~al/BusyBeaver.html+
 
 See study by H. Marxen in

\verb+http://turbotm.de/~heiner/BB/simAB3Y_c.html+

 See analysis by P. Michel in Section \ref{sec:tm33a}.

 \item
 Myron P. Souris found, in July 2005 (M.P. Souris said: actually in
 1995, but then  no one seemed to care), machines $M$ and $N$ with
 $s(M) = 544,884,219$ and $\sigma(N) = 36,089$.
 See study of $M$ by H. Marxen in

\verb+http://turbotm.de/~heiner/BB/simMS33_b.html+

and study of $N$ by H. Marxen in

\verb+http://turbotm.de/~heiner/BB/simMS33_a.html+

 See analysis of $M$ by P. Michel in Section \ref{sec:tm33c},

and analysis of $N$ by P. Michel in Section \ref{sec:tm33b}.

 \item
 Gr\'egory Lafitte and Christophe Papazian found, in August 2005,
 a machine $M$ with $s(M) = 4,939,345,068$ and $\sigma(M) = 107,900$.
 Eee  study by H. Marxen in

\verb+http://turbotm.de/~heiner/BB/simLaf33_b.html+

 See analysis by P. Michel in Section \ref{sec:tm33d}.

 \item
 Gr\'egory Lafitte and Christophe Papazian found, in September 2005,
 a machine $M$ with $s(M) = 987,522,842,126$ and $\sigma(M) = 1,525,688$.
 See  study by H. Marxen in

\verb+http://turbotm.de/~heiner/BB/simLaf33_d.html+

 See analysis by P. Michel in Section \ref{sec:tm33e}.

 \item
 Gr\'egory Lafitte and Christophe Papazian found, in April 2006,
 a machine $M$ with $s(M) = 4,144,465,135,614$ and $\sigma(M) = 2,950,149$.
 See  study by H. Marxen in

\verb+http://turbotm.de/~heiner/BB/simLaf33_e.html+

 See analysis by P. Michel in Section \ref{sec:tm33f}.

 \item
 Terry and Shawn Ligocki found, in August 2006,
 a machine $M$ with $s(M) = 4,345,166,620,336,565$
 and $\sigma(M) = 95,524,079$.
 See  study by H. Marxen in

\verb+http://turbotm.de/~heiner/BB/simLig33_a.html+

 See analysis in Section \ref{sec:tm33g}.

 \item
 Terry and Shawn Ligocki found, in November 2007,
 a machine $M$ with $s(M) = 119,112,334,170,342,540$
 and $\sigma(M) = 374,676,383$.
 See  study by H. Marxen in

\verb+http://turbotm.de/~heiner/BB/simLig33_b.html+

 See analysis by P. Michel in Section \ref{sec:tm33h}.

 It is the current record holder.

 \item
 Brady gives a list of machines with high values of $s(M)$ in

\verb+http://www.cse.unr.edu/~al/BusyBeaver.html+

\end{itemize}

\begin{center}
\begin{tabular}{|c|c|c|c|}
\hline
1966 & cited by Korfhage & $s$ = 57 & $\sigma'$ = 12\\
  \hline
October 2004 & Michel & $s$ = 40,737 & $\sigma$ = 208\\
\hline
November 2004 & Brady & $s$ = 29,403,894 & $\sigma$ = 5,600\\
\hline
December 2004 & Brady &  $s$ = 92,649,163 & $\sigma$ = 13,949\\
\hline
July 2005 & Souris & $s$ = 544,884,219 & $\sigma$ = 36,089\\
\hline
August 2005 & Lafitte, Papazian & $s$ = 4,939,345,068 & $\sigma$ = 107,900\\
\hline
September 2005 & Lafitte, Papazian & $s$ = 987,522,842,126 & $\sigma$ =
 1,525,688\\
\hline
April 2006 & Lafitte, Papazian & $s$ = 4,144,465,135,614 & $\sigma$ =
 2,950,149\\
\hline
August 2006 & T. and S. Ligocki & $s$ = 4,345,166,620,336,565 & $\sigma$
 = 95,524,079\\
\hline
November 2007 & T. and S. Ligocki & $s$ = 119,112,334,170,342,540 &
 $\sigma$ = 374,676,383\\
\hline
\end{tabular}
\end{center}

\bigskip

The record holder and some other good machines:

\smallskip

\begin{scriptsize}
\begin{tabular}{ccccccccccc}
A0  & A1  & A2  & B0  & B1  & B2  & C0  & C1  & C2  
 & $s(M)$ & $\sigma(M)$\\
1RB & 2LA & 1LC & 0LA & 2RB & 1LB & 1RH & 1RA & 1RC
 & 119,112,334,170,342,540 & 374,676,383\\ 
1RB & 2RC & 1LA & 2LA & 1RB & 1RH & 2RB & 2RA & 1LC
 & 4,345,166,620,336,565 & 95,524,079\\
1RB & 1RH & 2LC & 1LC & 2RB & 1LB & 1LA & 2RC & 2LA
 & 4,144,465,135,614 & 2,950,149\\
1RB & 2LA & 1RA & 1RC & 2RB & 0RC & 1LA & 1RH & 1LA
 & 987,522,842,126 & 1,525,688\\
1RB & 1RH & 2RB & 1LC & 0LB & 1RA & 1RA & 2LC & 1RC
 & 4,939,345,068 & 107,900\\
1RB & 2LA & 1RA & 1LB & 1LA & 2RC & 1RH & 1LC & 2RB
 & 1,808,669,066 & 43,925\\
1RB & 2LA & 1RA & 1LC & 1LA & 2RC & 1RH & 1LA & 2RB
 & 1,808,669,046 & 43,925\\
1RB & 1LB & 2LA & 1LA & 1RC & 1RH & 0LA & 2RC & 1LC
 & 544,884,219 & 32,213\\
1RB & 0LA & 1LA & 2RC & 1RC & 1RH & 2LC & 1RA & 0RC
 & 408,114,977 & 20,240\\
1RB & 2RA & 2RC & 1LC & 1RH & 1LA & 1RA & 2LB & 1LC
 & 310,341,163 & 36,089\\
1RB & 1RH & 2LC & 1LC & 2RB & 1LB & 1LA & 0RB & 2LA
 & 92,649,163 & 13,949\\
1RB & 2LA & 1LA & 2LA & 1RC & 2RB & 1RH & 0LC & 0RA
 & 51,525,774 & 7,205\\
1RB & 2RA & 1LA & 2LA & 2LB & 2RC & 1RH & 2RB & 1RB
 & 47,287,015 & 12,290\\
1RB & 2RA & 1LA & 2LC & 0RC & 1RB & 1RH & 2LA & 1RB
 & 29,403,894 & 5,600 
\end{tabular}
\end{scriptsize}

\smallskip

 (The first two machines were discovered by Terry and Shawn Ligocki,
 the next five ones were by Lafitte and Papazian,
 the next three ones were by Souris, and the last four ones were by Brady).

\subsection{Turing machines with 4 states and 3 symbols}

\begin{itemize}

 \item
 Terry and Shawn Ligocki found, in April 2005,
 a machine $M$ with $s(M) = 250,096,776$ and $\sigma(M) = 15,008$.
 See  study by H. Marxen in

\verb+http://turbotm.de/~heiner/BB/simLig43_a.html+

 \item
 This machine was superseded by the machines with 3 states
 and 3 symbols found in July 2005 by Myron P. Souris.

 \item
 Terry and Shawn Ligocki found, in October 2007,
 a machine $M$ with $s(M) > 1.5 \times  10^{1426}$ and
 $\sigma(M) > 1.1 \times 10^{713}$.
See  study by H. Marxen in

\verb+http://turbotm.de/~heiner/BB/simLig43_b.html+

 \item
 Terry and Shawn Ligocki found successively, in November 2007,
 machines $M$ with 

\begin{itemize}
   \item
    $s(M) > 7.7 \times 10^{1618}$ and
    $\sigma(M) > 1.6 \times 10^{809}$.
    See study by H. Marxen in

\verb+http://turbotm.de/~heiner/BB/simLig43_c.html+
 
   \item
    $s(M) > 3.7 \times 10^{1973}$ and
    $\sigma(M) > 8.0 \times 10^{986}$.
    See study by H. Marxen in

    \verb+http://turbotm.de/~heiner/BB/simLig43_d.html+
    
   \item
    $s(M) > 3.9 \times 10^{7721}$ and
    $\sigma(M) > 4.0 \times 10^{3860}$.
    See study by H. Marxen in

\verb+http://turbotm.de/~heiner/BB/simLig43_e.html+
  
   \item
    $s(M) > 3.9 \times 10^{9122}$ and
    $\sigma(M) > 2.5 \times 10^{4561}$.
    See study by H. Marxen in

\verb+http://turbotm.de/~heiner/BB/simLig43_f.html+

\end{itemize}

 \item
 Terry and Shawn Ligocki found successively, in December 2007,
 machines $M$ with
 
\begin{itemize}
   \item
    $s(M) > 7.9 \times 10^{9863}$ and
    $\sigma(M) > 8.9 \times 10^{4931}$.
    See  study by H. Marxen in

\verb+http://turbotm.de/~heiner/BB/simLig43_g.html+
 
   \item
    $s(M) > 5.3 \times 10^{12068}$ and
    $\sigma(M) > 4.2 \times 10^{6034}$.
    See study by H. Marxen in

\verb+http://turbotm.de/~heiner/BB/simLig43_h.html+

\end{itemize}

 \item
 Terry and Shawn Ligocki found, in January 2008,
 a machine $M$ with $s(M) > 1.0 \times 10^{14072}$ and
 $\sigma(M) > 1.3 \times 10^{7036}$.
   See study by H. Marxen in

\verb+http://turbotm.de/~heiner/BB/simLig43_i.html+

It is the current record holder.

\end{itemize}

\begin{center}
\begin{tabular}{|c|c|c|c|}
\hline
April 2005 & T. and S. Ligocki & $s$ = 250,096,776 & $\sigma$ = 15,008\\
\hline
July 2005 & Souris & \multicolumn{2}{c|}{superseded by a (3,3)-TM}\\
\hline
October 2007 & T. and S. Ligocki & $s > 1.5 \times 10^{1426}$ & $\sigma > 1.1 \times 10^{713}$\\
\hline
 & & $s > 7.7 \times 10^{1618}$ & $\sigma > 1.6 \times 10^{809}$\\
November 2007 & T. and S. Ligocki & $s > 3.7 \times 10^{1973}$ & $\sigma > 8.0 \times 10^{986}$\\
 & & $s > 3.9 \times 10^{7721}$ & $\sigma > 4.0 \times 10^{3860}$\\
 & & $s > 3.9 \times 10^{9122}$ & $\sigma > 2.5 \times 10^{4561}$\\
\hline
December 2007 & T. and S. Ligocki & $s > 7.9 \times 10^{9863}$ & $\sigma > 8.9 \times 10^{4931}$\\
 & & $s > 5.3 \times 10^{12068}$ & $\sigma > 4.2 \times 10^{6034}$\\
\hline
January 2008 & T. and S. Ligocki & $s > 1.0 \times 10^{14072}$ & $\sigma > 1.3 \times 10^{7036}$\\
\hline
\end{tabular}
\end{center}

\bigskip

The record holder and the past record holders:

\smallskip

\begin{tiny}
\begin{tabular}{cccccccccccccc}
A0  & A1  & A2  & B0  & B1  & B2  & C0  & C1  & C2  & D0  & D1  & D2  &
 $s(M)$ & $\sigma(M)$\\
1RB & 1RH & 2RC & 2LC & 2RD & 0LC & 1RA & 2RB & 0LB & 1LB & 0LD & 2RC &
 $> 1.0 \times 10^{14072}$ & $> 1.3 \times 10^{7036}$\\
1RB & 0LB & 1RD & 2RC & 2LA & 0LA & 1LB & 0LA & 0LA & 1RA & 0RA & 1RH &
 $> 5.3 \times 10^{12068}$ & $> 4.2 \times 10^{6034}$\\
1RB & 1LD & 1RH & 1RC & 2LB & 2LD & 1LC & 2RA & 0RD & 1RC & 1LA & 0LA &
 $> 7.9 \times 10^{9863}$ & $> 8.9 \times 10^{4931}$\\
1RB & 2LD & 1RH & 2LC & 2RC & 2RB & 1LD & 0RC & 1RC & 2LA & 2LD & 0LB &
 $> 3.9 \times 10^{9122}$ & $> 2.5 \times 10^{4561}$\\
1RB & 1LA & 1RD & 2LC & 0RA & 1LB & 2LA & 0LB & 0RD & 2RC & 1RH & 0LC &
 $> 3.9 \times 10^{7721}$ & $> 4.0 \times 10^{3860}$\\
1RB & 1RA & 0LB & 2LC & 1LB & 1RC & 0RD & 2LC & 1RA & 2RA & 1RH & 1RC &
 $> 3.7 \times 10^{1973}$ & $> 8.0 \times 10^{986}$\\
1RB & 2RC & 1RA & 2LC & 1LA & 1LB & 2LD & 0LB & 0RC & 0RD & 1RH & 0RA &
 $> 7.7 \times 10^{1618}$ & $> 1.6 \times 10^{809}$\\
1RB & 0LC & 1RH & 2LC & 1RD & 0LB & 2LA & 1LC & 1LA & 1RB & 2LD & 2RA &
 $> 1.5 \times 10^{1426}$ & $> 1.1 \times 10^{713}$\\
0RB & 1LD & 1RH & 1LA & 1RC & 1RD & 1RB & 2LC & 1RC & 1LA & 1LC & 2RB &
 250,096,776 & 15,008
\end{tabular}
\end{tiny}

\subsection{Turing machines with 2 states and 4 symbols}

\begin{itemize}

 \item
 Brady (1988) found a machine $M$ with $s(M) = 7,195$.

 \item
 This machine was found independently and analyzed by
 Michel (2004), who gave $\sigma(M) = 90$.
 See study by H. Marxen in

\verb+http://turbotm.de/~heiner/BB/simTM24_b.html+

 See analysis by P. Michel in Section \ref{sec:tm24b}.

 \item
 Terry and Shawn Ligocki found, in February 2005,
 a machine $M$ with $s(M) = 3,932,964$ and $\sigma(M) = 2,050$.
 See study by H. Marxen in

\verb+http://turbotm.de/~heiner/BB/simLig24_a.html+

 See analysis by P. Michel in Section \ref{sec:tm24a}.

 It is the current record holder. There is no machine
 between the first two ones (Ligocki, Brady). There is
 no machine such that $3,932,964 < s(M) < 200,000,000$
 (Ligocki, September 2005).

\end{itemize}

\begin{center}
\begin{tabular}{|c|c|c|c|}
\hline
1988 & Brady & $s$ = 7,195 & $\sigma$ = 90\\
\hline
February 2005 & T. and S. Ligocki & $s$ = 3,932,964 & $\sigma$ = 2,050\\
\hline
\end{tabular}
\end{center}

\bigskip

 The record holder and some other good machines:

\smallskip

\begin{tabular}{cccccccccc}
A0  & A1  & A2  & A3  & B0  & B1  & B2  & B3   & $s(M)$ & $\sigma(M)$\\
1RB & 2LA & 1RA & 1RA & 1LB & 1LA & 3RB & 1RH  & 3,932,964& 2,050\\
1RB & 3LA & 1LA & 1RA & 2LA & 1RH & 3RA & 3RB  &     7,195&    90\\
1RB & 3LA & 1LA & 1RA & 2LA & 1RH & 3LA & 3RB  &     6,445&    84\\
1RB & 3LA & 1LA & 1RA & 2LA & 1RH & 2RA & 3RB  &     6,445&    84\\
1RB & 2RB & 3LA & 2RA & 1LA & 3RB & 1RH & 1LB  &    2,351 &    60
\end{tabular}

\subsection{Turing machines with 3 states and 4 symbols}

\begin{itemize}

 \item
 Terry and Shawn Ligocki found, in April 2005,
 a machine $M$ with $s(M) = 262,759,288$ and
 $\sigma(M) = 17,323$.
 See study by H. Marxen in

 \verb+http://turbotm.de/~heiner/BB/simLig34_a.html+
 
 \item
 This machine was superseded by the machines with 3 states
 and 3 symbols found in July 2005 by Myron P. Souris.

 \item
 Terry and Shawn Ligocki found, in September 2007,
 a machine $M$ with $s(M) > 5.7 \times 10^{52}$ and
 $\sigma(M) > 2.4 \times 10^{26}$.
 See study by H. Marxen in

\verb+http://turbotm.de/~heiner/BB/simLig34_b.html+

 \item
 Terry and Shawn Ligocki found successively, in October 2007,
 machines $M$ with
\begin{itemize}
   \item
    $s(M) > 4.3 \times 10^{281}$ and
    $\sigma(M) > 6.0 \times 10^{140}$.
    See study by H. Marxen in

\verb+http://turbotm.de/~heiner/BB/simLig34_c.html+
 
   \item
    $s(M) > 7.6 \times 10^{868}$ and
    $\sigma(M) > 4.6 \times 10^{434}$.
    See study by H. Marxen in

\verb+http://turbotm.de/~heiner/BB/simLig34_d.html+

   \item
    $s(M) > 3.1 \times 10^{1256}$ and
    $\sigma(M) > 2.1 \times 10^{628}$.
    See study by H. Marxen in

\verb+http://turbotm.de/~heiner/BB/simLig34_e.html+

\end{itemize}

 \item
 Terry and Shawn Ligocki found successively, in November 2007,
 machines $M$ with 

\begin{itemize}
   \item
    $s(M) > 8.4 \times 10^{2601}$ and
    $\sigma(M) > 1.7 \times 10^{1301}$.
    See study by H. Marxen in

\verb+http://turbotm.de/~heiner/BB/simLig34_f.html+

   \item
    $s(M) > 3.4 \times 10^{4710}$ and
    $\sigma(M) > 1.4 \times 10^{2355}$.
    See study by H. Marxen in

\verb+http://turbotm.de/~heiner/BB/simLig34_g.html+

   \item
    $s(M) > 5.9 \times 10^{4744}$ and
    $\sigma(M) > 2.2 \times 10^{2372}$.
    See study by H. Marxen in

\verb+http://turbotm.de/~heiner/BB/simLig34_h.html+

\end{itemize}

 \item
 Terry and Shawn Ligocki found, in December 2007,
 a machine $M$ with $s(M) > 5.2 \times 10^{13036}$ and
 $\sigma(M) > 3.7 \times 10^{6518}$.
 See study by H. Marxen in

\verb+http://turbotm.de/~heiner/BB/simLig34_i.html+

 It is the current record holder.

\end{itemize}

\begin{center}
\begin{tabular}{|c|c|c|c|}
\hline
April 2005 & T. and S. Ligocki & $s$ = 262,759,288 & $\sigma$ = 17,323\\
\hline
July 2005 & Souris & \multicolumn{2}{c|}{superseded by a (3,3)-TM}\\
\hline
September 2007 & T. and S. Ligocki & $s > 5.7 \times 10^{52}$ & $\sigma > 2.4 \times 10^{26}$\\
\hline
 & & $s > 4.3 \times 10^{281}$ & $\sigma > 6.0 \times 10^{140}$\\
October 2007 & T. and S. Ligocki & $s > 7.6 \times 10^{868}$ & $\sigma > 4.6 \times 10^{434}$\\
 & & $s > 3.1 \times 10^{1256}$ & $\sigma > 2.1 \times 10^{628}$\\
\hline
 & & $s > 8.4 \times 10^{2601}$ & $\sigma > 1.7 \times 10^{1301}$\\
November 2007 & T. and S. Ligocki & $s > 3.4 \times 10^{4710}$ & $\sigma > 1.4 \times 10^{2355}$\\
 & & $s > 5.9 \times 10^{4744}$ & $\sigma > 2.2 \times 10^{2372}$\\
\hline
December 2007 & T. and S. Ligocki & $s > 5.2 \times 10^{13036}$ & $\sigma > 3.7 \times 10^{6518}$\\
\hline
\end{tabular}
\end{center}

\bigskip

 The record holder and the past record holders:

\smallskip

\begin{tiny}
\begin{tabular}{cccccccccccccc}
A0  & A1  & A2  & A3  & B0  & B1  & B2  & B3  & C0  & C1  & C2  & C3 &
 $s(M)$ & $ \sigma(M)$\\
1RB & 1RA & 2LB & 3LA & 2LA & 0LB & 1LC & 1LB & 3RB & 3RC & 1RH & 1LC &
 $> 5.2 \times 10^{13036}$ & $> 3.7 \times 10^{6518}$\\
1RB & 1RA & 1LB & 1RC & 2LA & 0LB & 3LC & 1RH & 1LB & 0RC & 2RA & 2RC &
 $> 5.9 \times 10^{4744}$ & $> 2.2 \times 10^{2372}$\\
1RB & 2LB & 2RA & 1LA & 2LA & 1RC & 0LB & 2RA & 1RB & 3LC & 1LA & 1RH &
 $> 3.4 \times 10^{4710}$ & $> 1.4 \times 10^{2355}$\\
1RB & 1LA & 3LA & 3RC & 2LC & 2LB & 1RB & 1RA & 2LA & 3LC & 1RH & 1LB &
 $> 8.4 \times 10^{2601}$ & $> 1.7 \times 10^{1301}$\\
1RB & 3LA & 3RC & 1RA & 2RC & 1LA & 1RH & 2RB & 1LC & 1RB & 1LB & 2RA &
 $> 3.1 \times 10^{1256}$ & $> 2.1 \times 10^{628}$\\
1RB & 0RB & 3LC & 1RC & 0RC & 1RH & 2RC & 3RC & 1LB & 2LA & 3LA & 2RB &
 $> 7.6 \times 10^{868}$ & $> 4.6 \times 10^{434}$\\
1RB & 3RB & 2LC & 3LA & 0RC & 1RH & 2RC & 1LB & 1LB & 2LA & 3RC & 2LC &
 $> 4.3 \times 10^{281}$ & $> 6.0 \times 10^{140}$\\
1RB & 1LA & 1LB & 1RA & 0LA & 2RB & 2LC & 1RH & 3RB & 2LB & 1RC & 0RC &
 $> 5.7 \times 10^{52}$ & $> 2.4 \times 10^{26}$\\
1RB & 3LC & 0RA & 0LC & 2LC & 3RC & 0RC & 1LB & 1RA & 0LB & 0RB & 1RH &
 262,759,288 & 17,323
\end{tabular}
\end{tiny}

\subsection{Turing machines with 2 states and 5 symbols}

\begin{itemize}

 \item
 Terry and Shawn Ligocki found, in February 2005,
 machines $M$ and $N$ with $s(M) = 16,268,767$ and
 $\sigma(N) = 4,099$.
 See study by H. Marxen in

\verb+http://turbotm.de/~heiner/BB/simLig25_a.html+

 \item
 Terry and Shawn Ligocki found, in April 2005,
 a machine $M$ with $s(M) = 148,304,214$ and
 $\sigma(M) = 11,120$.
 See study by H. Marxen in

\verb+http://turbotm.de/~heiner/BB/simLig25_c.html+

 \item
 Gr\'egory Lafitte and Christophe Papazian found, in August 2005,
 a machine $M$ with $s(M) = 8,619,024,596$ and
 $\sigma(M) = 90,604$.
 See study by H. Marxen in

\verb+http://turbotm.de/~heiner/BB/simLaf25_a.html+

 \item
 Gr\'egory Lafitte and Christophe Papazian found, in September 2005,
 a machine $M$ with $\sigma(M) = 97,104$ (and $s(M) = 7,543,673,517$).
 See study by H. Marxen in

\verb+http://turbotm.de/~heiner/BB/simLaf25_c.html+

 \item
 Gr\'egory Lafitte and Christophe Papazian found, in October 2005,
 a machine $M$ with $s(M) = 233,431,192,481$ and
 $\sigma(M) = 458,357$.
 See study by H. Marxen in

\verb+http://turbotm.de/~heiner/BB/simLaf25_d.html+

 \item
 Gr\'egory Lafitte and Christophe Papazian found, in October 2005,
 a machine $M$ with $s(M) = 912,594,733,606$ and $\sigma(M) = 1,957,771$.
 See study by H. Marxen in

\verb+http://turbotm.de/~heiner/BB/simLaf25_f.html+

 See analysis by P. Michel in Section \ref{sec:tm25a}.

 \item
 Gr\'egory Lafitte and Christophe Papazian found, in December 2005,
 a machine $M$ with $s(M) = 924,180,005,181$ (and $\sigma(M) = 1,137,477$).
 See study by H. Marxen in

\verb+http://turbotm.de/~heiner/BB/simLaf25_g.html+

 See analysis by P. Michel in Section \ref{sec:tm25b}.
 
 \item
 Gr\'egory Lafitte and Christophe Papazian found, in May 2006,
 a machine $M$ with $s(M) = 3,793,261,759,791$ and $\sigma(M) = 2,576,467$.
 See study by H. Marxen in

\verb+http://turbotm.de/~heiner/BB/simLaf25_h.html+

 See analysis by P. Michel in Section \ref{sec:tm25c}.

 \item
 Gr\'egory Lafitte and Christophe Papazian found, in June 2006,
 a machine $M$ with $s(M) = 14,103,258,269,249$ and $\sigma(M) = 4,848,239$.
 See study by H. Marxen in

\verb+http://turbotm.de/~heiner/BB/simLaf25_i.html+

 See analysis by P. Michel in Section \ref{sec:tm25d}.

 \item
 Gr\'egory Lafitte and Christophe Papazian found, in July 2006,
 a machine $M$ with $s(M) = 26,375,397,569,930$ (and $\sigma(M) = 143$).
 See study by H. Marxen in

\verb+http://turbotm.de/~heiner/BB/simLaf25_j.html+

 See comments in Section \ref{sec:ncl}.

 \item
 Terry and Shawn Ligocki found, in August 2006,
 a machine $M$ with $$s(M) = 7,069,449,877,176,007,352,687$$
 and $\sigma(M) = 172,312,766,455$.
 See study by H. Marxen in

\verb+http://turbotm.de/~heiner/BB/simLig25_j.html+

 See analysis in Section \ref{sec:tm25e}.

 \item
 Terry and Shawn Ligocki found, in October 2007,
 a machine $M$ with $s(M) > 5.2 \times 10^{61}$ and
 $\sigma(M) > 9.3 \times 10^{30}$.
 See study by H. Marxen in

\verb+http://turbotm.de/~heiner/BB/simLig25_k.html+

 \item
 Terry and Shawn Ligocki found, in October 2007,
 two machines $M$ and $N$ with $s(M) = s(N) > 1.6 \times 10^{211}$ and
 $\sigma(M) = \sigma(N) > 5.2 \times 10^{105}$.
    See study by H. Marxen of $M$ in

\verb+http://turbotm.de/~heiner/BB/simLig25_l.html+

    and study by H. Marxen of $N$ in

\verb+http://turbotm.de/~heiner/BB/simLig25_m.html+

 \item
 Terry and Shawn Ligocki found, in November 2007,
 a machine $M$ with $s(M) > 1.9 \times 10^{704}$ and
 $\sigma(M) > 1.7 \times 10^{352}$.
 See study by H. Marxen in

\verb+http://turbotm.de/~heiner/BB/simLig25_n.html+

 See analysis by P. Michel in Section \ref{sec:tm25h}.

 It is the current record holder.

\end{itemize}

\begin{center}
\begin{tabular}{|c|c|c|c|}
\hline
February 2005 & T. and S. Ligocki & $s$ = 16,268,767 & $\sigma$ = 4,099\\
\hline
April 2005 & T. and S. Ligocki & $s$ = 148,304,214 & $\sigma$ = 11,120\\
\hline
August 2005 & Lafitte, Papazian & $s$ = 8,619,024,596 & $\sigma$ =
 90,604\\
\hline
September 2005 & Lafitte, Papazian & & $\sigma$ = 97,104\\
\hline
October 2005 & Lafitte, Papazian & $s$ = 233,431,192,481 & $\sigma$ = 458,357\\
 & & $s$ = 912,594,733,606 & $\sigma$ = 1,957,771\\
\hline
December 2005 & Lafitte, Papazian & $s$ = 924,180,005,181 & \\
\hline
May 2006 & Lafitte, Papazian & $s$ = 3,793,261,759,791 & $\sigma$ = 2,576,467\\
\hline
June 2006 & Lafitte, Papazian & $s$ = 14,103,258,269,249 & $\sigma$ = 4,848,239\\
\hline
July 2006 & Lafitte, Papazian & $s$ = 26,375,397,569,930 & \\
\hline
August 2006 & T. and S. Ligocki & $s > 7.0 \times 10^{21}$ & $\sigma$ = 172,312,766,455\\
\hline
October 2007 & T. and S. Ligocki & $s > 5.2 \times 10^{61}$ & $\sigma > 9.3 \times 10^{30}$\\
 & & $s > 1.6 \times 10^{211}$ & $\sigma > 5.2 \times 10^{105}$\\
\hline
November 2007 & T. and S. Ligocki & $s > 1.9 \times 10^{704}$ &  $\sigma > 1.7 \times 10^{352}$\\
\hline
\end{tabular}
\end{center}

\bigskip

\noindent Note: Two machines were discovered by T. and S. Ligocki in
February 2005 with $s(M) = 16,268,767$, and two were in October 2007
with $s(M) > 1.6 \times 10^{211}$.

\bigskip
 
 The record holder and some other good machines:

\bigskip

\begin{tiny}
\begin{tabular}{cccccccccccc}
A0  & A1  & A2  & A3  & A4  & B0  & B1  & B2  & B3  & B4  &
 $s(M)$ & $\sigma(M)$\\
1RB & 2LA & 1RA & 2LB & 2LA & 0LA & 2RB & 3RB & 4RA & 1RH &
 $> 1.9 \times 10^{704}$ & $> 1.7 \times 10^{352}$\\
1RB & 2LA & 4RA & 2LB & 2LA & 0LA & 2RB & 3RB & 4RA & 1RH &
 $> 1.6 \times 10^{211}$ & $> 5.2 \times 10^{105}$\\
1RB & 2LA & 4RA & 2LB & 2LA & 0LA & 2RB & 3RB & 1RA & 1RH &
 $> 1.6 \times 10^{211}$ & $> 5.2 \times 10^{105}$\\
1RB & 2LA & 4RA & 1LB & 2LA & 0LA & 2RB & 3RB & 2RA & 1RH &
 $> 5.2 \times 10^{61}$ & $> 9.3 \times 10^{30}$\\
1RB & 0RB & 4RA & 2LB & 2LA & 2LA & 1LB & 3RB & 4RA & 1RH &
 $> 7.0 \times 10^{21}$ & 172,312,766,455\\
1RB & 3LA & 3LB & 0LB & 1RA & 2LA & 4LB & 4LA & 1RA & 1RH &
 339,466,124,499,007,251 & 1,194,050,967\\
1RB & 3RB & 3RA & 1RH & 2LB & 2LA & 4RA & 4RB & 2LB & 0RA &
 339,466,124,499,007,214 & 1,194,050,967\\
1RB & 1RH & 4LA & 4LB & 2RA & 2LB & 2RB & 3RB & 2RA & 0RB &
 91,791,666,497,368,316 & 620,906,587\\
1RB & 3LA & 1LA & 0LB & 1RA & 2LA & 4LB & 4LA & 1RA & 1RH &
 37,716,251,406,088,468 & 398,005,342\\
1RB & 2RA & 1LA & 3LA & 2RA & 2LA & 3RB & 4LA & 1LB & 1RH &
 9,392,084,729,807,219 & 114,668,733\\
1RB & 2RA & 1LA & 1LB & 3LB & 2LA & 3RB & 1RH & 4RA & 1LA &
 417,310,842,648,366 & 36,543,045
\end{tabular}
\end{tiny}

\smallskip

  (These machines were discovered by Terry and Shawn Ligocki).

\bigskip

 Previous record holders and some other good machines: 

\bigskip

\begin{tiny}
\begin{tabular}{cccccccccccc}
A0  & A1  & A2  & A3  & A4  & B0  & B1  & B2  & B3  & B4  & 
 $s(M)$ & $\sigma(M)$\\
1RB & 3LA & 1LA & 4LA & 1RA & 2LB & 2RA & 1RH & 0RA & 0RB &
 26,375,397,569,930 & 143\\
1RB & 3LB & 4LB & 4LA & 2RA & 2LA & 1RH & 3RB & 4RA & 3RB &
 14,103,258,269,249 & 4,848,239\\
1RB & 3RA & 4LB & 2RA & 3LA & 2LA & 1RH & 4RB & 4RB & 2LB &
 3,793,261,759,791 & 2,576,467\\
1RB & 3RA & 1LA & 1LB & 3LB & 2LA & 4LB & 3RA & 2RB & 1RH &
 924,180,005,181 & 1,137,477\\
1RB & 3LB & 1RH & 1LA & 1LA & 2LA & 3RB & 4LB & 4LB & 3RA &
 912,594,733,606 & 1,957,771\\
1RB & 2RB & 3LA & 2RA & 3RA & 2LB & 2LA & 3LA & 4RB & 1RH &
 469,121,946,086 & 668,420\\
1RB & 3RB & 3RB & 1LA & 3LB & 2LA & 3RA & 4LB & 2RA & 1RH &
 233,431,192,481 &458,357\\
1RB & 3LA & 1LB & 1RA & 3RA & 2LB & 3LA & 3RA & 4RB & 1RH &
 8,619,024,596 & 90,604\\
1RB & 2RB & 3RB & 4LA & 3RA & 0LA & 4RB & 1RH & 0RB & 1LB &
 7,543,673,517 & 97,104\\
1RB & 4LA & 1LA & 1RH & 2RB & 2LB & 3LA & 1LB & 2RA & 0RB &
 7,021,292,621 & 37\\
1RB & 2RB & 3LA & 2RA & 3RA & 2LB & 2LA & 1LA & 4RB & 1RH &
 4,561,535,055 & 64,665\\
1RB & 3LA & 4LA & 1RA & 1LA & 2LA & 1RH & 4RA & 3RB & 1RA &
 148,304,214 & 11,120\\
1RB & 3LA & 4LA & 1RA & 1LA & 2LA & 1RH & 1LA & 3RB & 1RA &
 16,268,767 & 3,685\\
1RB & 3RB & 2LA & 0RB & 1RH & 2LA & 4RB & 3LB & 2RB & 3RB &
 15,754,273 & 4,099
\end{tabular}
\end{tiny}

\smallskip

  (The first eleven machines were discovered by Lafitte and
  Papazian, and the last three ones were by T. and S. Ligocki).

\subsection{Turing machines with 2 states and 6 symbols}

\begin{itemize}

 \item
 Terry and Shawn Ligocki found, in February 2005,
 machines $M$ and $N$ with $s(M) = 98,364,599$ and
 $\sigma(N) = 10,574$.
 See study by H. Marxen in

 \verb+http://turbotm.de/~heiner/BB/simLig26_a.html+

 \item
 Terry and Shawn Ligocki found, in April 2005,
 a machine $M$ with $s(M) = 493,600,387$ and
 $\sigma(M) = 15,828$.
 See  study by H. Marxen in

\verb+http://turbotm.de/~heiner/BB/simLig26_c.html+

 \item
 This machine was superseded by the machine with 2 states
 and 5 symbols found in August 2005 by Gr\'egory Lafitte
 and Christophe Papazian.

 \item
 Terry and Shawn Ligocki found, in September 2007,
 a machine $M$ with $s(M) > 2.3 \times 10^{54}$ and
 $\sigma(M) > 1.9 \times 10^{27}$.
 See study by H. Marxen in

\verb+http://turbotm.de/~heiner/BB/simLig26_d.html+

 \item
 This machine was superseded by the machine with 2 states
 and 5 symbols found in October 2007 by Terry and Shawn
 Ligocki.

 \item
 Terry and Shawn Ligocki found successively, in November 2007,
 machines $M$ with

\begin{itemize}
   \item
    $s(M) > 4.9 \times 10^{1643}$ and
    $\sigma(M) > 8.6 \times 10^{821}$.
    See study by H. Marxen in

 \verb+http://turbotm.de/~heiner/BB/simLig26_e.html+
 
   \item
    $s(M) > 2.5 \times 10^{9863}$ and
    $\sigma(M) > 6.9 \times 10^{4931}$.
    See study by H. Marxen in

 \verb+http://turbotm.de/~heiner/BB/simLig26_f.html+
\end{itemize}

 \item
 Terry and Shawn Ligocki found, in January 2008,
 a machine $M$ with $s(M) > 2.4 \times 10^{9866}$ and
 $\sigma(M) > 1.9 \times 10^{4933}$.
    See study by H. Marxen in

 \verb+http://turbotm.de/~heiner/BB/simLig26_g.html+

 It is the current record holder.

\end{itemize}

\begin{center}
\begin{tabular}{|c|c|c|c|}
\hline
February 2005 & T. and S. Ligocki & $s$ = 98,364,599 & $\sigma$ =
 10,574\\
\hline
April 2005 & T. and S. Ligocki & $s$ = 493,600,387 & $\sigma$ = 15,828\\
\hline
August 2005 & Lafitte, Papazian & \multicolumn{2}{c|}{superseded by a (2,5)-TM}\\
\hline
September 2007 & T. and S. Ligocki & $s > 2.3 \times 10^{54}$ & $\sigma
 > 1.9 \times 10^{27}$\\
\hline
October 2007 & T. and S. Ligocki & \multicolumn{2}{c|}{superseded by a (2,5)-TM}\\
\hline
November 2007 & T. and S. Ligocki & $s > 4.9 \times 10^{1643}$ &
 $\sigma > 8.6 \times 10^{821}$\\
 & & $s > 2.5 \times 10^{9863}$ &  $\sigma > 6.9 \times 10^{4931}$\\
\hline
January 2008 & T. and S. Ligocki & $s > 2.4 \times 10^{9866}$ & $\sigma
 > 1.9 \times 10^{4933}$\\
\hline
\end{tabular}
\end{center}

\bigskip

 The record holder and the past record holders:

\smallskip

\begin{tiny}
\begin{tabular}{cccccccccccccc}
A0  & A1  & A2  & A3  & A4  & A5  & B0  & B1  & B2  & B3  & B4  & B5 
 & $s(M)$ & $\sigma(M)$\\
1RB & 2LA & 1RH & 5LB & 5LA & 4LB & 1LA & 4RB & 3RB & 5LB & 1LB & 4RA
 & $> 2.4  \times 10^{9866}$ & $> 1.9 \times 10^{4933}$\\
1RB & 1LB & 3RA & 4LA & 2LA & 4LB & 2LA & 2RB & 3LB & 1LA & 5RA & 1RH
 & $> 2.5  \times 10^{9863}$ & $> 6.9 \times 10^{4931}$ \\
1RB & 2LB & 4RB & 1LA & 1RB & 1RH & 1LA & 3RA & 5RA & 4LB & 0RA & 4LA
 & $> 4.9  \times 10^{1643}$ & $> 8.6 \times 10^{821}$\\
1RB & 0RB & 3LA & 5LA & 1RH & 4LB & 1LA & 2RB & 3LA & 4LB & 3RB & 3RA
 & $> 2.3 \times 10^{54}$ & $> 1.9 \times 10^{27}$ \\
1RB & 2LA & 1RA & 1RA & 5LB & 4LB & 1LB & 1LA & 3RB & 4LA & 1RH & 3LA
 & 493,600,387 & 15,828 \\
1RB & 3LA & 3LA & 1RA & 1RA & 3LB & 1LB & 2LA & 2RA & 4RB & 5LB & 1RH
 & 98,364,599 & 10,249 \\
1RB & 3LA & 4LA & 1RA & 3RB & 1RH & 2LB & 1LA & 1LB & 3RB & 5RA & 1RH
 & 94,842,383 & 10,574 
\end{tabular}
\end{tiny}

\section{Behaviors of busy beavers}

\subsection{Introduction}

How do good machines behave?
We give below the tricks that allow them to reach high scores.

A \emph{configuration} of the Turing machine $M$ is a description of the tape.
The position of the tape head and the state are indicated by writing
together between parentheses the state and the symbol currently read
by the tape head.

For example, the initial configuration on a blank tape is:
$$\ldots 0(A0)0 \ldots$$
We denote by $a^k$ the string $a \ldots a$, $k$ times.
We write $C \vdash(t)\ D$ if the next move function leads from
configuration $C$ to configuration $D$ in $t$ computation steps.

\subsection{Turing machines with 5 states and 2 symbols}

\subsubsection{Marxen and Buntrock's champion}\label{sec:tm52a}

This machine is the record holder in the Busy Beaver Competition for
machines with 5 states and 2 symbols, since 1990.

It was analyzed by Buro (1990) (p.\ 64-67), and independently
by Michel (1993).

\bigskip

\begin{tabular}{cc}
\begin{tabular}{c}
Marxen and Buntrock (1990)\\
$s(M) = 47,176,870 =$? $S(5,2)$\\
$\sigma(M) = 4098 =$? $\mathnormal{\Sigma}(5,2)$
\end{tabular}
 & 
\begin{tabular}{c|cc|}
  &  0  &  1 \\ 
\hline
A & 1RB & 1LC\\
B & 1RC & 1RB\\
C & 1RD & 0LE\\
D & 1LA & 1LD\\
E & 1RH & 0LA\\
\hline
\end{tabular}
\end{tabular}

\bigskip

Let $C(n) = \ldots 0(A0)1^n0 \ldots$.

Then we have, for all $k \ge 0$,
$$\begin{array}{ccc}
C(3k)     & \vdash (5k^2 + 19k + 15) & C(5k + 6)\\
C(3k + 1) & \vdash (5k^2 + 25k + 27) & C(5k + 9)\\
C(3k + 2) & \vdash (6k + 12)         & \ldots 01(H0)1(001)^{k+1}10 \ldots
\end{array}$$

So we have:
$$\begin{array}{rl}
\multicolumn{2}{c}{\ldots 0(A0)0 \ldots =}\\
C( 0 ) & \vdash( 15 )\\
C( 6 ) & \vdash ( 73 )\\
C( 16 ) & \vdash ( 277 )\\
C( 34 ) & \vdash ( 907 )\\
C( 64 ) & \vdash ( 2,757 )\\
C( 114 ) & \vdash ( 7,957 )\\
C( 196 ) & \vdash ( 22,777 )\\
C( 334 ) & \vdash ( 64,407 )\\
C( 564 ) & \vdash ( 180,307 )\\
C( 946 ) & \vdash ( 504,027 )\\
C( 1,584 ) & \vdash ( 1,403,967 )\\
C( 2,646 ) & \vdash ( 3,906,393 )\\
C( 4,416 ) & \vdash ( 10,861,903 )\\
C( 7,366 ) & \vdash ( 30,196,527 )\\
C( 12,284 ) & \vdash ( 24,576 )\\
\multicolumn{2}{c}{\ldots 01(H0)1(001)^{4095}10 \ldots} 
\end{array}$$

\subsubsection{Marxen and Buntrock's runner-up}\label{sec:tm52b}

\begin{tabular}{cc}
\begin{tabular}{c}
Marxen and Buntrock (1990)\\
$s(M) = 23,554,764$\\
$\sigma(M) = 4097$
\end{tabular}
 & 
\begin{tabular}{c|cc|}
  &  0  &  1\\
\hline
A & 1RB & 0LD\\
B & 1LC & 1RD\\
C & 1LA & 1LC\\
D & 1RH & 1RE\\
E & 1RA & 0RB\\
\hline
\end{tabular}
\end{tabular}

\bigskip

Let $C(n) = \ldots 0(A0)1^n0 \ldots$.

Then we have, for all $k\ge 0$,
$$\begin{array}{ccc}
C(3k)     & \vdash ( 10k^2 + 10k + 4 )  & C(5k + 3)\\
C(3k + 1) & \vdash ( 3k + 3 )           & \ldots 01(110)^k11(H0)0 \ldots\\
C(3k + 2) & \vdash ( 10k^2 + 26k + 12 ) & C(5k + 7) 
\end{array}$$

So we have:

$$\begin{array}{rl}
\multicolumn{2}{c}{\ldots 0(A0)0 \ldots =}\\
C( 0 ) & \vdash ( 4 )\\
C( 3 ) & \vdash ( 24 )\\
C( 8 ) & \vdash ( 104 )\\
C( 17 ) & \vdash ( 392 )\\
C( 32 ) & \vdash ( 1,272 )\\
C( 57 ) & \vdash ( 3,804 )\\
C( 98 ) & \vdash ( 11,084 )\\
C( 167 ) & \vdash ( 31,692 )\\
C( 282 ) & \vdash ( 89,304 )\\
C( 473 ) & \vdash ( 250,584 )\\
C( 792 ) & \vdash ( 699,604 )\\
C( 1,323 ) & \vdash ( 1,949,224 )\\
C( 2,208 ) & \vdash ( 5,424,324 )\\
C( 3,683 ) & \vdash ( 15,087,204 )\\
C( 6,142 ) & \vdash ( 6,144 )\\
\multicolumn{2}{c}{\ldots 01(110)^{2047}11(H0)0 \ldots}
\end{array}$$

\subsection{Turing machines with 6 states and 2 symbols}

\subsubsection{Kropitz's machine found in May 2022}
\label{sec:tm62i}

This machine is the record holder in the Busy Beaver Competition
for machines with 6 states and 2 symbols, since May 2022.

\bigskip

\begin{tabular}{cc}
\begin{tabular}{c}
Kropitz (2022)\\  
$s(M)$ and $S(6,2) > 10^{\wedge\wedge}15$\\
$\sigma(M)$ and $\mathnormal{\Sigma}(6,2) > 10^{\wedge\wedge}15$
\end{tabular}
 &
\begin{tabular}{c|cc|}
 & 0 & 1\\
\hline
A & 1RB & 0LD\\
B & 1RC & 0RF\\
C & 1LC & 1LA\\
D & 0LE & 1RH\\
E & 1LF & 0RB\\
F & 0RC & 0RE\\
\hline
\end{tabular}
\end{tabular}

\bigskip

{\bf Analysis adapted from Pavel Kropitz and Shawn Ligocki}:

Let $C(n) = \ldots 010^n110^5(C0)0 \ldots$.

Then we have, for all $k \ge 0$,
$$\begin{array}{ccc}
\ldots 0(A0)0 \ldots & \vdash ( 45 ) & C(5)\\
C(4k) & \vdash ( (3 \times 9^{k+3} - 80 \times 3^{k+3} + 712k + 261)/32 )
 & \ldots01(H0)1^{n-1}0\ldots\\
C(4k + 1) & \vdash ( (3 \times 9^{k+3} - 64 \times 3^{k+3} + 712k + 981)/32 )
 & C((3^{k+3} - 11)/2)\\
C(4k + 2) & \vdash ( (3 \times 9^{k+3} - 64 \times 3^{k+3} + 712k + 1045)/32 )
 & C((3^{k+3} - 11)/2)\\
C(4k + 3) & \vdash ( (3 \times 9^{k+3} - 64 \times 3^{k+3} + 712k + 1749)/32 )
 & C((3^{k+3} + 1)/2)
\end{array}$$

with $n = (3^{k+3} - 11)/2$.

So we have (the final configuration is reached in 17 transitions):
$$\begin{array}{rl}
\ldots 0(A0)0 \ldots & \vdash ( 45 )\\
C( 5 )               & \vdash ( 506 )\\
C( 35 )              & \vdash ( 2,941,620,277 )\\
C( 88574 )          & \vdash ( )\\
\multicolumn{2}{c}{\cdots}\\
\multicolumn{2}{c}{\ldots 01(H0)1^n0 \ldots}
\end{array}$$
with $n > 10^{\wedge\wedge}15$.

\subsubsection{Kropitz's machine found in June 2010}
\label{sec:tm62h}

This machine was the record holder in the Busy Beaver Competition
for machines with 6 states and 2 symbols, from June 2010
to May 2022.

\bigskip

\begin{tabular}{cc}
\begin{tabular}{c}
Kropitz (2010)\\  
$s(M)$ and $S(6,2) > 7.4 \times 10^{36534}$\\
$\sigma(M)$ and $\mathnormal{\Sigma}(6,2) > 3.5 \times 10^{18267}$
\end{tabular}
 &
\begin{tabular}{c|cc|}
 & 0 & 1\\
\hline
A & 1RB & 1LE\\
B & 1RC & 1RF\\
C & 1LD & 0RB\\
D & 1RE & 0LC\\
E & 1LA & 0RD\\
F & 1RH & 1RC\\
\hline
\end{tabular}
\end{tabular}

\bigskip

Let $C(n) = \ldots 0(A0)1^n0 \ldots$.

Then we have, for all $k \ge 0$,
$$\begin{array}{ccc}
\ldots 0(A0)0 \ldots & \vdash ( 29 ) & C(9)\\
C(3k + 1) & \vdash ( 3k + 3 ) & \ldots 0111(011)^k(H0)0 \ldots\\
C(9k + 9) & \vdash ( (125 \times 16^{k+2} + 325 \times 4^{k+2} + 228k - 2289)/27 )
 & C((50 \times 4^{k+1} - 11)/3)\\
C(9k + 12) & \vdash ( (125 \times 16^{k+2} + 325 \times 4^{k+2} + 228k - 912)/27 )
 & C((50 \times 4^{k+1} + 1)/3)
\end{array}$$

So we have:
$$\begin{array}{rl}
\ldots 0(A0)0 \ldots & \vdash ( 29 )\\
C( 9 )               & \vdash ( 1293 )\\
C( 63 )              & \vdash ( 19,884,896,677 )\\
C( 273063 )          & \vdash ( 125 \times 16^{30341} + 325 \times 4^{30341} + 6,916,380)/27 )\\
C( 50 \times 4^{30340} + 1)/3 )     & \vdash ( 50 \times 4^{30340} + 7)/3 )\\
\multicolumn{2}{c}{\ldots 0111(011)^p(H0)0 \ldots}
\end{array}$$
with $p = (50 \times 4^{30340} - 2)/9$.

So the total time is
$s(M) = (125 \times 16^{30341} + 1750 \times 4^{30340} + 15)/27 + 19,885,154,163$,
and the final number of 1 is
$\sigma(M) = (25 \times 4^{30341} + 23)/9$.

\smallskip

Some configurations take a long time to halt.
For example, $C(2) \vdash(t)$ END
with $t > 10^{10^{10^{10^{18,705,352}}}}$. A proof of this fact is given by
``Cloudy176'' in

\noindent
\verb+http://googology.wikia.com/wiki/User_blog:Cloudy176/Proving_the_bound_for_S(7)+

This property was used by ``Wythagoras'', in March 2014, to define
a (7,2)-TM that extends the present (6,2)-TM and enters this
configuration C(2) in two steps. See

\noindent
\verb+http://googology.wikia.com/wiki/User_blog:Wythagoras/A_good_bound_for_S(7)%3F+

See detailed analysis in Michel (2015), Section 6.

\subsubsection{Kropitz's machine found in May 2010}
\label{sec:tm62g}

This machine was the record holder in the Busy Beaver Competition
for machines with 6 states and 2 symbols, from May 2010 to June 2010.

\bigskip

\begin{tabular}{cc}
\begin{tabular}{c}
Kropitz (2010)\\  
$s(M) > 3.8 \times 10^{21132}$\\
$\sigma(M) > 3.1 \times 10^{10566}$
\end{tabular}
 &
\begin{tabular}{c|cc|}
 & 0 & 1\\
\hline
A & 1RB & 0LD\\
B & 1RC & 0RF\\
C & 1LC & 1LA\\
D & 0LE & 1RH\\
E & 1LA & 0RB\\
F & 0RC & 0RE\\
\hline
\end{tabular}
\end{tabular}

\bigskip

{\bf Analysis adapted from Shawn Ligocki}:

Let $C(n, k) = \ldots 010^n1(C1)1^{3k}0 \ldots$.

Then we have, for all $k \ge 0$, all $n \ge 0$,
$$\begin{array}{ccc}
\ldots 0(A0)0 \ldots & \vdash(47) & C(5, 2)\\
C(0, k)     & \vdash(3)                  & \ldots 01(H0)1^{3k+1}0 \dots\\
C(1, k)     & \vdash(3k + 37)            & C(3k + 2, 2)\\
C(2, k)     & \vdash(12k + 44)           & C(4, k + 2)\\
C(3, k)     & \vdash(3k + 57)            & C(3k + 8, 2)\\
C(n + 4, k) & \vdash(27k^2 + 105k + 112) & C(n, 3k + 5)
\end{array}$$

\bigskip

So we have (the final configuration is reached in 22158 transitions):
$$\begin{array}{rl}
\ldots 0(A0)0 \ldots & \vdash(47)\\
C( 5, 2 ) & \vdash(430)\\
C( 1, 11 ) & \vdash(70)\\
C( 35, 2 ) & \vdash(430)\\
C( 31, 11 ) & \vdash(4,534)\\
C( 27, 38 ) & \vdash(43,090)\\
C( 23, 119 ) & \vdash(394,954)\\
C( 19, 362 ) & \vdash(3,576,310)\\
C( 15, 1091 ) & \vdash(32,252,254)\\
C( 11, 3278 ) & \vdash(290,466,970)\\
C( 7, 9839 ) & \vdash(2,614,793,074)\\
C( 3, 29522 ) & \vdash(88,623)\\
C( 88574, 2 ) & \vdash(430)\\
C( 88570, 11 ) & \vdash(4,534)\\
C( 88566, 38 ) & \vdash(43,090)\\
\multicolumn{2}{c}{\cdots}
\end{array}$$

Note that $C(4n + r, 2) \vdash(t_n)\ C(r, u_n)$,
with $u_n = (3^{n+2} - 5)/2$,
and $t_n = (3 \times 9^{n+3} - 80 \times 3^{n+3} + 584 n - 27)/32$.

\smallskip

Some configurations take a long time to halt.
For example, $C(1, 9) \vdash(t)$ END
with $t > 10^{10^{10^{10^{10^{3520}}}}}$.

See detailed analysis in Michel (2015), Section 7.

\subsubsection{Ligockis' machine found in December 2007}
\label{sec:tm62e}

This machine was the record holder in the Busy Beaver Competition for
machines with 6 states and 2 symbols, from December 2007 to May 2010.

\bigskip

\begin{tabular}{cc}
\begin{tabular}{c}
Terry and Shawn Ligocki (2007)\\
$s(M) > 2.5 \times 10^{2879}$\\
$\sigma(M) > 4.6 \times 10^{1439}$
\end{tabular}
 & 
\begin{tabular}{c|cc|}
  &  0  &  1\\
\hline
A & 1RB & 0LE\\
B & 1LC & 0RA\\
C & 1LD & 0RC\\
D & 1LE & 0LF\\
E & 1LA & 1LC\\
F & 1LE & 1RH\\
\hline
\end{tabular}
\end{tabular}

\bigskip

Let $C(n, p) = \ldots 0(A0)(10)^nR(\bin(p))0 \ldots$, where $R(\bin(p))$ is
the number $p$ written in binary in reverse order, so that
$C(n, 4m + 1) = C(n + 1, m)$.
The number of transitions between configurations $C(n, p)$ is infinite,
but only 18 transitions are used in the computation on a blank tape.
For all $m \ge 0$, all $k \ge 0$,
$$\begin{array}{ccc}
C(k, 4m + 3) & \vdash ( 4k + 6 ) & C(k + 2, m)\\
C(2k + 1, 4m) & \vdash ( 6k^2 + 52k + 98 ) & C(3k + 8, m)\\
C(4k, 4m) & \vdash ( 24k^2 + 36k + 13 ) & C(6k + 2, 2m + 1)\\
C(4k + 2, 4m) & \vdash ( 24k^2 + 60k + 27 ) & C(6k + 2, 128m + 86)\\
C(k, 8m + 2) & \vdash ( 4k + 14 ) & C(k + 2, 2m + 1)\\
C(2k + 1, 32m + 22) & \vdash ( 6k^2 + 64k + 160 ) & C(3k + 10, 2m + 1)\\
C(4k, 32m + 22) & \vdash ( 24k^2 + 36k + 29 ) & C(6k + 4, m)\\
C(4k + 2, 32m + 22) & \vdash ( 24k^2 + 60k + 43 ) & C(6k + 2, 1024m + 342)\\
C(k, 64m + 46) & \vdash ( 4k + 30 ) & C(k + 4, m)\\
C(k + 1, 128m + 6) & \vdash ( 8k + 66 ) & C(k + 6, 2m + 1)\\
C(2k, 256m + 14) & \vdash ( 6k^2 + 64k + 172 ) & C(3k + 11, m)\\
C(4k + 1, 256m + 14) & \vdash ( 24k^2 + 84k + 89 ) & C(6k + 8, 2m + 1)\\
C(4k + 3, 256m + 14) & \vdash ( 24k^2 + 108k + 127 ) & C(6k + 8, 128m + 86)\\
C(4k, 512m + 30) & \vdash ( 24k^2 + 156k + 173 ) & C(6k + 11, m)\\
C(4k + 2, 512m + 30) & \vdash ( 24k^2 + 60k + 57 ) & C(6k + 2, 16384m + 11134)\\
C(4k + 2, 131072m + 11134) & \vdash ( 24k^2 + 60k + 89 ) & C(6k + 2, 4194304m + 2848638)\\
C(4k, 131072m + 96126) & \vdash ( 24k^2 + 36k + 109 ) & C(6k + 10, m)\\
C(k + 1, 512m + 94) & \vdash ( 2k + 61 ) & . . . 0(10)^k1(H0)1110110101R(\bin(m))0 . . . 
\end{array}$$

\bigskip

So we have (the final configuration is reached in 11026 transitions):
$$\begin{array}{rl}
\multicolumn{2}{c}{\ldots 0(A0)0 \ldots =}\\
C( 0, 0 ) & \vdash ( 13 )\\
C( 3, 0 ) & \vdash ( 156 )\\
C( 11, 0 ) & \vdash ( 508 )\\
C( 23, 0 ) & \vdash ( 1396 )\\
C( 41, 0 ) & \vdash ( 3538 )\\
C( 68, 0 ) & \vdash ( 7,561 )\\
C( 105, 0 ) & \vdash ( 19,026 )\\
C( 164, 0 ) & \vdash ( 41,833 )\\
C( 249, 0 ) & \vdash ( 98,802 )\\
C( 380, 0 ) & \vdash ( 220,033 )\\
C( 573, 0 ) & \vdash ( 505,746 )\\
C( 866, 0 ) & \vdash ( 1,132,731 )\\
C( 1298, 86 ) & \vdash ( 2,538,907 )\\
C( 1946, 2390 ) & \vdash ( 5,697,907 )\\
C( 2918, 76118 ) & \vdash ( 12,798,367 )\\
C( 4376, 2435414 ) & \vdash ( 1,034,066,333 )\\
C( 6568, 76106 ) & \vdash ( 26,286 )\\
C( 6570, 19027 ) & \vdash ( 26,286 )\\
C( 6572, 4756 ) & \vdash ( 64,845,937 )\\
C( 9860, 2379 ) & \vdash ( 39,446 )\\
C( 9862, 594 ) & \vdash ( 39,462 )\\
C( 9867, 2 ) & \vdash ( 39,482 )\\
\multicolumn{2}{c}{\cdots}
\end{array}$$

\subsubsection{Ligockis' machine found in November 2007}
\label{sec:tm62f}

This machine was the record holder in the Busy Beaver Competition for
machines with 6 states and 2 symbols, from November to December 2007.

\bigskip

\begin{tabular}{cc}
\begin{tabular}{c}
Terry and Shawn Ligocki (2007)\\
$s(M) > 8.9 \times 10^{1762}$\\
$\sigma(M) > 2.5 \times 10^{881}$
\end{tabular}
 & 
\begin{tabular}{c|cc|}
  &  0  &  1\\
\hline
A & 1RB & 0RF\\
B & 0LB & 1LC\\
C & 1LD & 0RC\\
D & 1LE & 1RH\\
E & 1LF & 0LD\\
F & 1RA & 0LE\\
\hline
\end{tabular}
\end{tabular}

\bigskip

Let $C(n, p) = \ldots 0(F0)(10)^nR(\bin(p))0 \ldots$, where $R(\bin(p))$ is
the number $p$ written in binary in reverse order, so that
$C(n, 4m + 1) = C(n + 1, m)$.
The number of transitions between configurations $C(n, p)$ is infinite,
but only 12 transitions are used in the computation on a blank tape.
For all $m \ge 0$, all $k \ge 0$,
$$\begin{array}{ccc}
\ldots 0(A0)0\ldots & \vdash ( 6 ) & C(0, 15)\\
C(k, 4m + 3) & \vdash ( 4k + 6 ) & C(k + 2, m)\\
C(2k, 4m) & \vdash ( 30k^2 + 20k + 15 ) & C(5k + 2, 2m + 1)\\
C(2k + 1, 4m) & \vdash ( 30k^2 + 40k + 25 ) & C(5k + 2, 32m + 20)\\
C(k, 8m + 2) & \vdash ( 8k + 20 ) & C(k + 3, 2m + 1)\\
C(2k, 16m + 6) & \vdash ( 30k^2 + 40k + 23 ) & C(5k + 2, 32m + 20)\\
C(2k + 1, 16m + 6) & \vdash ( 30k^2 + 80k + 63 ) & C(5k + 7, 2m + 1)\\
C(k, 32m + 14) & \vdash ( 4k + 18 ) & C(k + 3, 2m + 1)\\
C(2k, 128m + 94) & \vdash ( 30k^2 + 40k + 39 ) & C(5k + 2, 256m + 84)\\
C(2k + 1, 128m + 94) & \vdash ( 30k^2 + 80k + 79 ) & C(5k + 9, m)\\
C(k, 256m + 190) & \vdash ( 4k + 34 ) & C(k + 5, m)\\
C(k, 512m + 30) & \vdash ( 2k + 43 ) & \ldots 0(10)^k1(H0)10100101R(\bin(m))0 \ldots 
\end{array}$$

\bigskip

So we have (the final configuration is reached in 3346 transitions):
$$\begin{array}{rl}
\ldots 0(A0)0\ldots & \vdash ( 6 )\\
C( 0, 15 ) & \vdash ( 6 )\\
C( 2, 3 ) & \vdash ( 14 )\\
C( 4, 0 ) & \vdash ( 175 )\\
C( 13, 0 ) & \vdash ( 1,345 )\\
C( 32, 20 ) & \vdash ( 8,015 )\\
C( 82, 11 ) & \vdash ( 334 )\\
C( 84, 2 ) & \vdash ( 692 )\\
C( 88, 0 ) & \vdash ( 58,975 )\\
C( 223, 0 ) & \vdash ( 374,095 )\\
C( 557, 20 ) & \vdash ( 2,329,665 )\\
C( 1392, 180 ) & \vdash ( 14,546,415 )\\
C( 3482, 91 ) & \vdash ( 13,934 )\\
C( 3484, 22 ) & \vdash ( 91,106,623 )\\
C( 8712, 52 ) & \vdash ( 569,329,215 )\\
C( 21782, 27 ) & \vdash ( 87,134 )\\
C( 21784, 6 ) & \vdash ( 3,559,505,623 )\\
C( 54462, 20 ) & \vdash ( 22,246,365,465 )\\
C( 136157, 11 ) & \vdash ( 544,634 )\\
C( 136159, 2 ) & \vdash ( 1,089,292 )\\
C( 136163, 0 ) & \vdash ( 139,053,400,095 )\\
C( 340407, 20 ) & \vdash ( 869,078,644,415 )\\
\multicolumn{2}{c}{\cdots}
\end{array}$$

See detailed analysis in Michel (2015), Section 8.

\subsubsection{Marxen and Buntrock's machine found in March 2001}
\label{sec:tm62a}

This machine was the record holder in the Busy Beaver Competition for
machines with 6 states and 2 symbols, from March 2001 to November 2007.

\bigskip

\begin{tabular}{cc}
\begin{tabular}{c}
Marxen and Buntrock (2001)\\
$s(M) > 3.0 \times 10^{1730}$\\
$\sigma(M) > 1.2 \times 10^{865}$
\end{tabular}
 & 
\begin{tabular}{c|cc|}
  &  0  &  1\\
\hline
A & 1RB & 0LF\\
B & 0RC & 0RD\\
C & 1LD & 1RE\\
D & 0LE & 0LD\\
E & 0RA & 1RC\\
F & 1LA & 1RH\\
\hline
\end{tabular}
\end{tabular}

\bigskip

Let $C(n, p) = \ldots 0(A0)(01)^nR(\bin(p))0 \ldots$, where $R(\bin(p))$ is
the number $p$ written in binary in reverse order, so that
$C(n, 4m + 2) = C(n + 1, m)$.
The number of transitions between configurations $C(n, p)$ is infinite,
but only 20 transitions are used in the computation on a blank tape.
For all $m \ge 0$, all $k \ge 0$,
$$\begin{array}{ccc}
C(2k, 4m) & \vdash (9k^2 + 25k + 9) & C(3k + 1, 2m + 1)\\
C(2k, 16m + 1) & \vdash (9k^2 + 25k + 17) & C(3k + 2, 2m + 1)\\
C(2k, 4m + 3) & \vdash (9k^2 + 25k + 9) & C(3k + 1, 2m)\\
C(2k, 64m + 53) & \vdash (9k^2 + 25k + 25) & C(3k + 3, 2m)\\
C(2k, 256m + 9) & \vdash (9k^2 + 25k + 29) & C(3k + 4, 2m + 1)\\
C(2k, 1024m + 57) & \vdash (9k^2 + 25k + 33) & C(3k + 2, 128m + 104)\\
C(2k, 1024m + 85) & \vdash (9k^2 + 25k + 41) & C(3k + 5, 2m + 1)\\
C(2k + 1, 16m) & \vdash (9k^2 + 25k + 21) & C(3k + 3, 2m + 1)\\
C(2k + 1, 4m + 1) & \vdash (9k^2 + 25k + 13) & C(3k + 1, 8m + 4)\\
C(2k + 1, 64m + 4) & \vdash (9k^2 + 25k + 29) & C(3k + 4, 2m + 1)\\
C(2k + 1, 64m + 3) & \vdash (9k^2 + 25k + 25) & C(3k + 1, 128m + 104)\\
C(2k + 1, 1024m + 104) & \vdash (9k^2 + 43k + 75) & C(3k + 7, 2m + 1)\\
C(2k + 1, 16m + 12) & \vdash (9k^2 + 25k + 21) & C(3k + 3, 2m)\\
C(2k + 1, 16m + 7) & \vdash (9k^2 + 25k + 17) & C(3k + 1, 32m + 16)\\
C(2k + 1, 256m + 15) & \vdash (9k^2 + 25k + 29) & C(3k + 1, 512m + 416)\\
C(2k + 1, 64m + 52) & \vdash (9k^2 + 25k + 29) & C(3k + 4, 2m)\\
C(2k + 1, 256m + 20) & \vdash (9k^2 + 25k + 37) & C(3k + 5, 2m + 1)\\
C(2k + 1, 4096m + 420) & \vdash (9k^2 + 43k + 89) & C(3k + 8, 2m + 1)\\
C(2k + 1, 256m + 211) & \vdash (9k^2 + 25k + 33) & C(3k + 1, 512m + 168)\\
C(2k + 1, 16m + 11) & \vdash (9k^2 + 13k + 10) & \ldots 0(10)^{3k+1}11(H0)10R(\bin(m))0 \ldots
\end{array}$$

\bigskip

So we have (the final configuration is reached in 4911 transitions):
$$\begin{array}{rl}
\multicolumn{2}{c}{\ldots 0(A0)0 \ldots =}\\
C( 0, 0 ) & \vdash ( 9 )\\
C( 1, 1 ) & \vdash ( 13 )\\
C( 1, 4 ) & \vdash ( 29 )\\
C( 4, 1 ) & \vdash ( 103 )\\
C( 8, 1 ) & \vdash ( 261 )\\
C( 14, 1 ) & \vdash ( 633 )\\
C( 23, 1 ) & \vdash ( 1,377 )\\
C( 34, 4 ) & \vdash ( 3,035 )\\
C( 52, 3 ) & \vdash ( 6,743 )\\
C( 79, 0 ) & \vdash ( 14,685 )\\
C( 120, 1 ) & \vdash ( 33,917 )\\
C( 182, 1 ) & \vdash ( 76,821 )\\
C( 275, 1 ) & \vdash ( 172,359 )\\
C( 412, 4 ) & \vdash ( 387,083 )\\
C( 619, 3 ) & \vdash ( 867,079 )\\
C( 928, 104 ) & \vdash ( 1,949,273 )\\
C( 1393, 53 ) & \vdash ( 4,377,157 )\\
C( 2089, 108 ) & \vdash ( 9,835,545 )\\
C( 3135, 12 ) & \vdash ( 22,138,597 )\\
C( 4704, 0 ) & \vdash ( 49,845,945 )\\
C( 7057, 1 ) & \vdash ( 112,109,269 )\\
C( 10585, 4 ) & \vdash ( 252,179,705 )\\
\multicolumn{2}{c}{\cdots}
\end{array}$$

\bigskip

\noindent{\bf Note}: Clive Tooth posted an analysis of this machine on Google Groups
(sci.math$>$The Turing machine known as \#r), on June 28, 2002.
He used the configurations $S(n, x) = \ldots 0101(B1)010(01)^nx0 \ldots$
His analysis can be easily connected to the present one, by noting that
$$C(n, p) \vdash ( 15 )\ S(n - 2, R(\bin(p))).$$ 

\subsubsection{Marxen and Buntrock's second machine}\label{sec:tm62b}

This machine was the record holder in the Busy Beaver Competition for
machines with 6 states and 2 symbols, from October 2000 to March 2001.

\bigskip

\begin{tabular}{cc}
\begin{tabular}{c}
Marxen and Buntrock (2000)\\
$s(M) > 6.1 \times 10^{925}$\\
$\sigma(M) > 6.4 \times 10^{462}$
\end{tabular}
 & 
\begin{tabular}{c|cc|}
  &  0  &  1\\
\hline
A & 1RB & 0LB\\
B & 0RC & 1LB\\
C & 1RD & 0LA\\
D & 1LE & 1LF\\
E & 1LA & 0LD\\
F & 1RH & 1LE\\
\hline
\end{tabular}
\end{tabular}

\bigskip

Let $C(n) = \ldots 01^n(B0)0 \ldots$.

Then we have, for all $k \ge 0$,
$$\begin{array}{ccc}
\ldots 0(A0)0 \ldots & \vdash ( 1 ) & C(1)\\
C(3k)     & \vdash ( 54 \times 4^{k+1} - 27 \times 2^{k+3} + 26k + 86 ) & C(9 \times 2^{k+1} - 8)\\
C(3k + 1) & \vdash ( 2048 \times (4^k-1)/3 - 3 \times 2^{k+7} + 26k + 792 ) & C(2^{k+5} - 8)\\
C(3k + 2) & \vdash ( 3k + 8 ) & \ldots 01(H1)(011)^k(0101)0 \ldots
\end{array}$$

So we have:
$$\begin{array}{rl}
\ldots 0(A0)0 \ldots & \vdash ( 1 )\\
C( 1 )               & \vdash ( 408 )\\
C( 24 )              & \vdash ( 14,100,774 )\\
C( 4600 )            & \vdash ( 2048 \times (4^{1533} - 1)/3 - 3 \times 2^{1540} + 40650 )\\
C( 2^{1538} - 8 )     & \vdash ( 2^{1538} - 2 )\\
\multicolumn{2}{c}{\ldots 01(H1)(011)^p(0101)0 \ldots}
\end{array}$$
with $p = (2^{1538} - 10)/3$.

So the total time is $T = 2048 \times (4^{1533} - 1)/3 - 11 \times 2^{1538} + 14141831$,
and the final number of 1 is $2 \times (2^{1538} - 10)/3 + 4$.

Note that
$$C(6k + 1) \vdash (\ )\ C(3m) \vdash (\ )\ C(6p + 4) \vdash (\ )\ C(3q + 2) \vdash (\ )\mbox{ END},$$
with $m = (2^{2k+5} - 8)/3$,
$p = 3 \times 2^m - 2$,
$q = (2^{2p+6} - 10)/3$.

So all configurations $C(n)$ lead to a halting configuration.
Those taking the most time are $C(6k + 1)$.
For example: $$C(7) \vdash ( t )\mbox{ END\quad with}\quad t > 10^{3.9 \times 10^{12}}.$$
More generally: $$C(6k + 1) \vdash (t(k))\mbox{ END\quad with}\quad t(k) > 10^{10^{10^{(3k + 2)/5}}}.$$ 

\noindent See also the analyses by Robert Munafo: the short one in

\verb+http://mrob.com/pub/math/ln-notes1-4.html#mb6q+

\noindent and the detailed one in

\verb+http://mrob.com/pub/math/ln-mb6q.html+

\noindent See detailed analysis in Michel (2015), Section 9.

\subsubsection{Marxen and Buntrock's third machine}\label{sec:tm62c}

\begin{tabular}{cc}
\begin{tabular}{c}
Marxen and Buntrock (2000)\\
$s(M) > 6.1 \times 10^{119}$\\
$\sigma(M) > 1.4 \times 10^{60}$
\end{tabular}
 & 
\begin{tabular}{c|cc|}
  &  0  &  1\\
\hline
A & 1RB & 0LC\\
B & 1LA & 1RC\\
C & 1RA & 0LD\\
D & 1LE & 1LC\\
E & 1RF & 1RH\\
F & 1RA & 1RE\\
\hline
\end{tabular}
\end{tabular}

\bigskip

Let $C(n, x) = \ldots 0(E0)1000(10)^nx0 \ldots$, so that $C(n, 10y) = C(n + 1, y)$.
The number of transitions between configurations $C(n, x)$ is infinite,
but only 9 transitions are used in the computation on a blank tape.
For all $k \ge 0$,
$$\begin{array}{ccc}
\ldots 0(A0)0 \ldots & \vdash ( 18 )              & C(1, 01)\\
C(2k, 01^n)          & \vdash ( 6k^2 + 22k + 15 ) & C(3k + 1, 01^{n+1})\\
C(2k, 11)            & \vdash ( 6k^2 + 34k + 41 ) & C(3k + 4, 01)\\
C(2k, 111)           & \vdash ( 6k^2 + 34k + 45 ) & C(3k + 5, 01)\\
C(2k, 1111)          & \vdash ( 6k^2 + 28k + 25 ) & \ldots 01^{6k+11}(H0)0 \ldots\\
C(2k + 1, 0)         & \vdash ( 6k^2 + 34k + 43 ) & C(3k + 4, 0)\\
C(2k + 1, 01)        & \vdash ( 6k^2 + 22k + 27 ) & C(3k + 4, 01)\\
C(2k + 1, 01^{n+2})   & \vdash ( 6k^2 + 22k + 23 ) & C(3k + 4, 1^n)\\
C(2k + 1, 1^{n+2})    & \vdash ( 6k^2 + 34k + 41 ) & C(3k + 4, 01^n)\\
\end{array}$$

So we have (the final configuration is reached in 337 transitions):
$$\begin{array}{rl}
\ldots 0(A0)0 \ldots & \vdash ( 18 )\\
C( 1, 01 ) & \vdash ( 27 )\\
C( 4, 01 ) & \vdash ( 83 )\\
C( 7, 011 ) & \vdash ( 143 )\\
C( 13, 0 ) & \vdash ( 463 )\\
C( 22, 0 ) & \vdash ( 983 )\\
C( 34, 01 ) & \vdash ( 2,123 )\\
C( 52, 011 ) & \vdash ( 4,643 )\\
C( 79, 0111 ) & \vdash ( 10,007 )\\
C( 122, 0 ) & \vdash ( 23,683 )\\
C( 184, 01 ) & \vdash ( 52,823 )\\
C( 277, 011 ) & \vdash ( 117,323 )\\
C( 418, 0 ) & \vdash ( 266,699 )\\
C( 628, 01 ) & \vdash ( 598,499 )\\
C( 943, 011 ) & \vdash ( 1,341,431 )\\
C( 1417, 0 ) & \vdash ( 3,031,699 )\\
C( 2128, 0 ) & \vdash ( 6,815,999 )\\
C( 3193, 01 ) & \vdash ( 15,318,435 )\\
C( 4792, 01 ) & \vdash ( 34,497,623 )\\
C( 7189, 011 ) & \vdash ( 77,580,107 )\\
C( 10786, 0 ) & \vdash ( 174,625,355 )\\
\multicolumn{2}{c}{\cdots}
\end{array}$$

Note that, if $C(n, m) = \ldots 0(E0)1000(10)^nR(\bin(m))0 \ldots$,
where $R(\bin(m))$ is the number $m$ written in binary in reverse order,
so that $C(n, 4m + 1) = C(n + 1, m)$, then we have also, for all $k$, $m \ge 0$,
$$\begin{array}{ccc}
\ldots 0(A0)0 \ldots & \vdash ( 18 )              & C(1, 2)\\
C(2k, 2m)            & \vdash ( 6k^2 + 22k + 15 ) & C(3k + 1, 4m + 2)\\
C(2k, 32m + 3)       & \vdash ( 6k^2 + 34k + 41 ) & C(3k + 4, 4m + 2)\\
C(2k, 128m + 7)      & \vdash ( 6k^2 + 34k + 45 ) & C(3k + 5, 4m + 2)\\
C(2k, 32m + 15)      & \vdash ( 6k^2 + 28k + 25 ) & \ldots 01^{6k+11}(H0)R(\bin(m))0 \ldots\\
C(2k + 1, 4m)        & \vdash ( 6k^2 + 34k + 43 ) & C(3k + 4, 2m)\\
C(2k + 1, 32m + 2)   & \vdash ( 6k^2 + 22k + 27 ) & C(3k + 4, 4m + 2)\\
C(2k + 1, 8m + 6)    & \vdash ( 6k^2 + 22k + 23 ) & C(3k + 4, m)\\
C(2k + 1, 4m + 3)    & \vdash ( 6k^2 + 34k + 41 ) & C(3k + 4, 2m)
\end{array}$$

\subsubsection{Another Marxen and Buntrock's machine}\label{sec:tm62d}

This machine was discovered in January 1990, and was published
on the web (Google groups) on September 3, 1997.
It was the record holder in the Busy Beaver Competition for
machines with 6 states and 2 symbols up to July 2000.

\bigskip

\begin{tabular}{cc}
\begin{tabular}{c}
Marxen and Buntrock (1997)\\
$s(M) = 8,690,333,381,690,951$\\
$\sigma(M) = 95,524,079$
\end{tabular}
 & 
\begin{tabular}{c|cc|}
  &  0   &  1\\
\hline
A & 1RB & 1RA\\
B & 1LC & 1LB\\
C & 0RF & 1LD\\
D & 1RA & 0LE\\
E & 1RH & 1LF\\
F & 0LA & 0LC\\
\hline
\end{tabular}
\end{tabular}

\bigskip

Note the likeness to the machine $N$ with 3 states and 3 symbols discovered,
in August 2006, by Terry and Shawn Ligocki, and studied in Section \ref{sec:tm33g}.
For this machine $N$, we have $s(N) = 4,345,166,620,336,565$ and
$\sigma(N)= 95,524,079$, that is, same value of $\sigma$, and almost
half the value of $s$. See analysis of this similarity in Section \ref{sec:sb}. 

\bigskip

\noindent{\bf Analysis by Robert Munafo}:

Let $C(n) = \ldots 0(D0)1^n0 \ldots$. 

Then we have, for all $k \ge 0$,
$$\begin{array}{ccc}
\ldots 0(A0)0 \ldots  & \vdash ( 3 )                & C(2)\\
C(4k)                 & \vdash ( 8k + 6 )           & \ldots 01(H0)(10)^{2k}110 \ldots\\ 
C(4k + 1)             & \vdash ( 20k^2 + 56k + 30 ) & C(10k + 9)\\
C(4k + 2)             & \vdash ( 20k^2 + 56k + 33 ) & C(10k + 9)\\
C(4k + 3)             & \vdash ( 20k^2 + 68k + 51 ) & C(10k + 12)
\end{array}$$

So we have:
$$\begin{array}{rl}
\ldots 0(A0)0 \ldots & \vdash ( 3 )\\
C( 2 ) & \vdash ( 33 )\\
C( 9 ) & \vdash ( 222 )\\
C( 29 ) & \vdash ( 1,402 )\\
C( 79 ) & \vdash ( 8,563 )\\
C( 202 ) & \vdash ( 52,833 )\\
C( 509 ) & \vdash ( 329,722 )\\
C( 1,279 ) & \vdash ( 2,056,963 )\\
C( 3,202 ) & \vdash ( 12,844,833 )\\
C( 8,009 ) & \vdash ( 80,272,222 )\\
C( 20,029 ) & \vdash ( 501,681,402 )\\
C( 50,079 ) & \vdash ( 3,135,358,563 )\\
C( 125,202 ) & \vdash ( 19,595,552,833 )\\
C( 313,009 ) & \vdash ( 122,471,892,222 )\\
C( 782,529 ) & \vdash ( 765,448,543,902 )\\
C( 1,956,329 ) & \vdash ( 4,784,051,443,102 )\\
C( 4,890,829 ) & \vdash ( 29,900,316,628,602 )\\
C( 12,227,079 ) & \vdash ( 186,876,942,247,563 )\\
C( 30,567,702 ) & \vdash ( 1,167,980,782,060,333 )\\
C( 76,419,259 ) & \vdash ( 7,299,879,658,619,323 )\\
C( 191,048,152 ) & \vdash ( 382,096,310 )\\
\multicolumn{2}{c}{\ldots 01(H0)(10)^{95524076}110 \ldots}
\end{array}$$

\subsection{Turing machines with 3 states and 3 symbols}

\subsubsection{Ligockis' champion}\label{sec:tm33h}

This machine is the record holder in the Busy Beaver Competition
for machines with 3 states and 3 symbols, since November 2007.

\bigskip

\begin{tabular}{cc}
\begin{tabular}{c}
Terry and Shawn Ligocki (2007)\\
$s(M) = 119,112,334,170,342,540 =$? $S(3,3)$\\
$\sigma(M) = 374,676,383 =$? $\mathnormal{\Sigma}(3,3)$
\end{tabular}
 & 
\begin{tabular}{c|ccc|}
  &  0  &  1  &  2\\
\hline
A & 1RB & 2LA &	1LC\\
B & 0LA & 2RB &	1LB\\
C & 1RH & 1RA &	1RC\\
\hline
\end{tabular}
\end{tabular}

\bigskip

Let $C(n) = \ldots 0(A0)2^n0 \ldots$.

Then we have, for all $k \ge 0$,
$$\begin{array}{ccc}
\ldots 0(A0)0 \ldots &   \vdash ( 3 )      & C(1)\\
C(8k + 1) & \vdash ( 112k^2 + 116k + 13 )  & C(14k + 3)\\
C(8k + 2) & \vdash ( 112k^2 + 144k + 38 )  & C(14k + 7)\\
C(8k + 3) & \vdash ( 112k^2 + 172k + 54 )  & C(14k + 8)\\
C(8k + 4) & \vdash ( 112k^2 + 200k + 74 )  & C(14k + 9)\\
C(8k + 5) & \vdash ( 112k^2 + 228k + 97 )  & \ldots  01(H1)2^{14k+9}0 \ldots\\
C(8k + 6) & \vdash ( 112k^2 + 256k + 139 ) & C(14k + 14)\\
C(8k + 7) & \vdash ( 112k^2 + 284k + 169 ) & C(14k + 15)\\
C(8k + 8) & \vdash ( 112k^2 + 312k + 203 ) & C(14k + 16)
\end{array}$$

So we have (in 34 transitions):
$$\begin{array}{rl}
 \ldots  0(A0)0 \ldots & \vdash ( 3 )\\
C( 1 ) & \vdash ( 13 )\\
C( 3 ) & \vdash ( 54 )\\
C( 8 ) & \vdash ( 203 )\\
C( 16 ) & \vdash ( 627 )\\
C( 30 ) & \vdash ( 1915 )\\
\multicolumn{2}{c}{\cdots}\\
C( 122,343,306 ) & \vdash ( 26,193,799,261,043,238 )\\
C( 214,100,789 ) & \vdash ( 80,218,511,093,348,089 )\\
\multicolumn{2}{c}{\ldots 01(H1)2^{374676381}0 \ldots} 
\end{array}$$

See detailed analysis in Michel (2015), Section 3.

\subsubsection{Ligockis' machine found in August 2006}\label{sec:tm33g}

This machine was the record holder in the Busy Beaver Competition for
machines with 3 states and 3 symbols, from August 2006 to November 2007.

\bigskip

\begin{tabular}{cc}
\begin{tabular}{c}
Terry and Shawn Ligocki (2006)\\
$s(M) = 4,345,166,620,336,565$\\
$\sigma(M) = 95,524,079$
\end{tabular}
 & 
\begin{tabular}{c|ccc|}
  &	 0  &	 1  &  	 2\\
\hline
A &	1RB &	2RC &	1LA\\
B &	2LA &	1RB &	1RH\\
C &	2RB &	2RA &	1LC\\
\hline
\end{tabular}
\end{tabular}

\bigskip

Note the likeness to the machine $N$ with 6 states and 2 symbols discovered,
in January 1990, by Heiner Marxen and J\"{u}rgen Buntrock, and studied
in Section \ref{sec:tm62d}.
For this machine $N$, we have $s(N) = 8,690,333,381,690,951$ and
$\sigma(N)= 95,524,079$, that is, same value of $\sigma$,
and almost twice the value of $s$.
See analysis of this similarity in Section \ref{sec:sb}.

\bigskip

\noindent{\bf Analysis by Shawn Ligocki}:

Let $C(n, 0) = \ldots 0(A0)12^n0 \ldots$,\\
and $C(n, 1) = \ldots 0(C0)12^n0 \ldots$. 

Then we have, for all $k \ge 0$,
$$\begin{array}{ccc}
C(2k, 0)     & \vdash ( 40k^2 + 32k + 5 )  & C(5k + 1, 1)\\
C(2k + 1, 0) & \vdash ( 40k^2 + 82k + 42 ) & \ldots  01^{10k+9}(H0)0  \ldots\\ 
C(2k + 1, 1) & \vdash ( 40k^2 + 52k + 19 ) & C(5k + 3, 1)\\
C(2k + 2, 1) & \vdash ( 40k^2 + 92k + 53 ) & C(5k + 5, 0)
\end{array}$$

So we have:
$$\begin{array}{rl}
\multicolumn{2}{c}{\ldots 0(A0)0 \ldots =}\\
C( 0, 0 ) & \vdash ( 5 )\\
C( 1, 1 ) & \vdash ( 19 )\\
C( 3, 1 ) & \vdash ( 111 )\\
C( 8, 1 ) & \vdash ( 689 )\\
C( 20, 0 ) & \vdash ( 4,325 )\\
C( 51, 1 ) & \vdash ( 26,319 )\\
C( 128, 1 ) & \vdash ( 164,609 )\\
C( 320, 0) & \vdash ( 1,029,125 )\\
C( 801, 1 ) & \vdash ( 6,420,819 )\\
C( 2003, 1 ) & \vdash ( 40,132,111 )\\
C( 5008, 1 ) & \vdash ( 250,830,689 )\\
C( 12520, 0 ) & \vdash ( 1,567,704,325 )\\
C( 31301, 1 ) & \vdash ( 9,797,713,819 )\\
C( 78253, 1 ) & \vdash ( 61,235,789,611 )\\
C( 195633, 1 ) & \vdash ( 382,723,880,691 )\\
C( 489083, 1 ) & \vdash ( 2,392,024,743,391 )\\
C( 1222708, 1 ) & \vdash ( 14,950,155,868,889 )\\
C( 3056770, 0 ) & \vdash ( 93,438,477,237,325 )\\
C( 7641926, 1 ) & \vdash ( 582,990,375,746,317 )\\
C( 19104815, 0 ) & \vdash ( 3,649,939,963,043,376 )\\
\multicolumn{2}{c}{\ldots 01^{95524079}(H0)0 \ldots} 
\end{array}$$

\subsubsection{Lafitte and Papazian's machine found in April 2006}
\label{sec:tm33f}

This machine was the record holder in the Busy Beaver Competition
for machines with 3 states and 3 symbols, from April to August 2006.

\bigskip

\begin{tabular}{cc}
\begin{tabular}{c}
Lafitte and Papazian (2006)\\
$s(M) = 4,144,465,135,614$\\
$\sigma(M) = 2,950,149$
\end{tabular}
 & 
\begin{tabular}{c|ccc|}
  &  0  &  1  &  2\\
\hline
A & 1RB & 1RH &	2LC\\
B & 1LC & 2RB &	1LB\\
C & 1LA & 2RC &	2LA\\
\hline
\end{tabular}
\end{tabular}

\bigskip

Let $C(n, 0) = \ldots 0(A0)1^n0 \ldots$,\\
and $C(n, 1) = \ldots 0(A0)1^n210 \ldots$. 

Then we have, for all $k \ge 0$ (note the likeness to Brady's machine
of Section \ref{sec:tm33a}),
$$\begin{array}{ccc}
\ldots  0(A0)0 \ldots & \vdash ( 16 )      & C(6, 0)\\
C(2k + 1, 0) &     \vdash ( 4k + 5 )       & \ldots 01(H2)2^{2k}10 \ldots\\ 
C(2k + 2, 0) & \vdash ( 10k^2 + 27k + 23 ) & C(5k + 6, 1)\\
C(2k, 1)     & \vdash ( 10k^2 + 27k + 18 ) & C(5k + 5, 1)\\
C(2k + 1, 1) & \vdash ( 10k^2 + 51k + 60 ) & C(5k + 12, 0)
\end{array}$$

So we have:
$$\begin{array}{rl}
\ldots  0(A0)0 \ldots & \vdash ( 16 )\\
C( 6, 0 ) & \vdash ( 117 )\\
C( 16, 1 ) & \vdash ( 874 )\\
C( 45, 1 ) & \vdash ( 6,022 )\\
C( 122, 0 ) & \vdash ( 37,643 )\\
C( 306, 1 ) & \vdash ( 238,239 )\\
C( 770, 1 ) & \vdash ( 1,492,663 )\\
C( 1930, 1 ) & \vdash ( 9,338,323 )\\
C( 4830, 1 ) & \vdash ( 58,387,473 )\\
C( 12080, 1 ) & \vdash ( 364,979,098 )\\
C( 30205, 1 ) & \vdash ( 2,281,474,302 )\\
C( 75522, 0 ) & \vdash ( 14,259,195,543 )\\
C( 188806, 1 ) & \vdash ( 89,121,812,989 )\\
C( 472020, 1 ) & \vdash ( 557,013,573,288 )\\
C( 1180055, 1 ) & \vdash ( 3,481,348,698,727 )\\
C( 2950147, 0 ) & \vdash ( 5,900,297 )\\
\multicolumn{2}{c}{\ldots  01(H2)2^{2950146}10 \ldots} 
\end{array}$$

\bigskip

Note that we have also, for all $k \ge 0$,
$$\begin{array}{ccc}
\ldots  0(A0)0 \ldots &    \vdash ( 133 )     & C(16, 1)\\
C(2k, 1)     & \vdash ( 10k^2 + 27k + 18 )    & C(5k + 5, 1)\\
C(4k + 1, 1) & \vdash ( 290k^2 + 737k + 468 ) & C(25k + 31, 1)\\
C(4k + 3, 1) & \vdash ( 40k^2 + 162k + 158 )  & \ldots  01(H2)2^{10k+16}10 \ldots
\end{array}$$

\subsubsection{Lafitte and Papazian's machine found in September 2005}
\label{sec:tm33e}

This machine was the record holder in the Busy Beaver Competition for
machines with 3 states and 3 symbols, from September 2005 to April 2006.

\bigskip

\begin{tabular}{cc}
\begin{tabular}{c}
Lafitte and Papazian (2005)\\
$s(M) = 987,522,842,126$\\
$\sigma(M) = 1,525,688$
\end{tabular}
 & 
\begin{tabular}{c|ccc|}
  &  0  &  1  &  2\\
\hline
A & 1RB & 2LA &	1RA\\
B & 1RC & 2RB &	0RC\\
C & 1LA & 1RH &	1LA\\
\hline
\end{tabular}
\end{tabular}

\bigskip

Let $C(n, 0) = \ldots 0(A0)2^n0 \ldots$,\\
and $C(n, 1) = \ldots 0(A0)2^n10 \ldots$. 

Then we have, for all $k \ge 0$,
$$\begin{array}{ccc}
C(4k, 0)     & \vdash ( 14k^2 + 16k + 5 )  & 	C(7k + 2, 1)\\
C(4k + 1, 0) & \vdash ( 14k^2 + 30k + 15 ) & 	C(7k + 5, 0)\\
C(4k + 2, 0) & \vdash ( 14k^2 + 30k + 15 ) & 	C(7k + 5, 0)\\
C(4k + 3, 0) & \vdash ( 14k^2 + 44k + 35 ) & 	C(7k + 9, 1)\\
C(2k + 1, 1) & \vdash ( 4k + 3 )           & \ldots  01(12)^k01(H0)0  \ldots\\ 
C(4k, 1)     & \vdash ( 14k^2 + 26k + 11 ) & 	C(7k + 4, 0)\\
C(4k + 2, 1) & \vdash ( 14k^2 + 40k + 29 ) & 	C(7k + 8, 1)
\end{array}$$

So we have:
$$\begin{array}{rl}
\multicolumn{2}{c}{\ldots 0(A0)0 \ldots =}\\
C( 0, 0 ) & \vdash (5 )\\
C( 2, 1 ) & \vdash (29 )\\
C( 8, 1 ) & \vdash (119 )\\
C( 18, 0 ) & \vdash (359 )\\
C( 33, 0 ) & \vdash (1,151 )\\
C( 61, 0 ) & \vdash (3,615 )\\
C( 110, 0 ) & \vdash (11,031 )\\
C( 194, 0 ) & \vdash (33,711 )\\
C( 341, 0 ) & \vdash (103,715 )\\
C( 600, 0 ) & \vdash (317,405 )\\
C( 1052, 1 ) & \vdash (975,215 )\\
C( 1845, 0 ) & \vdash (2,989,139 )\\
C( 3232, 0 ) & \vdash (9,153,029 )\\
C( 5658, 1 ) & \vdash (28,048,133 )\\
C( 9906, 1 ) & \vdash (85,927,133 )\\
C( 17340, 1 ) & \vdash (263,203,871 )\\
C( 30349, 0 ) & \vdash (806,103,591 )\\
C( 53114, 0 ) & \vdash (2,468,672,331 )\\
C( 92951, 0 ) & \vdash (7,560,436,829 )\\
C( 162668, 1 ) & \vdash (23,154,325,799 )\\
C( 284673, 0 ) & \vdash (70,910,514,191 )\\
C( 498181, 0 ) & \vdash (217,164,134,715 )\\
C( 871820, 0 ) & \vdash (665,064,835,635 )\\
C( 1525687, 1 ) & \vdash (3,051,375 )\\
\multicolumn{2}{c}{\ldots  01(12)^{762843}01(H0)0 \ldots} 
\end{array}$$

\subsubsection{Lafitte and Papazian's machine found in August 2005}
\label{sec:tm33d}

This machine was the record holder in the Busy Beaver Competition for
machines with 3 states and 3 symbols, from August to September 2005.

\bigskip

\begin{tabular}{cc}
\begin{tabular}{c}
Lafitte and Papazian (2005)\\
$s(M) = 4,939,345,068$\\
$\sigma(M) = 107,900$
\end{tabular}
 & 
\begin{tabular}{c|ccc|}
  &  0  &  1  &  2\\
\hline
A & 1RB & 1RH &	2RB\\
B & 1LC & 0LB &	1RA\\
C & 1RA & 2LC &	1RC\\
\hline
\end{tabular}
\end{tabular}

\bigskip

Let $C(n, 0) = \ldots 0(C0)2^n0 \ldots$,\\
and $C(n, 1) = \ldots 0(C0)2^n10  \ldots$. 

Then we have, for all $k \ge 0$,
$$\begin{array}{ccc}
\ldots 0(A0)0 \ldots &    \vdash ( 3 )     & C(1, 1)\\
C(4k, 0)     & \vdash ( 14k^2 + 16k + 5 )  & C(7k + 2, 1)\\
C(4k + 1, 0) & \vdash ( 14k^2 + 22k + 7 )  & C(7k + 3, 0)\\
C(4k + 2, 0) & \vdash ( 14k^2 + 30k + 15 ) & C(7k + 5, 0)\\
C(4k + 3, 0) & \vdash ( 14k^2 + 36k + 23 ) & C(7k + 7, 1)\\
C(2k, 1)     & \vdash ( 2k + 2 )           & \ldots  01(21)^k1(H0)0  \ldots\\ 
C(4k + 1, 1) & \vdash ( 14k^2 + 20k + 9 )  & C(7k + 3, 1)\\
C(4k + 3, 1) & \vdash ( 14k^2 + 34k + 21 ) & C(7k + 6, 0)
\end{array}$$

So we have:
$$\begin{array}{rl}
 \ldots  0(A0)0 \ldots & \vdash ( 3 )\\
C( 1, 1 ) & \vdash ( 9 )\\
C( 3, 1 ) & \vdash ( 21 )\\
C( 6, 0 ) & \vdash ( 59 )\\
C( 12, 0 ) & \vdash ( 179 )\\
C( 23, 1 ) & \vdash ( 541 )\\
C( 41, 0 ) & \vdash ( 1,627 )\\
C( 73, 0 ) & \vdash ( 4,939 )\\
C( 129, 0 ) & \vdash ( 15,047 )\\
C( 227, 0 ) & \vdash ( 45,943 )\\
C( 399, 1 ) & \vdash ( 140,601 )\\
C( 699, 0 ) & \vdash ( 430,151 )\\
C( 1225, 1 ) & \vdash ( 1,317,033 )\\
C( 2145, 1 ) & \vdash ( 4,032,873 )\\
C( 3755, 1 ) & \vdash ( 12,349,729 )\\
C( 6572, 0 ) & \vdash ( 37,818,579 )\\
C( 11503, 1 ) & \vdash ( 115,816,521 )\\
C( 20131, 0 ) & \vdash ( 354,675,511 )\\
C( 35231, 1 ) & \vdash ( 1,086,184,945 )\\
C( 61655, 0 ) & \vdash ( 3,326,402,857 )\\
C( 107898, 1 ) & \vdash ( 107,900 )\\
\multicolumn{2}{c}{\ldots 01(21)^{53949}1(H0)0 \ldots} 
\end{array}$$

\subsubsection{Souris's machine for {\boldmath $S(3,3)$}}
\label{sec:tm33c}

This machine was the record holder in the Busy Beaver Competition
for $S(3,3)$, from July to August 2005.

\begin{tabular}{cc}
\begin{tabular}{c}
Souris (2005)\\
$s(M) = 544,884,219$\\
$\sigma(M) = 32,213$
\end{tabular}
 & 
\begin{tabular}{c|ccc|}
  &  0  &  1  &  2\\
\hline
A & 1RB & 1LB &	2LA\\
B & 1LA & 1RC &	1RH\\
C & 0LA & 2RC &	1LC\\
\hline
\end{tabular}
\end{tabular}

\bigskip

Let $C(n, 0) = \ldots 0(A0)1^n0 \ldots$,\\
and $C(n, 1) = \ldots 0(A0)1^n20  \ldots$. 

Then we have, for all $k \ge 0$,
$$\begin{array}{ccc}
\ldots  0(A0)0 \ldots & \vdash ( 4 ) & 	C(3, 0)\\
C(3k + 2, 0) & \vdash ( 21k^2 + 43k + 19 ) & \ldots  011(H2)2^{7k + 1}0  \ldots\\
C(3k + 3, 0) & \vdash ( 21k^2 + 43k + 24 ) & C(7k + 7, 0)\\
C(3k + 4, 0) & \vdash ( 21k^2 + 43k + 26 ) & C(7k + 7, 1)\\
C(3k + 1, 1) & \vdash ( 21k^2 + 61k + 35 ) & \ldots  011(H2)2^{7k + 3}0  \ldots\\ 
C(3k + 2, 1) & \vdash ( 21k^2 + 61k + 42 ) & C(7k + 9, 0)\\
C(3k + 3, 1) & \vdash ( 21k^2 + 61k + 46 ) & C(7k + 9, 1)
\end{array}$$

So we have:
$$\begin{array}{rl}
\ldots  0(A0)0 \ldots & \vdash ( 4 )\\
C( 3, 0 ) & \vdash ( 24 )\\
C( 7, 0 ) & \vdash ( 90 )\\
C( 14, 1 ) & \vdash ( 622 )\\
C( 37, 0 ) & \vdash ( 3,040 )\\
C( 84, 1 ) & \vdash ( 17,002 )\\
C( 198, 1 ) & \vdash ( 92,736 )\\
C( 464, 1 ) & \vdash ( 507,472 )\\
C( 1087, 0 ) & \vdash ( 2,752,290 )\\
C( 2534, 1 ) & \vdash ( 15,010,582 )\\
C( 5917, 0 ) & \vdash ( 81,666,440 )\\
C( 13804, 1 ) & \vdash ( 444,833,917 )\\
\multicolumn{2}{c}{\ldots  011(H2)2^{32210}0 \ldots} 
\end{array}$$

\subsubsection{Souris's machine for {\boldmath $\mathnormal{\Sigma}(3,3)$}}
\label{sec:tm33b}

This machine was the record holder in the Busy Beaver Competition
for $\mathnormal{\Sigma}(3,3)$, from July to August 2005.

\begin{tabular}{cc}
\begin{tabular}{c}
Souris (2005)\\
$s(M) = 310,341,163$\\
$\sigma(M) = 36,089$
\end{tabular}
 & 
\begin{tabular}{c|ccc|}
  &  0  & 1  & 	 2\\
\hline
A & 1RB & 2RA &	2RC\\
B & 1LC & 1RH &	1LA\\
C & 1RA & 2LB &	1LC\\
\hline
\end{tabular}
\end{tabular}

\bigskip

Let $C(n, 0) = \ldots 0(C0)1^n0 \ldots$,\\
and $C(n, 1) = \ldots 0(C0)1^n210 \ldots$. 

Then we have, for all $k \ge 0$,

$$\begin{array}{ccc}
\ldots  0(A0)0 \ldots &   \vdash ( 4 )    & C(1, 1)\\
C(2k + 2, 0) & \vdash ( 5k^2 + 32k + 17 ) & C(5k + 5, 0)\\
C(2k + 3, 0) & \vdash ( 5k^2 + 32k + 21 ) & C(5k + 4, 1)\\
C(2k + 1, 1) & \vdash ( 5k^2 + 32k + 15 ) & C(5k + 4, 0)\\
C(2k + 2, 1) & \vdash ( 5k^2 + 37k + 30 ) & \ldots 012^{5k + 5}1(H2)10 \ldots 
\end{array}$$

So we have:

$$\begin{array}{rl}
\ldots 0(A0)0 \ldots & \vdash ( 4 )\\
C( 1, 1 ) & \vdash ( 15 )\\
C( 4, 0 ) & \vdash ( 54 )\\
C( 10, 0 ) & \vdash ( 225 )\\
C( 25, 0 ) & \vdash ( 978 )\\
C( 59, 1 ) & \vdash ( 5,148 )\\
C( 149, 0 ) & \vdash ( 29,002 )\\
C( 369, 1 ) & \vdash ( 175,183 )\\
C( 924, 0 ) & \vdash ( 1,077,374 )\\
C( 2310, 0 ) & \vdash ( 6,695,525 )\\
C( 5775, 0 ) & \vdash ( 41,737,353 )\\
C( 14434, 1 ) & \vdash ( 260,620,302 )\\
\multicolumn{2}{c}{\ldots 012^{36085}1(H2)10 \ldots}
\end{array}$$

\bigskip

Note that we have also, for all $k \ge 0$,

$$\begin{array}{ccc}
\ldots 0(A0)0 \ldots & \vdash ( 19 ) & C(4, 0)\\
C(2k + 2, 0) & \vdash ( 5k^2 + 32k + 17 ) & C(5k + 5, 0)\\
C(4k + 3, 0) & \vdash ( 145k^2 + 299k + 93 ) & \ldots  012^{25k + 10}1(H2)10 \ldots\\ 
C(4k + 5, 0) & \vdash ( 145k^2 + 444k + 281 ) & C(25k + 24, 0)
\end{array}$$

\subsubsection{Brady's machine}\label{sec:tm33a}

This machine was the record holder in the Busy Beaver Competition for
machines with 3 states and 3 symbols, from December 2004 to July 2005.

\bigskip

\begin{tabular}{cc}
\begin{tabular}{c}
Brady (2004)\\
$s(M) = 92,649,163$\\
$\sigma(M) = 13,949$
\end{tabular}
 & 
\begin{tabular}{c|ccc|}
  &  0  &  1  &  2\\
\hline
A & 1RB & 1RH &	2LC\\
B & 1LC & 2RB &	1LB\\
C & 1LA & 0RB &	2LA\\
\hline
\end{tabular}
\end{tabular}

\bigskip

Let $C(n, 0) = \ldots 0(A0)1^n0  \ldots$,\\
and $C(n, 1) = \ldots 0(A0)1^n210  \ldots$.

Then we have, for all $k \ge 0$,
$$\begin{array}{ccc}
\ldots  0(A0)0 \ldots & \vdash ( 6 ) & C(0, 1)\\
C(2k + 1, 0) & \vdash ( 4k + 5 ) & \ldots  01(H2)2^{2k}10 \ldots\\ 
C(2k + 2, 0) & \vdash ( 10k^2 + 15k + 10 ) & C(5k + 3, 1)\\
C(2k, 1)     & \vdash ( 10k^2 + 27k + 18 ) & C(5k + 5, 1)\\
C(2k + 1, 1) & \vdash ( 10k^2 + 51k + 60 ) & C(5k + 12, 0)
\end{array}$$

So we have:
$$\begin{array}{rl}
\ldots 0(A0)0 \ldots & \vdash ( 6 )\\
C( 0, 1 ) & \vdash ( 18 )\\
C( 5, 1 ) & \vdash ( 202 )\\
C( 22, 0 ) & \vdash ( 1,160 )\\
C( 53, 1 ) & \vdash ( 8,146 )\\
C( 142, 0 ) & \vdash ( 50,060 )\\
C( 353, 1 ) & \vdash ( 318,796 )\\
C( 892, 0 ) & \vdash ( 1,986,935 )\\
C( 2228, 1 ) & \vdash ( 12,440,056 )\\
C( 5575, 1 ) & \vdash ( 77,815,887 )\\
C( 13947, 0 ) & \vdash ( 27,897 )\\
\multicolumn{2}{c}{\ldots 01(H2)2^{13946}10 \ldots} 
\end{array}$$

\bigskip

Note that we have also, for all $k \ge 0$,
$$\begin{array}{ccc}
\ldots  0(A0)0 \ldots & \vdash ( 6 ) & 	C(0, 1)\\
C(2k, 1)     & \vdash ( 10k^2 + 27k + 18 ) & C(5k + 5, 1)\\
C(4k + 1, 1) & \vdash ( 290k^2 + 677k + 395 ) & C(25k + 28, 1)\\
C(4k + 3, 1) & \vdash ( 40k^2 + 162k + 158 ) & \ldots 01(H2)2^{10k+16}10 \ldots 
\end{array}$$

\subsection{Turing machines with 2 states and 4 symbols}

\subsubsection{Ligockis' champion}\label{sec:tm24a}

This machine is the record holder in the Busy Beaver Competition
for machines with 2 states and 4 symbols, since February 2005.

\bigskip

\begin{tabular}{cc}
\begin{tabular}{c}
Terry and Shawn Ligocki (2005)\\
$s(M) = 3,932,964 =$? $S(2,4)$\\
$\sigma(M) = 2,050 =$? $\mathnormal{\Sigma}(2,4)$
\end{tabular}
 & 
\begin{tabular}{c|cccc|}
  &  0  &  1  &	 2  &  3\\
\hline
A & 1RB & 2LA &	1RA & 1RA\\
B & 1LB & 1LA &	3RB & 1RH\\
\hline
\end{tabular}
\end{tabular}

\bigskip

Let $C(n, 1) = \ldots 0(A0)2^n10 \ldots$,\\
and $C(n, 2) = \ldots 0(A0)2^n110 \ldots$.

Then we have, for all $k \ge 0$,
$$\begin{array}{ccc}
\ldots 0(A0)0 \ldots &    \vdash ( 6 )     & C(1, 2)\\
C(3k, 1)     & \vdash ( 15k^2 + 9k + 3 )   & C(5k + 1, 1)\\
C(3k + 1, 1) & \vdash ( 15k^2 + 24k + 13 ) & \ldots 013^{5k+2}1(H1)0 \ldots\\
C(3k + 2, 1) & \vdash ( 15k^2 + 29k + 17 ) & C(5k + 4, 2)\\
C(3k, 2)     & \vdash ( 15k^2 + 11k + 3 )  & C(5k + 1, 2)\\
C(3k + 1, 2) & \vdash ( 15k^2 + 21k + 7 )  & C(5k + 3, 1)\\
C(3k + 2, 2) & \vdash ( 15k^2 + 36k + 23 ) & \ldots 013^{5k+4}1(H1)0 \ldots
\end{array}$$

So we have:
$$\begin{array}{rl}
\ldots 0(A0)0 \ldots & \vdash ( 6 )\\
C( 1, 2 ) & \vdash ( 7 )\\
C( 3, 1 ) & \vdash ( 27 )\\
C( 6, 1 ) & \vdash ( 81 )\\
C( 11, 1 ) & \vdash ( 239 )\\
C( 19, 2 ) & \vdash ( 673 )\\
C( 33, 1 ) & \vdash ( 1,917 )\\
C( 56, 1 ) & \vdash ( 5,399 )\\
C( 94, 2 ) & \vdash ( 15,073 )\\
C( 158, 1 ) & \vdash ( 42,085 )\\
C( 264, 2 ) & \vdash ( 117,131 )\\
C( 441, 2 ) & \vdash ( 325,755 )\\
C( 736, 2 ) & \vdash ( 905,527 )\\
C( 1228, 1 ) & \vdash ( 2,519,044 )\\
\multicolumn{2}{c}{\ldots 013^{2047}1(H1)0 \ldots}
\end{array}$$

See detailed analysis in Michel (2015), Section 4.

\subsubsection{Brady's runner-up}\label{sec:tm24b}

This machine was the record holder in the Busy Beaver Competition
for machines with 2 states and 4 symbols, from 1988 to February 2005.

\bigskip

\begin{tabular}{cc}
\begin{tabular}{c}
Brady (1988)\\
$s(M) = 7,195$\\
$\sigma(M) = 90$
\end{tabular}
 & 
\begin{tabular}{c|cccc|}
  &  0  &  1  &	 2  &  3\\
\hline
A & 1RB & 3LA &	1LA & 1RA\\
B & 2LA & 1RH &	3RA & 3RB\\
\hline
\end{tabular}
\end{tabular}

\bigskip

Let $C(n, 0) = \ldots 0(A0)3^n0 \ldots$,\\
and $C(n, 1) = \ldots 0(A0)3^n20 \ldots$.

Then we have, for all $k \ge 0$,
$$\begin{array}{ccc}
C(3k, 0) 	 & \vdash ( 15k^2 + 7k + 3 )   &  C(5k + 1, 1)\\
C(3k + 1, 0) 	 & \vdash ( 15k^2 + 22k + 11 ) &  \ldots 013^{5k+1}1(H0)0 \ldots\\
C(3k + 2, 0) 	 & \vdash ( 15k^2 + 27k + 13 ) &  C(5k + 4, 0)\\
C(3k, 1) 	 & \vdash ( 15k^2 + 28k + 16 ) &  \ldots 013^{5k+3}1(H0)0 \ldots\\
C(3k + 1, 1) 	 & \vdash ( 15k^2 + 33k + 19 ) &  C(5k + 5, 0)\\
C(3k + 2, 1) 	 & \vdash ( 15k^2 + 43k + 33 ) &  C(5k + 7, 1)
\end{array}$$

So we have:
$$\begin{array}{rl}
\multicolumn{2}{c}{\ldots 0(A0)0 \ldots=}\\
C( 0, 0 ) & \vdash ( 3 )\\
C( 1, 1 ) & \vdash ( 19 )\\
C( 5, 0 ) & \vdash ( 55 )\\
C( 9, 0 ) & \vdash ( 159 )\\
C( 16, 1 ) & \vdash ( 559 )\\
C( 30, 0 ) & \vdash ( 1,573 )\\
C( 51, 1 ) & \vdash ( 4,827 )\\
\multicolumn{2}{c}{\ldots 013^{88}1(H0)0 \ldots}
\end{array}$$

\subsection{Turing machines with 2 states and 5 symbols}

\subsubsection{Ligockis' champion}\label{sec:tm25h}

This machine is the record holder in the Busy Beaver Competition
for machines with 2 states and 5 symbols, since November 2007.

\bigskip

\begin{tabular}{cc}
\begin{tabular}{c}
Terry and Shawn Ligocki (2007)\\
$s(M)$ and $S(2,5) > 1.9 \times 10^{704}$\\
$\sigma(M)$ and $\mathnormal{\Sigma}(2,5) > 1.7 \times 10^{352}$
\end{tabular}
 & 
\begin{tabular}{c|ccccc|} 
  &	 0  &	 1  & 	 2  & 	 3  &	 4\\
\hline
A &	1RB &	2LA &	1RA &	2LB &	2LA\\
B &	0LA &	2RB &	3RB &	4RA &	1RH\\
\hline
\end{tabular}
\end{tabular}

\bigskip

Let $C(n, 1) = \ldots 013^n(B0)0 \ldots$,\\
and $C(n, 2) = \ldots 023^n(B0)0 \ldots$,\\
and $C(n, 3) = \ldots 03^n(B0)0 \ldots$,\\
and $C(n, 4) = \ldots 04113^n(B0)0 \ldots$,\\
and $C(n, 5) = \ldots 04123^n(B0)0 \ldots$,\\
and $C(n, 6) = \ldots 0413^n(B0)0 \ldots$,\\
and $C(n, 7) = \ldots 0423^n(B0)0 \ldots$,\\
and $C(n, 8) = \ldots 043^n(B0)0 \ldots$.

Then we have, for all $k \ge 0$,
$$\begin{array}{ccc}
\ldots 0(A0)0\ldots & \vdash (1)          &  C(0,1)\\
C(2k, 1)     & \vdash ( 3k^2 + 8k + 4 )   & C(3k + 1, 1)\\
C(2k + 1, 1) & \vdash ( 3k^2 + 8k + 4 )   & C(3k + 1, 2)\\
C(2k, 2)     & \vdash ( 3k^2 + 14k + 9 )  & C(3k + 2, 1)\\
C(2k + 1, 2) & \vdash ( 3k^2 + 8k + 4 )   & C(3k + 2, 3)\\
C(2k, 3)     & \vdash ( 3k^2 + 8k + 2 )   & C(3k, 1)\\
C(2k + 1, 3) & \vdash ( 3k^2 + 8k + 22 )  & C(3k + 1, 4)\\
C(2k, 4)     & \vdash ( 3k^2 + 8k + 8 )   & C(3k + 3, 1)\\
C(2k + 1, 4) & \vdash ( 3k^2 + 8k + 4 )   & C(3k + 1, 5)\\
C(2k, 5)     & \vdash ( 3k^2 + 14k + 13 ) & C(3k + 4, 1)\\
C(2k + 1, 5) & \vdash ( 3k^2 + 8k + 4 )   & C(3k + 2, 6)\\
C(2k, 6)     & \vdash ( 3k^2 + 8k + 6 )   & C(3k + 2, 1)\\
C(2k + 1, 6) & \vdash ( 3k^2 + 8k + 4 )   & C(3k + 1, 7)\\
C(2k, 7)     & \vdash ( 3k^2 + 14k + 11 ) & C(3k + 3, 1)\\
C(2k + 1, 7) & \vdash ( 3k^2 + 8k + 4 )   & C(3k + 2, 8)\\
C(2k, 8)     & \vdash ( 3k^2 + 8k + 4 )   & C(3k + 1, 1)\\
C(2k + 1, 8) & \vdash ( 3k^2 + 5k + 3 )   & \ldots 01(H2)2^{3k}0 \ldots
\end{array}$$

So we have:
$$\begin{array}{rl}
\ldots 0(A0)0 \ldots & \vdash ( 1 )\\
C( 0, 1 ) & \vdash ( 4 )\\
C( 1, 1 ) & \vdash ( 4 )\\
C( 1, 2 ) & \vdash ( 4 )\\
C( 2, 3 ) & \vdash ( 13 )\\
C( 3, 1 ) & \vdash ( 15 )\\
C( 4, 2 ) & \vdash ( 49 )\\
C( 8, 1 ) & \vdash ( 84 )\\
C( 13, 1 ) & \vdash ( 160 )\\
C( 19, 2 ) & \vdash ( 319 )\\
C( 29, 3 ) & \vdash ( 722 )\\
C( 43, 4 ) & \vdash ( 1495 )\\
C( 64, 5 ) & \vdash ( 3533 )\\
C( 100, 1 ) & \vdash ( 7904 )\\
\multicolumn{2}{c}{\cdots}
\end{array}$$

See detailed analysis in Michel (2015), Section 5.

\subsubsection{Ligockis' machine found in August 2006}
\label{sec:tm25e}

This machine was the record holder in the Busy Beaver Competition
for machines with 2 states and 5 symbols, from August 2006 to October 2007.

\bigskip

\begin{tabular}{cc}
\begin{tabular}{c}
Terry and Shawn Ligocki (2006)\\
$s(M) = 7,069,449,877,176,007,352,687$\\
$\sigma(M) = 172,312,766,455$
\end{tabular}
 & 
\begin{tabular}{c|ccccc|} 	
  &  0  &  1  &	 2  &  3  &  4\\
\hline
A & 1RB & 0RB &	4RA & 2LB & 2LA\\
B & 2LA & 1LB &	3RB & 4RA & 1RH\\
\hline
\end{tabular}
\end{tabular}

\bigskip

\noindent{\bf Analysis by Shawn Ligocki}:

Let $C(n, 1) = \ldots 03^n(B0)0 \ldots$,\\
and $C(n, 2) = \ldots 013^n(B0)0 \ldots$,\\
and $C(n, 3) = \ldots 01403^n(B0)0 \ldots$,\\
and $C(n, 4) = \ldots 01413^n(B0)0 \ldots$.

Then we have, for all $k \ge 0$,
$$\begin{array}{ccc}
\ldots 0(A0)0 \ldots 	 & \vdash ( 1 )   &  C(0, 2)\\
C(2k, 1)     & \vdash ( 5k^2 + 14k + 3 )  &  C(5k + 1, 2)\\
C(2k + 1, 1) & \vdash ( 5k^2 + 14k + 7 )  &  C(5k + 3, 2)\\
C(2k, 2)     & \vdash ( 5k^2 + 14k + 3 )  &  C(5k + 1, 1)\\
C(2k + 1, 2) & \vdash ( 5k^2 + 14k + 11 ) &  C(5k + 2, 3)\\
C(2k, 3)     & \vdash ( 5k^2 + 14k + 3 )  &  C(5k + 1, 4)\\
C(2k + 1, 3) & \vdash ( 5k^2 + 14k + 9 )  &  C(5k + 4, 1)\\
C(2k, 4)     & \vdash ( 5k^2 + 14k + 3 )  &  C(5k + 1, 3)\\
C(2k + 1, 4) & \vdash ( 5k^2 + 9k + 4 )   & \ldots 011(H1)2^{5k+2}0 \ldots
\end{array}$$

So we have (in 30 transitions):
$$\begin{array}{rl}
\ldots 0(A0)0 \ldots & \vdash ( 1 )\\
C( 0, 2 ) & \vdash ( 3 )\\
C( 1, 1 ) & \vdash ( 7 )\\
C( 3, 2 ) & \vdash ( 30 )\\
C( 7, 3 ) & \vdash ( 96 )\\
C( 19, 1 ) & \vdash ( 538 )\\
\multicolumn{2}{c}{\cdots}\\
C( 4411206821, 1 ) & \vdash ( 24,323,432,041,896,588,247 )\\
C( 11028017053, 2 ) & \vdash ( 152,021,450,201,199,582,755 )\\
C( 27570042632, 3 ) & \vdash ( 950,134,063,605,862,157,707 )\\
C( 68925106581, 4 ) & \vdash ( 5,938,337,896,640,612,100,114 )\\
\multicolumn{2}{c}{\ldots 011(H1)2^{172312766452}0 \ldots}
\end{array}$$

\subsubsection{Lafitte and Papazian's machine found in June 2006}
\label{sec:tm25d}

This machine was the record holder in the Busy Beaver Competition
for $\mathnormal{\Sigma}(2,5)$, from June to August 2006.

\bigskip

\begin{tabular}{cc}
\begin{tabular}{c}
G. Lafitte and C. Papazian (2006)\\
$s(M) = 14,103,258,269,249$\\
$\sigma(M) = 4,848,239$
\end{tabular}
 & 
\begin{tabular}{c|ccccc|} 	
  &	 0  & 	 1  & 	 2  & 	 3  &	 4\\
\hline
A &	1RB &	3LB &	4LB &	4LA &	2RA\\
B &	2LA &	1RH &	3RB &	4RA &	3RB\\
\hline
\end{tabular}
\end{tabular}

\bigskip

Let $C(n, 1) = \ldots 0132^n33(B0)0 \ldots$,\\
and $C(n, 2) = \ldots 01342^n33(B0)0 \ldots$,\\
and $C(n, 3) = \ldots 0142^n33(B0)0 \ldots$,\\
and $C(n, 4) = \ldots 012^n33(B0)0 \ldots$.

Then we have, for all $k \ge 0$,
$$\begin{array}{ccc}
\ldots 0(A0)0 \ldots 	 & \vdash ( 10 )  &  C(0, 1)\\
C(2k, 1)     & \vdash ( 3k^2 + 12k + 15 ) &  C(3k + 2, 2)\\
C(2k + 1, 1) & \vdash ( 3k^2 + 12k + 11 ) &  C(3k + 2, 3)\\
C(2k, 2)     & \vdash ( 3k^2 + 12k + 9 )  &  C(3k + 2, 1)\\
C(2k + 1, 2) & \vdash ( 3k^2 + 18k + 30 ) &  C(3k + 5, 2)\\
C(2k, 3)     & \vdash ( 3k^2 + 12k + 9 )  &  C(3k + 2, 4)\\
C(2k + 1, 3) & \vdash ( 3k^2 + 18k + 28 ) &  C(3k + 4, 2)\\
C(2k, 4)     & \vdash ( 3k^2 + 12k + 13 ) &  C(3k + 1, 2)\\
C(2k + 1, 4) & \vdash ( 3k^2 + 9k + 5 )   &  \ldots 01(H4)4^{3k+2}20 \ldots
\end{array}$$

So we have (in 36 transitions):
$$\begin{array}{rl}
\ldots 0(A0)0 \ldots & \vdash ( 10 )\\
C( 0, 1 ) & \vdash ( 15 )\\
C( 2, 2 ) & \vdash ( 24 )\\
C( 5, 1 ) & \vdash ( 47 )\\
C( 8, 3 ) & \vdash ( 105 )\\
C( 14, 4 ) & \vdash ( 244 )\\
\multicolumn{2}{c}{\cdots}\\
C( 957674, 2 ) & \vdash ( 687,860,363,760 )\\
C( 1436513, 1 ) & \vdash ( 1,547,683,663,691 )\\
C( 2154770, 3 ) & \vdash ( 3,482,288,243,304 )\\
C( 3232157, 4 ) & \vdash ( 7,835,138,850,959 )\\
\multicolumn{2}{c}{\ldots 01(H4)4^{4848236}20 \ldots}
\end{array}$$

\subsubsection{Lafitte and Papazian's machine found in May 2006}
\label{sec:tm25c}

This machine was the record holder in the Busy Beaver Competition
for machines with 2 states and 5 symbols, from May to June 2006.

\bigskip

\begin{tabular}{cc}
\begin{tabular}{c}
G. Lafitte and C. Papazian (2006)\\
$s(M) = 3,793,261,759,791$\\
$\sigma(M) = 2,576,467$
\end{tabular}
 & 
\begin{tabular}{c|ccccc|} 	
  &	 0  & 	 1  & 	 2  & 	 3  &	4\\
\hline
A &	1RB &	3RA &	4LB &	2RA &	3LA\\
B &	2LA &	1RH &	4RB &	4RB &	2LB\\
\hline
\end{tabular}
\end{tabular}

\bigskip

Let $C(n, 1) = \ldots 014^n(B0)0 \ldots$,\\
and $C(n, 2) = \ldots 034^n(B0)0 \ldots$.

Then we have, for all $k \ge 0$,
$$\begin{array}{ccc}
\ldots 0(A0)0 \ldots 	 & \vdash ( 1 )   &  C(0, 1)\\
C(3k, 1)     & \vdash ( 4k^2 + 17k + 11 ) &  C(4k + 3, 1)\\
C(3k + 1, 1) & \vdash ( 4k^2 + 25k + 20 ) &  C(4k + 4, 1)\\
C(3k + 2, 1) & \vdash ( 4k^2 + 17k + 13 ) &  C(4k + 3, 2)\\
C(3k, 2)     & \vdash ( 4k^2 + 17k + 11 ) &  C(4k + 3, 1)\\
C(3k + 1, 2) & \vdash ( 4k^2 + 25k + 20 ) &  C(4k + 4, 1)\\
C(3k + 2, 2) & \vdash ( 4k^2 + 21k + 24 ) &  \ldots 01(H2)23^{4k+3}20 \ldots
\end{array}$$

So we have (in 45 transitions):
$$\begin{array}{rl}
\ldots 0(A0)0 \ldots & \vdash ( 1 )\\
C( 0, 1 ) & \vdash ( 11 )\\
C( 3, 1 ) & \vdash ( 32 )\\
C( 7, 1 ) & \vdash ( 86 )\\
C( 12, 1 ) & \vdash ( 143 )\\
C( 19, 1 ) & \vdash ( 314 )\\
\multicolumn{2}{c}{\cdots}\\
C( 815207, 1 ) & \vdash ( 295,364,260,408 )\\
C( 1086943, 2 ) & \vdash ( 525,094,796,254 )\\
C( 1449260, 1 ) & \vdash ( 933,496,546,059 )\\
C( 1932347, 2 ) & \vdash ( 1,659,550,059,339 )\\
\multicolumn{2}{c}{\ldots 01(H2)23^{2576463}20 \ldots}
\end{array}$$

Note that we have also, for all $k \ge 0$,
$$\begin{array}{ccc}
\ldots 0(A0)0 \ldots 	 & \vdash ( 1 )       & C(0, 1)\\
C(3k, 1)     & \vdash ( 4k^2 + 17k + 11 )     & C(4k + 3, 1)\\
C(3k + 1, 1) & \vdash ( 4k^2 + 25k + 20 )     & C(4k + 4, 1)\\
C(9k + 2, 1) & \vdash ( 100k^2 + 151k + 45 )  & C(16k + 7, 1)\\
C(9k + 5, 1) & \vdash ( 100k^2 + 239k + 120 ) & C(16k + 12, 1)\\
C(9k + 8, 1) & \vdash ( 100k^2 + 279k + 186 ) & \ldots 01(H2)23^{16k+15}20 \ldots
\end{array}$$

Note: The machine obtained by replacing B4 $\rightarrow$ 2LB by B4 $\rightarrow$ 3LB has
the same behavior but final configuration $\ldots 01(H3)3^{2576464}20 \ldots$.

\subsubsection{Lafitte and Papazian's machine found in December 2005}
\label{sec:tm25b}

This machine was the record holder in the Busy Beaver Competition
for $S(2,5)$, from December 2005 to May 2006.

\bigskip

\begin{tabular}{cc}
\begin{tabular}{c}
G. Lafitte and C. Papazian (2005)\\
$s(M) = 924,180,005,181$\\
$\sigma(M) = 1,137,477$
\end{tabular}
 & 
\begin{tabular}{c|ccccc|} 	
  &  0  &  1  &  2  &  3  &  4\\
\hline
A & 1RB & 3RA &	1LA & 1LB & 3LB\\
B & 2LA & 4LB &	3RA & 2RB & 1RH\\
\hline
\end{tabular}
\end{tabular}

\bigskip

Let $C(n, 1) = \ldots 012^n(B0)0 \ldots$,\\
and $C(n, 2) = \ldots 032^n(B0)0 \ldots$.

Then we have, for all $k \ge 0$,
$$\begin{array}{ccc}
\ldots 0(A0)0 \ldots 	 & \vdash ( 69 )   & C(8, 1)\\
C(2k + 1, 1) & \vdash ( 15k^2 + 37k + 31 ) & \ldots 01221(H1)1^{5k+1}20 \ldots\\
C(2k + 2, 1) & \vdash ( 15k^2 + 32k + 19 ) & C(5k + 3, 2)\\
C(2k, 2)     & \vdash ( 15k^2 + 32k + 19 ) & C(5k + 3, 1)\\
C(2k + 1, 2) & \vdash ( 15k^2 + 62k + 70 ) & C(5k + 9, 1)
\end{array}$$

So we have:
$$\begin{array}{rl}
\ldots 0(A0)0 \ldots & \vdash ( 69 )\\
C( 8, 1 ) & \vdash ( 250 )\\
C( 18, 2 ) & \vdash ( 1,522 )\\
C( 48, 1 ) & \vdash ( 8,690 )\\
C( 118, 2 ) & \vdash ( 54,122 )\\
C( 298, 1 ) & \vdash ( 333,315 )\\
C( 743, 2 ) & \vdash ( 2,087,687 )\\
C( 1864, 1 ) & \vdash ( 13,031,226 )\\
C( 4658, 2 ) & \vdash ( 81,438,162 )\\
C( 11648, 1 ) & \vdash ( 508,796,290 )\\
C( 29118, 2 ) & \vdash ( 3,179,933,122 )\\
C( 72798, 1 ) & \vdash ( 19,873,380,815 )\\
C( 181993, 2 ) & \vdash ( 124,209,722,062 )\\
C( 454989, 1 ) & \vdash ( 776,311,217,849 )\\
\multicolumn{2}{c}{\ldots 01221(H1)1^{1137471}20 \ldots}
\end{array}$$

Note that we have also, for all $k \ge 0$,
$$\begin{array}{ccc}
\ldots 0(A0)0 \ldots 	 & \vdash ( 69 )      & C(8, 1)\\
C(2k + 1, 1) & \vdash ( 15k^2 + 37k + 31 )    & \ldots 01221(H1)1^{5k+1}20 \ldots\\
C(4k + 2, 1) & \vdash ( 435k^2 + 524k + 166 ) & C(25k + 14, 1)\\
C(4k + 4, 1) & \vdash ( 435k^2 + 884k + 453 ) & C(25k + 23, 1)
\end{array}$$

\subsubsection{Lafitte and Papazian's machine found in October 2005}
\label{sec:tm25a}

This machine was the record holder in the Busy Beaver Competition
for $\mathnormal{\Sigma}(2,5)$, from October 2005 to May 2006.

\bigskip

\begin{tabular}{cc}
\begin{tabular}{c}
G. Lafitte and C. Papazian (2005)\\
$s(M) = 912,594,733,606$\\
$\sigma(M) = 1,957,771$
\end{tabular}
 & 
\begin{tabular}{c|ccccc|}
  &  0  &  1  &  2  &  3  &  4\\
\hline
A & 1RB & 3LB &	1RH & 1LA & 1LA\\
B & 2LA & 3RB &	4LB & 4LB & 3RA\\
\hline
\end{tabular}
\end{tabular}

\bigskip

Let $C(n, 1) = \ldots 0(A0)1^n20 \ldots$,\\
and $C(n, 2) = \ldots 0(A0)1^n40 \ldots$,\\
and $C(n, 3) = \ldots 0(A0)1^n320 \ldots$.

Then we have, for all $k \ge 0$,
$$\begin{array}{ccc}
\ldots 0(A0)0 \ldots 	 & \vdash ( 11 )  & C(3, 1)\\
C(2k + 1, 1) & \vdash ( 5k^2 + 28k + 26 ) & C(5k + 6, 1)\\
C(2k + 2, 1) & \vdash ( 5k^2 + 18k + 11 ) & C(5k + 3, 2)\\
C(2k, 2)     & \vdash ( 5k^2 + 18k + 11 ) & C(5k + 3, 1)\\
C(2k + 1, 2) & \vdash ( 5k^2 + 18k + 13 ) & C(5k + 3, 3)\\
C(2k + 1, 3) & \vdash ( 5k^2 + 18k + 9 )  & C(5k + 3, 1)\\
C(2k + 2, 3) & \vdash ( 5k^2 + 23k + 17 ) & \ldots 013^{5k+4}1(H0)0 \ldots
\end{array}$$

So we have:
$$\begin{array}{rl}
\ldots 0(A0)0 \ldots & \vdash ( 11 )\\
C( 3, 1 ) & \vdash ( 59 )\\
C( 11, 1 ) & \vdash ( 291 )\\
C( 31, 1 ) & \vdash ( 1,571 )\\
C( 81, 1 ) & \vdash ( 9,146 )\\
C( 206, 1 ) & \vdash ( 53,867 )\\
C( 513, 2 ) & \vdash ( 332,301 )\\
C( 1283, 3 ) & \vdash ( 2,065,952 )\\
C( 3208, 1 ) & \vdash ( 12,876,910 )\\
C( 8018, 2 ) & \vdash ( 80,432,578 )\\
C( 20048, 1 ) & \vdash ( 502,483,070 )\\
C( 50118, 2 ) & \vdash ( 3,140,218,478 )\\
C( 125298, 1 ) & \vdash ( 19,624,987,195 )\\
C( 313243, 2 ) & \vdash ( 122,653,507,396 )\\
C( 783108, 3 ) & \vdash ( 766,577,764,781 )\\
\multicolumn{2}{c}{\ldots 013^{1957769}1(H0)0 \ldots}
\end{array}$$

Note that we have also, for all $k \ge 0$,
$$\begin{array}{ccc}
\ldots 0(A0)0 \ldots  &   \vdash ( 11 )       & C(3, 1)\\
C(2k + 1, 1) & \vdash ( 5k^2 + 28k + 26 )     & C(5k + 6, 1)\\
C(2k + 2, 1) & \vdash ( 5k^2 + 18k + 11 )     & C(5k + 3, 2)\\
C(2k, 2)     & \vdash ( 5k^2 + 18k + 11 )     & C(5k + 3, 1)\\
C(4k + 1, 2) & \vdash ( 145k^2 + 176k + 45 )  & C(25k + 8, 1)\\
C(4k + 3, 2) & \vdash ( 145k^2 + 321k + 167 ) & \ldots 013^{25k+19}1(H0)0 \ldots 
\end{array}$$

\subsection{Collatz-like problems}\label{sec:clp}

Sameness of behaviors of the Turing machines above is striking.
Their behaviors depend on transitions in the following form:
$$C(ak + b) \vdash (\ )\ C(ck + d),$$
where $a$, $c$ are fixed, and $b = 0, \ldots , a-1$.
Sometimes, another parameter is added: $C(ak + b, p)$.

 These transitions can be compared to the following problem.
 Let $T$ be defined by
$$T(x) = \left\{\begin{array}{ll}
x/2        & \mbox{ if } x \mbox{ is even},\\
(3x + 1)/2 & \mbox{ if } x \mbox{ is odd}.
\end{array}\right.$$
 This can also be written
$$\begin{array}{rcl}
T(2m)     & = & m\\
T(2m + 1) & = & 3m + 2
\end{array}$$
 When $T$ is iterated over positive integers, do we always
 reach the loop: $T(2) = 1$, $T(1) = 2$?
 This question is a famous open problem in mathematics, called
 3$x$ + 1 problem, or \emph{Collatz problem}.

 A similar question can be asked about iterating transitions of
 configurations $C(ak + b, p)$ on positive integers.
 Do the iterated transitions always reach a halting configuration?
 For all the machines above
 (except for the machine with 6 states and 2 symbols in Section
 \ref{sec:tm62b}),
 this question is presently an open problem in mathematics.
 Because of likeness to Collatz problem, these problems are
 called \emph{Collatz-like problems}.
 Thus, for each machine above
 (except for the machine with 6 states and 2 symbols in Section
 \ref{sec:tm62b}),
 the halting problem (that is, on what inputs does this machine stop?)
 depends on an open Collatz-like problem.

\subsection{Non-Collatz-like behaviors}\label{sec:ncl}

Some Turing machines run a large number of steps on a small piece of tape.
Such machines do not seem to be Collatz-like.
We list below some interesting machines with this sort of behavior.

\subsubsection{Turing machines with 3 states and 3 symbols}

\begin{tabular}{cc}
\begin{tabular}{c}
A. H. Brady (November 2004)\\
$s(M) = 2,315,619$\\
$\sigma(M) = 31$
\end{tabular}
 & 
\begin{tabular}{c|ccc|}
  &  0 	&  1  &	 2\\
\hline
A & 1RB & 2LB &	1LC\\
B & 1LA & 2RB &	1RB\\
C & 1RH & 2LA &	0LC\\
\hline
\end{tabular}
\end{tabular}

\bigskip

Brady called this machine ``Surprise-in-a-Box''.

\bigskip

See also the simulation by Heiner Marxen:

\verb+http://turbotm.de/~heiner/BB/simAB3Y_SB.html+

\subsubsection{Turing machines with 2 states and 5 symbols}

\noindent{\bf (a) First machine}

\bigskip

\begin{tabular}{cc}
\begin{tabular}{c}
G. Lafitte and C. Papazian (July 2006)\\
$s(M) = 26,375,397,569,930$\\
$\sigma(M) = 143$
\end{tabular}
 & 
\begin{tabular}{c|ccccc|}
  &   0  &  1  &  2  &	 3  &  4\\
\hline
A &  1RB & 3LA & 1LA &	4LA & 1RA\\
B &  2LB & 2RA & 1RH &	0RA & 0RB\\
\hline
\end{tabular}
\end{tabular}

\bigskip

This machine was the record holder for $S(2,5)$, from July to August 2006.

\bigskip

See also the simulation by Heiner Marxen:

\verb+http://turbotm.de/~heiner/BB/simLaf25_j.html+

\bigskip

\noindent{\bf (b) Second machine}

\bigskip

\begin{tabular}{cc}
\begin{tabular}{c}
G. Lafitte and C. Papazian (July 2006)\\
$s(M) = 7,021,292,621$\\
$\sigma(M) = 37$
\end{tabular}
 & 
\begin{tabular}{c|ccccc|}
  &  0 	&  1  &	 2  &  3  &  4\\
\hline
A & 1RB & 4LA &	1LA & 1RH & 2RB\\
B & 2LB & 3LA &	1LB & 2RA & 0RB\\
\hline
\end{tabular}
\end{tabular}

\subsection{Turing machines in distinct classes with similar behaviors}
\label{sec:sb}

In this section, we give examples of machines that have similar behaviors,
but not the same numbers of states and symbols.

\subsubsection{(2,4)-TM and (3,3)-TM}

\begin{tabular}{cc}
\begin{tabular}{c}
Terry and Shawn Ligocki (2005)\\
$s(M) = 3,932,964 =$? $S(2,4)$\\
$\sigma(M) = 2,050 =$? $\mathnormal{\Sigma}(2,4)$
\end{tabular}
 & 
\begin{tabular}{c|cccc|}
  &  0 	&  1  &	 2  &  3\\
\hline
A & 1RB & 2LA &	1RA & 1RA\\
B & 1LB & 1LA &	3RB & 1RH\\
\hline
\end{tabular}
\end{tabular}

\bigskip

This machine is the record holder in the Busy Beaver Competition
for machines with 2 states and 4 symbols, since February 2005.

\bigskip

\begin{tabular}{cc}
\begin{tabular}{c}
A. H. Brady (2004)\\
$s(M) = 3,932,964$\\
$\sigma(M) = 2,050$
\end{tabular}
 & 
\begin{tabular}{c|ccc|}
  &  0  &  1  &	 2\\
\hline
A & 1RB & 1LC &	1RH\\
B & 1LA & 1LC &	2RB\\
C & 1RB & 2LC &	1RC\\
\hline
\end{tabular}
\end{tabular}

\bigskip

There is a step-by-step correspondence between the
configurations of these machines.

\subsubsection{(6,2)-TM and (3,3)-TM}

\begin{tabular}{cc}
\begin{tabular}{c}
Marxen and Buntrock (1997)\\
$s(M) = 8,690,333,381,690,951$\\
$\sigma(M) = 95,524,079$
\end{tabular}
 & 
\begin{tabular}{c|cc|}
  &  0  &  1\\
\hline
A & 1RB & 1RA\\
B & 1LC & 1LB\\
C & 0RF & 1LD\\
D & 1RA & 0LE\\
E & 1RH & 1LF\\
F & 0LA & 0LC\\
\hline
\end{tabular}
\end{tabular}

\bigskip

This machine was discovered in January 1990, and was published
on the web (Google groups) on September 3, 1997.
It was the record holder in the Busy Beaver Competition
for machines with 6 states and 2 symbols up to July 2000.

\bigskip

\begin{tabular}{cc}
\begin{tabular}{c}
Terry and Shawn Ligocki (2006)\\
$s(M) = 4,345,166,620,336,565$\\
$\sigma(M) = 95,524,079$
\end{tabular}
 & 
\begin{tabular}{c|ccc|}
  &  0 	&  1  &  2\\
\hline
A & 1RB & 2RC & 1LA\\
B & 2LA & 1RB &	1RH\\
C & 2RB & 2RA &	1LC\\
\hline
\end{tabular}
\end{tabular}

\bigskip

This machine was the record holder in the Busy Beaver Competition
for machines with 3 states and 3 symbols, from August 2006 to November 2007.

\bigskip

Note that these machines have same $\sigma$ value, and the $s$ value
of the first one is almost twice the $s$ value of the second one.

\bigskip

The behaviors of these machines can be related as follows.

Given the analyses of the (6,2)-TM in Section \ref{sec:tm62d}
and the (3,3)-TM in Section \ref{sec:tm33g},
the following functions $f$ and $g$ can be defined:

$$\begin{array}{cc}
\left\{ \begin{array}{rcl}
f(4k)     &   & \mbox{undefined,}\\
f(4k + 1) & = & 10k + 9,\\
f(4k + 2) & = & 10k + 9,\\
f(4k + 3) & = & 10k + 12.
\end{array}\right.
  &
\left\{ \begin{array}{rcl}
g(2k, 0)     & = & (5k + 1, 1),\\
g(2k + 1, 0) &   & \mbox{undefined,}\\
g(2k, 1)     & = & (5k, 0),\\
g(2k + 1, 1) & = & (5k + 3, 1).
\end{array}\right.
\end{array}$$

Now, let $h$ be defined by
$$\begin{array}{c}
h(n, 0) = 10n + 2,\\
h(n, 1) = 10n - 1. 
\end{array}$$

Then: $h \circ g = f \circ h$.

There is no step-by-step correspondence between these machines,
but there is a phase correspondence, according to functions $f$ and $g$.

\section{Properties of the busy beaver functions}

\subsection{Growth properties}

\begin{itemize}
\item
  Rado (1962) defined $S(n)$ and $\mathnormal{\Sigma}(n)$,
  that are denoted in this article
  $S(n,2)$ and $\mathnormal{\Sigma}(n,2)$.

\item
  These functions grow faster than any computable function.
  Formally, for any computable function $f$, there is an
  integer $N$ such that, for any integer $n > N$,

  $$S(n) > \mathnormal{\Sigma}(n) > f(n)$$

  This was proved by Rado (1962) who defined these functions
  in order to get noncomputable functions.

\item
  It is easy to prove that the two variables functions
  $S(n,k)$ and $\mathnormal{\Sigma}(n,k)$ are increasing with
  the number $n$ of states if the number $k$ of symbols is constant.
  Formally, for any integer $k \ge 2$, if $n > m$, then

  $$S(n,k) > S(m,k)\qquad \mbox{and}\qquad \mathnormal{\Sigma}(n,k) > \mathnormal{\Sigma}(m,k)$$

\item
  As Harland (2022) noticed, the same result for the
  number of symbols, with a constant number of states,
  is far from obvious, and still unproven. Petersen (2017)
  proved that functions $S(n,k)$ and $\mathnormal{\Sigma}(n,k)$
  are increasing with the number $k$ of symbols if the
  number $n$ of states is sufficiently large. The proof uses
  \emph{introspective encoding}, a tool developped by
  Ben-Amram and Petersen (2002).
\end{itemize}

\subsection{Relations between the busy beaver functions}
 
\begin{itemize}
 \item
 Rado (1962) proved that
 
$$S(n) < (n+1)\mathnormal{\Sigma}(5n) 2^{\mathnormal{\Sigma}(5n)}.$$

\item
  Julstrom (1993) proved that

  $$S(n) < \mathnormal{\Sigma}(28n).$$
 
 \item
 Julstrom (1992) proved that

$$S(n) < \mathnormal{\Sigma}(20n).$$
 
 \item
 Wang and Xu (1995) proved that

$$S(n) < \mathnormal{\Sigma}(10n).$$

\item
  In an unpublished technical report in German, Buro (1990) (p.\ 5-6)
  proved that

$$S(n) < \mathnormal{\Sigma}(9n).$$
 
 \item
 Yang, Ding and Xu (1997) proved that

$$S(n) < \mathnormal{\Sigma}(8n),$$
 
 and that there is a constant $c$ such that

$$S(n) < \mathnormal{\Sigma}(3n+c).$$

 \item
 Ben-Amram, Julstrom and Zwick (1996) proved that

$$S(n) < \mathnormal{\Sigma}(3n+6),$$

 and

$$S(n) < (2n-1)\mathnormal{\Sigma}(3n+3).$$

 \item
 Ben-Amram and Petersen (2002) proved that there is a constant $c$ such that

$$S(n) < \mathnormal{\Sigma}(n +8n/\log_2 n + c).$$

\end{itemize}

\bigskip

\section{Variants of busy beavers}

\subsection{Busy beavers defined by 4-tuples}

The Turing machines used for regular busy beavers are based on 5-tuples.
For example, the initial transition is
\begin{quote}
 (A,0) $\longrightarrow$ (1,R,B)
\end{quote}

\noindent and generally a transition is
\begin{quote}
 (state, scanned symbol) $\longrightarrow$ (new written symbol, move of the head, new state)
\end{quote}

Instead of both writing a symbol and moving the head in one transition,
these actions can be split up into two transitions, in the form of a 4-tuple:
\begin{quote}
 (state, scanned symbol) $\longrightarrow$ (new written symbol {\bf or} move of the head, new state)
\end{quote}

This alternative definition was introduced by Post in 1947
(Recursive unsolvability of a problem of Thue,
{\em The Journal of Symbolic Logic}, Vol. 12, 1-11). So Turing machines
defined by 4-tuples are also called {\em Post machines}, or
{\em Post-Turing machines}.

A busy beaver competition for such machines was studied by Oberschelp,
Schmidt-G\"{o}ttsch and Todt (1988), who defined two busy beaver functions,
for the number of non-blank symbols, and for the number of steps,
and gave some values and lower bounds for these functions.

The busy beaver competition for such machines are also studied by
P.\ Machado and F.\ Pereira, see

\verb+http://fmachado.dei.uc.pt/publications+

\noindent and B. van Heuveln and his team, see

\verb+http://www.cogsci.rpi.edu/~heuveb/Research/BB/index.html+

\smallskip

\noindent In their book, Boolos and Jeffrey (1974) used the
4-tuples variant to display the busy beaver problem.

\smallskip

\noindent Harland (2022) tackled 4-tuples machines and
he gave a proof of the following theorem:

\noindent {\bf Theorem}. For any $n$-state, $m$-symbol, 4-tuples machine $M$, halting on a blank tape,
there exists a $n$-state, $m$-symbol, 5-tuples machine $N$, halting on a blank tape,
such that $\sigma(N) = \sigma(M)$, that is, with the same number of
non-blank symbols written on the tape when it halts.

Moreover, the proof provides a simple algorithm that transforms a 4-tuples
machine into an equivalent 5-tuples machine. So Harland concludes that
searching for 5-tuples machines subsumes searching for 4-tuples machines.

\subsection{Busy beavers whose head can stand still}

In the definition of the Turing machines used for regular busy beavers,
the tape head has to move one cell right or left at each step, and cannot
stand still.
If we allow the tape head to stand still, new machines come into the competition,
and they can beat the current champions.

So Norbert B\'atfai found,
in August 2009, a Turing machine M with 5 states and 2 symbols with
$s(M)$ = 70,740,810 and $\sigma(M)$ = 4098. See

\verb+http://arxiv.org/abs/0908.4013+

This machine beats the current champion
for the number of steps ($s$ = 47,176,870).
It seems that relaxing this condition on moves does not allow us to obtain
machines with behaviors different from those of regular busy beavers.
But the study is still to be done.

\subsection{Busy beavers on a one-way infinite tape}

In the definition of the Turing machines used for regular busy beavers,
the tape is infinite on both left and right sides. Walsh (1982)
considered Turing machines with one-way infinite tape. Initially,
the tape head scans the first (leftmost) tape cell. A Turing machine halts
either by entering a halting state or by falling off the left end of the tape,
that is, moving left from cell 1. If a Turing machine $M$ halts
when it starts from a blank tape, its score is defined to be $k$ if the
rightmost tape cell ever visited by $M$'s head is the $k$th cell from
the left. $\mathnormal{\Sigma}(n,m)$ is defined as the largest score of all halting
$n$-state, $m$-symbol Turing machines. Walsh proved that, with this definition,
$\mathnormal{\Sigma}(2,3) = 6$.

\subsection{Two-dimensional busy beavers}

The Turing machines used for regular busy beavers have a one-dimensional tape.
Turing machines with two-dimensional or higher-dimensional tapes were first defined by
Hartmanis and Stearns in 1965 (On the computational complexity of algorithms,
{\em Transactions of the AMS}, Vol. 117, 285-306).

Brady (1988) launched the busy beaver competition for two-dimensional Turing
machines. He also defined, first, ``TurNing machines'', where
the head reorients itself at each step, and, second, machines that work
on a triangular grid.  

Tim Hutton resumed the search for two-dimensional busy beavers. See

\noindent \verb+https://github.com/GollyGang/ruletablerepository/wiki/TwoDimensionalTuringMachines+

\noindent He gave the following results:

\bigskip

For $S_2(k,n)$: ($k$ states, $n$ symbols)

\bigskip

\begin{tabular}{|r|c|c|c|c|}
\hline
3 symbols &   38     &    ?     &          &\\
\hline
2 symbols &    6     &   32     &  4632 ?  & 25,772,988,638 ?\\
\hline
          & 2 states & 3 states & 4 states & 5 states\\
\hline
\end{tabular}

\bigskip

For $\mathnormal{\Sigma}_2(k,n)$: ($k$ states, $n$ symbols)

\bigskip

\begin{tabular}{|r|c|c|c|c|}
\hline
3 symbols &   10     &    ?     &          &\\
\hline
2 symbols &    4     &   11     &  244 ?   & 935,508,401 ?\\
\hline
          & 2 states & 3 states & 4 states & 5 states\\ 
\hline
\end{tabular}

\bigskip

Note that
$$S_2(3,2) = 32 > S(3,2) = 21,$$
and
$$\mathnormal{\Sigma}_2(3,2) = 11 > \mathnormal{\Sigma}(3,2) = 6.$$

\bigskip

Tim Hutton also studied  higher-dimensional machines and found that,
for all $n > 0$, $S_n(2,2) = 6$ and $\mathnormal{\Sigma}_n(2,2) = 4$.

He also studied one-dimensional and higher-dimensional Turing machines
with {\em relative movements}, that is, where the head has an orientation
and reorients itself at each step.

\subsection{Beeping busy beavers}

In a survey article about busy beavers, Scott Aaronson (2020) defined
the beeping busy beaver function. While the classical busy beaver functions
are as complex as the halting problem for Turing machines, this new function
is as complex as the halting problem for Turing machines with an oracle
for the halting problem for Turing machines.

Let $M$ be a Turing machine with $k$ states and $n$ symbols, without
halting state, and let $q$ be a state of machine $M$. Machine $M$ is launched
on a blank tape, and beeps when it is in state $q$.
Let $s(M,q)$ be the last time that machine $M$ beeps on state $q$
(and $s(M,q)$ is infinite if $M$ beeps on state $q$ infinitely often).
Then the {\em beeping busy beaver function} is defined by

\bigskip

$BBB(k,n)$ = 1 + max\{$s(M,q)$ : $M$ is a Turing machine with $k$ states
and $n$ symbols, $q$ is a state of $M$ and $s(M,q)$ is finite\}

\bigskip

(The ``1 +'' is added for technical reasons).

\bigskip

The following results are known:
\begin{itemize}
\item $BBB(1,2) = 1$ (Scott Aaronson)
\item $BBB(2,2) = 6$ (Scott Aaronson)
\item $BBB(3,2) \ge 55$ (Scott Aaronson)
\item $BBB(4,2) \ge 32,779,478$ (Nicholas Drozd, July 2021).
  See analysis by Shawn Ligocki in
  
  \verb+https://www.sligocki.com/2021/07/17/bb-collatz.html+

\item $BBB(5,2) > 2.1 \times 10^{18}$ (Nicholas Drozd, January 2022)\\
      $BBB(5,2) > 1.7 \times 10^{502}$ (Shawn Ligocki, February 2022).
      See analysis by Shawn Ligocki in

  \verb+https://www.sligocki.com/2022/02/22/collatz-complex.html+\\

      $BBB(5,2) > 7.4 \times 10^{4079}$ (Nicholas Drozd, March 2022)\\
      $BBB(5,2) > 10^{14006}$ (Shawn Ligocki, April 2022)\\
      $BBB(5,2) > 10^{10^{286574}}$ (Shawn Ligocki, April 2022).
      See analysis by Shawn ligocki in

  \verb+https://www.sligocki.com/2022/04/03/mother-of-giants.html+\\
  
\item $BBB(2,3) \ge 59$ (Nicholas Drozd, September 2020)
\item $BBB(3,3) > 7.2 \times 10^{62}$ (Shawn Ligocki, February 2022)
      See analysis by Shawn Ligocki in
  
  \verb+https://www.sligocki.com/2022/02/27/bb-recurrence-relations.html+
  
\item $BBB(2,4) > 1.3 \times 10^{12}$ (Nicholas Drozd, January 2022)\\
      $BBB(2,4) > 6.7 \times 10^{16}$ (Nicholas Drozd, January 2022)\\
      $BBB(2,4) > 2.0 \times 10^{23}$ (Shawn Ligocki, February 2022)
\end{itemize}

\bigskip

Nicholas Drozd said that it is not difficult to prove that
$BBB(3,2) = 55$ and $BBB(2,3) = 59$.

\section{The methods}

The machines presented in this paper were discovered by means of computer programs. 
These programs contain procedures that achieve the following tasks:
\begin{enumerate}
\item To enumerate Turing machines without repetition.
\item To simulate Turing machines efficiently.
\item To recognize non-halting Turing machines.
\end{enumerate}

Note that these procedures are often mixed together in
real programs as follows: A tree of transition tables is generated,
and, as soon as some transitions are defined, the corresponding
Turing machine is simulated.
If the definition of a new transition is necessary, the tree is extended.
If the computation seems to loop, a proof of this fact is provided.

If the purpose is to prove a value for the busy beaver functions,
then all Turing machines in a class have to be studied. The machines
that pass through the three procedures above are either halting machines,
from which the better one is selected, or holdouts waiting for better
programs or for hand analyses.

If the purpose is to find lower bounds, a systematic enumeration
of machines is not necessary. Terry and Shawn Ligocki said they
used simulated annealing to find some of their machines.

The following references can be consulted for more information:
\begin{itemize}
\item Brady (1983) and Machlin and Stout (1990) for (4,2)-TM,
\item Marxen and Buntrock (1990) and Hertel (2009) for (5,2)-TM,
\item Lafitte and Papazian (2007) for (2,3)-TM,
\item Page about Macro Machines on Marxen's website. See

  \verb+http://turbotm.de/~heiner/BB/macro.html+

\item Harland (2016) and Harland (2022). 
  
\end{itemize}

\bigskip

\section{Busy beavers and unprovability}

\subsection{The result}

Let $S(n) = S(n,2)$ be Rado's busy beaver function. We know that $S(2) = 6$,
$S(3) = 21$, $S(4) = 107$, and we can hope to prove that $S(5) = 47,176,870$.
As we will see below, the fact that the busy beaver function $S$
is not computable implies that it is not possible to prove
that, for any natural number $n$, $S(n)$ has its true value.

Formally, we have the following theorem.

\smallskip

\noindent {\bf Theorem}. {\em Let $T$ be a well-known mathematical theory such as Peano arithmetic (PA)
or Zermelo-Fraenkel set theory with axiom of choice (ZFC). Then there exist numbers
$N$ and $L$ such that $S(N) = L$, but the sentence ``$S(N) = L$'' is not provable in $T$}.

\smallskip

This theorem is an easy consequence of the following proposition.

\smallskip

\noindent {\bf Proposition}. {\em Let $T$ be a well-known mathematical theory such as PA or ZFC.
Then there exists a Turing machine with two symbols $M$ that does not stop when
it is launched on a blank tape, but the fact that it does not stop is not provable in $T$}.

\smallskip

{\it Proof of the theorem from the proposition}. Let $M$ be the Turing machine given
by the proposition, let $N$ be the number of states of $M$, and let $L = S(N)$.
Then, to prove that ``$S(N) = L$'', we have to prove that $M$ does not stop.
But, by the proposition, such a proof does not exist.

\smallskip

Note that, if ``$S(N) = L$'' is a true sentence unprovable in theory $T$,
then, for all $m > L$, ``$S(N) < m$'' is also a true sentence unprovable
in theory $T$.

In the following, we consider many kinds of proofs of the proposition
and of the theorem.

\subsection{A direct proof}

This proposition is well-known and a one line proof can be given, as follows.

\smallskip

\noindent{\it Proof}. If all non-halting machines were provably non-halting,
then an algorithm that gives simultaneously the computable enumeration
of the halting machines and the computable enumeration of the provably
non-halting machines would solve the halting problem on a blank tape.

\smallskip

We give a detailed proof for nonspecialist readers.

\smallskip

\noindent{\it Detailed proof}. Let $M_1, M_2,\ldots$ be a computably
enumerable sequence of all Turing machines with two symbols.
Such a sequence can be obtained as follows: we list machines according to their
number of states, and, inside the set of machines with $n$ states, we list
the machines according to the alphabetical order of their transition tables.

Let $T_1, T_2,\ldots$ be a computably enumerable sequence of the
theorems of the theory $T$.
The existence of such a sequence is the main requirement that theory $T$ has to
satisfy in order that the proposition holds, and of course such a sequence exists
for well-known mathematical theories such as PA or ZFC.

Now consider the following algorithms $A$ and $B$.

\noindent{\bf Algorithm A}. We launch the machines $M_i$ on the blank tape as follows:
\begin{itemize}
\item one step of computation of $M_1$,
\item 2 steps of computation of $M_1$, 2 steps of computation of $M_2$, 
\item 3 steps of computation of $M_1$, 3 steps of computation of $M_2$,
3 steps of computation of $M_3$,
\item $\ldots$
\end{itemize}
When a machine $M_i$ stops, we add it to a list of machines that stop when
they are launched on a blank tape.

Note that, given a machine $M$, by running Algorithm $A$ we will know that $M$ stops if
$M$ stops, but we will never know that $M$ doesn't stop if $M$ doesn't stop.

\smallskip

\noindent{\bf Algorithm B}. We launch the algorithm that provides the computably enumerable
sequence of theorems of theory $T$, and each time we get a theorem $T_i$,
we look and see if this is a theorem of the form ``The Turing machine $M$ does not
stop when it is launched on a blank tape''. If that is the case, we add
$M$ to a list of Turing machines that provably do not stop on a blank tape.

\smallskip

Note that, given a machine $M$, by running Algorithm $B$ we will know that $M$ is
provably non-halting if $M$ is provably non-halting, but we will never know that $M$
is not provably non-halting if $M$ is not provably non-halting.

\smallskip

Now we have two algorithms, $A$ and $B$, and
\begin{itemize}
\item
Algorithm $A$ gives us a computably enumerable list of the Turing machines that stop
when they are launched on a blank tape.
\item
Algorithm $B$ gives us a computably enumerable list of the Turing machines that
provably do not stop on a blank tape.
\end{itemize}

We mix together these two algorithms, by a procedure called dovetailing, to get
Algorithm $C$, as follows.

\smallskip

\noindent{\bf Algorithm C}.
\begin{itemize}
\item one step of Algorithm $A$, one step of Algorithm $B$,
\item 2 steps of Algorithm $A$, 2 steps of Algorithm $B$,
\item 3 steps of Algorithm $A$, 3 steps of Algorithm $B$,
\item $\ldots$
\end{itemize}

Algorithm $C$ gives us simultaneously both the computably enumerable lists
provided by Algorithm $A$ and Algorithm $B$.

So Algorithm $C$ gives us both the list of halting Turing machines and the list of
provably non-halting Turing machines (on a blank tape).

Now we are ready to prove the proposition.
If all non-halting Turing machines were provably non-halting, then Algorithm $C$
would give us the list of halting Turing machines and the list of non-halting
Turing machines (on a blank tape).
So, given a Turing machine $M$, by running Algorithm $C$, we would see $M$ appearing
in one of the lists, and we could settle the halting problem for machine $M$
on a blank tape. So Algorithm $C$ would give us a computable procedure to
settle the halting problem on a blank tape.
But it is known that such a computable procedure does not exist.
Thus, there exists a non-halting Turing machine that is not provably
non-halting on a blank tape.

\subsection{The proposition as a special case of a general result}

The proposition is a special case of the following theorem.

\smallskip

\noindent{\bf Theorem}. {\em Let $A$ be a set of natural numbers that is computably enumerable
but not computable, and let $T$ be a well-known mathematical theory such as PA or ZFC.
Then there exists a natural number $n$ such that the sentence ``$n$ is not a member of $A$'' is
true but not provable in theory $T$}.

\smallskip

\noindent{\it Proof}. Since $A$ is computably enumerable, there exists an algorithm
that enumerates the natural numbers in $A$. If all natural numbers not in $A$
were provably not in $A$, then, by enumerating the proofs of theorems of theory
$T$, we would get an algorithm that enumerates the natural numbers not in $A$.
By running simultaneously both these algorithms, we could get a procedure that
decides membership in $A$, contradicting the fact that $A$ is not computable.  

\smallskip

The proposition is obtained from this theorem by numbering the list of Turing
machines, and by defining $A$ as the set of numbers of Turing machines
that stop on a blank tape.

\subsection{Some theoretical examples of Turing machines that satisfy the proposition}

Consider the Turing machine $M$ given by the proposition: $M$ does not stop when
it is launched on a blank tape, but this fact is not provable in theory $T$.
Can we get an idea of what such a machine $M$ looks like?
We give below some examples of such a Turing machine.

\subsubsection{Example 1: Using G\"odel's Second Incompleteness Theorem}

Let $M$ be a machine that enumerates the theorems of theory $T$,
and stops when it finds a contradiction (such as 0 = 1 if $T$ is Peano
arithmetic).

Then a proof within theory $T$ that $M$ does not stop would be a proof within
theory $T$ of the consistency of $T$, which is impossible by G\"odel's
Second Incompleteness Theorem (if theory $T$ is consistent).

\subsubsection{Example 2: Using G\"odel's First Incompleteness Theorem}

Another example can be given using G\"odel's First Incompleteness Theorem.
If $T$ is PA or ZFC, supposed to be consistent, the proof of this theorem provides
a formula $F$ that asserts its own unprovability. Thus $F$ is true, but
unprovable within theory $T$.

Consider the machine $M$ that enumerates the theorems of theory $T$, and stops when
it finds formula $F$.
Machine $M$ does not stop, since $F$ is unprovable, but a proof that it does not stop
would be a proof that $F$ is unprovable, so, since $F$ {\it is} ``$F$ is
unprovable'', a proof of $F$, which is impossible, since $F$ is unprovable. 

\subsubsection{Example 3: Using the Recursion Theorem}

As a third example, consider the machine $M$ that enumerates the theorems
of theory $T$ (PA or ZFC, supposed to be consistent), and stops when it finds
a formula $F$ that says that $M$ itself does not stop. Such a machine can be proved
to exist by applying the Recursion Theorem to the function $f$ such that
machine $M_{f(x)}$ stops if it finds a proof that machine $M_x$
does not stop.

Then $F$ is true, because, if $F$ were false, then $M$ would stop, so $F$ would be
a theorem of $T$, so $F$ would be true. But $F$ is unprovable, because since $F$ is true,
$M$ does not stop, so $F$ is not a theorem of theory $T$.
So the fact that $M$ does not stop is true and unprovable.

\subsection{Some explicit examples of Turing machines that satisfy the proposition}

Since May 2016, there are explicit constructions of Turing machines
whose behaviors are independent of ZFC. These machines never halt
on a blank tape, but this fact cannot be proved in ZFC.

\subsubsection{Example 1: Yedidia and Aaronson's machine}

Adam Yedidia and Scott Aaronson gave, in May 2016, a Turing machine
with 7910 states and two symbols such as it cannot be proved
in ZFC that it never halts. They note that enumerating
the theorems of ZFC would need a big number of states.
They use a graph theoretic statement that Harvey Friedman proved to be
equivalent to the consistency of a theory that implies the consistency of ZFC.
By using a new high-level language that is easily compiled down to
Turing machine description, they build a machine that would halt
if it finds a counterexample to Friedman's statement. 
See Yedidia and Aaronson (2016).

\subsubsection{Example 2: O'Rear's machine} 

S. O'Rear improved the number of states to 1919, in September 2016.
He improved later the number of states to 748.
His machines enumerate the theorems of a formal system
which has the same power as ZFC. See

\verb+https://github.com/sorear/metamath-turing-machines+

\noindent For a general presentation, see also Scott Aaronson's blog,
available at

\verb+http://www.scottaaronson.com/blog/?p=2725+

\subsection{A proof using Kolmogorov complexity}

There is another proof of unprovability, based on Kolmogorov complexity.
The Kolmogorov complexity of a number is the length of the shortest program
from which a universal Turing machine can output this number.
By Chaitin's Incompleteness Theorem, for any well-known mathematical
theory $T$, there exists a number $n(T)$ such that, for all numbers of
complexity greater than $n(T)$, the fact that they have complexity greater than
$n(T)$ is true but unprovable within theory $T$.

\smallskip

Chaitin's theorem also applies to the complexity defined as follows:
The complexity of a number $k$ is the smallest number $n$ of states
of a Turing machine with $n$ states and two symbols that outputs this number $k$,
written as a string of $k$ symbols 1, when the machine is launched on a blank tape.

So there exists a number $n(T)$ such that, for any number $k$ of
complexity greater than $n(T)$, the sentence ``the complexity of $k$ is
greater than $n(T)$'' is true but unprovable within theory $T$.
But ``$k > \mathnormal{\Sigma}(n(T))$'' implies
``the complexity of $k$ is greater than $n(T)$'',
so, for any number $k > \mathnormal{\Sigma}(n(T))$, the sentence
``$k > \mathnormal{\Sigma}(n(T))$'' is true but unprovable within theory $T$. 

For more details, see  Chaitin (1987),
Boolos, Burgess and Jeffrey (2002), p.\ 230, who note that
$n(T) < 10\uparrow\uparrow 10$, a stack of 10 powers of 10,
and Lafitte (2009).

\bigskip

\begin{Large}
{\bf Acknowledgments}
\end{Large}

This historical survey would not have been possible without the contributions
of many people. In particular, many thanks to Scott Aaronson,
Norbert B\'atfai, Allen Brady, Daniel Briggs, Nicholas Drozd,
Serge Grigorieff, James Harland, Tim Hutton, Pavel Kropitz,
Gr\'egory Lafitte, Shawn Ligocki, Terry Ligocki, Maurice Margenstern,
Heiner Marxen, Armando Matos, Robert Munafo,
Christophe Papazian, Holger Petersen, Myron Souris, H.J.M. Wijers,
and Hector Zenil.

\bigskip

\begin{Large}
{\bf References}
\end{Large}

\begin{enumerate}

\item
Aaronson S. (2020)\\
The busy beaver frontier\\
{\it ACM SIGACT News \/} {\bf 51} (3), September 2020, 32--54.

\item
Ben-Amram A.M., Julstrom B.A. and Zwick U. (1996)\\
A note on busy beavers and other creatures\\
{\it Mathematical Systems Theory \/} {\bf 29} (4), July-August 1996, 375--386.

\item
Ben-Amram A.M. and Petersen H. (2002)\\
Improved bounds for functions related to busy beavers\\
{\it Theory of Computing Systems \/} {\bf 35} (1), January-February 2002, 1--11. 

\item
Boolos G.S., Burgess J.P. and Jeffrey R.C. (2002)\\
Computability and Logic, 4th Ed.\\
Cambridge, 2002.

\item
Boolos G.S.\ and Jeffrey R.C. (1974)\\
Computability and Logic\\
Cambridge, 1974.

\item
Brady A.H. (1964)\\
Solutions to restricted cases of the halting problem \\
Ph.D. Thesis, Oregon State University, Corvallis, December 1964.

\item
Brady A.H. (1965)\\
Solutions of restricted cases of the halting problem used to
determine particular values of a non-computable function (abstract)\\
{\it Notices of the AMS \/} {\bf 12} (4), June 1965, 476--477.

\item
Brady A.H. (1966)\\
The conjectured highest scoring machines for Rado's $\mathnormal{\Sigma}(k)$
for the value $k = 4$\\
{\it IEEE Transactions on Electronic Computers\/}, EC-{\bf 15}, October 1966, 802--803.

\item
Brady A.H. (1974)\\
{\it UNSCC Technical Report\/} 11-74-1, November 1974.

\item
Brady A.H. (1975)\\
The solution to Rado's busy beaver game is now decided for
$k$ = 4 (abstract)\\
{\it Notices of the AMS \/} {\bf 22} (1), January 1975, A-25.

\item
Brady A.H. (1983)\\
The determination of the value of Rado's noncomputable function
$\mathnormal{\Sigma}(k)$ for four-state Turing machines\\
{\it Mathematics of Computation \/} {\bf 40} (162), April 1983, 647--665.

\item
Brady A.H. (1988)\\
The busy beaver game and the meaning of life\\
in: The Universal Turing Machine: A Half-Century Survey,
R. Herken (Ed.), Oxford University Press, 1988, 259--277.

\item
Buro M. (1990)\\
Ein Beitrag zur Bestimmung von Rados $\mathnormal{\Sigma}(5)$ oder Wie f\"{a}ngt man flei{\ss}ige Biber? (German)\\
{[}A contribution to the determination of Rado's $\mathnormal{\Sigma}(5)$, or: How to catch a busy beaver?]\\
Technical Report 146, Rheinisch-Westf\"{a}lische Technische Hochschule, Aachen, November 1990\\
(available at \verb+https://skatgame.net/mburo/ps/diploma.pdf+).

\item
Chaitin G. (1987)\\
Computing the busy beaver function\\
in: Open Problems in Communication and Computation, Springer, 1987, 108-112.\\
Reprinted in: Information, Randomness and Incompleteness, World Scientifics, 2nd Ed., 1990\\
(available at \verb+http://logic.amu.edu.pl/images/6/6c/Bobr.pdf+).

\item
Dewdney A.K. (1984a)\\
Computer recreations\\
{\it Scientific American \/} {\bf 251} (2), August 1984, 10--17.

\item
Dewdney A.K. (1984b)\\
Computer recreations\\
{\it Scientific American \/} {\bf 251} (5), November 1984, 27.

\item
Dewdney A.K. (1985a)\\
Computer recreations\\
{\it Scientific American \/} {\bf 252} (3), March 1985, 19.

\item
Dewdney A.K. (1985b)\\
Computer recreations\\
{\it Scientific American \/} {\bf 252} (4), April 1985, 16.

\item
Green M.W. (1964)\\
A lower bound on Rado's sigma function for binary Turing machines\\
in: {\it Proceedings of the 5th IEEE Annual Symposium on Switching Circuit Theory
and Logical Design\/}, November 1964, 91--94.

\item
  Harland J.\ (2013)\\
  Busy beaver machines and the observant otter heuristic\\
  in: {\it Proceedings of the 19th Computing: The Australasian Theory
  Symposium (CATS 2013)}, January 2013.

\item
Harland J.\ (2016)\\
Busy beaver machines and the observant otter heuristic
(or how to tame dreadful dragons)\\
{\it Theoretical Computer Science} {\bf 646}, 20 September 2016, 61--85\\
(preprint available at 
\verb+https://arxiv.org/abs/1602.03228+).

\item
Harland J.\ (2022)\\
Generating candidate busy beaver machines (or how to build the zany zoo)\\
{\it Theoretical Computer Science} {\bf 922}, 24 June 2022, 368--394\\
(preprint available at
\verb+https://arxiv.org/abs/1610.03184+).

\item
Hertel J. (2009)\\
Computing the uncomputable Rado sigma function\\
{\it The Mathematica Journal\/} {\bf 11} (2), 2009, 270--283.

\item
Jones J.P. (1974)\\
Recursive undecidability -- an exposition\\
{\it The American Mathematical Monthly\/} {\bf 81} (7), September 1974, 724--738.

\item
Julstrom B.A. (1992)\\
A bound on the shift function in terms of the busy beaver function\\
{\it SIGACT News \/} {\bf 23} (3), Summer 1992, 100--106.

\item
Julstrom B.A. (1993)\\
Noncomputability and the busy beaver problem\\
{\it The UMAP Journal \/} {\bf 14} (1), Spring 1993, 39--74.

\item
Kopp R.J. (1981)\\
The busy beaver problem\\
M.A. Thesis, State University of New York, Binghamton, 1981.

\item
Korfhage R.R.\ (1966)\\
Logic and Algorithms: With applications to the computer and information sciences\\
Wiley, 1966.

\item
Lafitte G. (2009)\\
Busy beavers gone wild\\
in: The Complexity of Simple Programs 2008, EPTCS 1, 2009, 123-129.\\
(available at \verb+https://arxiv.org/abs/0906.3257v1+).

\item
Lafitte G. and Papazian C. (2007)\\
The fabric of small Turing machines\\
in: {\it Computation and Logic in the Real World,
Proceedings of the Third Conference on Computabiliy in Europe\/},
June 2007, 219--227
(available at
\verb+https://citeseerx.ist.psu.edu/+\\
\verb+viewdoc/download?doi=10.1.1.104.3021&rep=rep1&type=pdf#page=231+).

\item
Lee C. Y. (1963)\\
Busy beaver examples\\
in: {\it Automata Theory: Advanced Concepts in Information Processing Systems\/}.
Lectures given at the University of Michigan, Summer 1963, 18--21.

\item
Lin S. (1963)\\
Computer studies of Turing machine problems\\
Ph.D. Thesis, The Ohio State University, Columbus, 1963.

\item
Lin S. and Rado T. (1965)\\
Computer studies of Turing machine problems\\
{\it Journal of the ACM \/} {\bf 12} (2), April 1965, 196--212.

\item
Ludewig J.\ , Schult U.\ and Wankm\"{u}ller F.\ (1983)\\
Chasing the busy beaver. Notes and observations on a competition to find the 5-state busy beaver\\
Forschungsberichte des Fachbereiches Informatik Nr 159, Universit\"{a}t Dortmund, 1983 (63 pages)\\
(available at
\verb+http://elib.uni-stuttgart.de/bitstream/11682/8466/1/lud20.pdf+).

\item
Lynn D.S. (1972)\\
New results for Rado's sigma function for binary Turing machines\\
{\it IEEE Transactions on Computers \/} C-{\bf 21} (8), August 1972, 894--896.

\item
Machlin R. and Stout Q.F. (1990)\\
The complex behavior of simple machines\\
{\it Physica D \/} {\bf 42}, June 1990, 85--98.

\item
Marxen H. and Buntrock J. (1990)\\
Attacking the Busy Beaver 5\\
{\it Bulletin of the EATCS \/} No 40, February 1990, 247--251\\
(available at
\verb+https://turbotm.de/~heiner/BB/mabu90.html+).

\item
Michel P. (1993)\\
Busy beaver competition and Collatz-like problems\\
{\it Archive for Mathematical Logic \/} {\bf 32} (5), 1993, 351--367.

\item
Michel P. (2004)\\
Small Turing machines and generalized busy beaver competition\\
{\it Theoretical Computer Science \/} {\bf 326} (1-3), October 2004, 45--56.

\item
Michel P. (2010)\\
Homology of groups and third busy beaver function\\
{\it International Journal of Algebra and Computation\/} {\bf 20} (6), 2010, 769--791.

\item
Michel P. (2015)\\
Problems in number theory from busy beaver competition\\
{\it Logical Methods in Computer Science \/} {\bf 11} (4:10), 2015, 1--35\\
(available at
\verb+https://arxiv.org/pdf/1311.1029.pdf+).

\item
Nabutovsky A. and Weinberger S. (2007)\\
Betti numbers of finitely presented groups and very rapidly growing functions\\
{\it Topology\/} {\bf 46} (2), March 2007, 211--223.

\item
Oberschelp A., Schmidt-G\"{o}ttsch K. and Todt G. (1988)\\
Castor quadruplorum\\
{\it Archive for Mathematical Logic\/} {\bf 27} (1), 1988, 35--44.

\item
Petersen H. (2006)\\
Computable lower bounds for busy beaver Turing machines\\
{\it Studies in Computational Intelligence} {\bf 25}, 2006, 305--319.

\item
  Petersen H. (2017)\\
  Busy beaver scores and alphabet size\\
  in: {\it Proceedings of the 21st International Symposium
    on Fundamentals of Computation Theory, FCT 2017},
  Lecture Notes in Computer Science Vol.\ 10472, Springer, 2017,
  409--417\\
  (available at
  \verb+https://arxiv.org/abs/1704.08752+).

\item
Rado T. (1962)\\
On non-computable functions\\
{\it Bell System Technical Journal\/} {\bf 41} (3), May 1962, 877--884.

\item
Rado T. (1963)\\
On a simple source for non-computable functions\\
in: {\it Proceedings of the Symposium on Mathematical
Theory of Automata, April 1962\/}, Polytechnic Institute of Brooklyn,
April 1963, 75--81.

\item
Soare R.I. (1996)\\
Computability and recursion\\
{\it Bulletin of Symbolic Logic\/} {\bf 2}, 1996, 284--321.

\item
Soare R.I. (2007)\\
Computability and incomputability\\
in: {\it Proceedings of the 3rd CiE}, LNCS 4497, Springer, 2007, 705--715.

\item
Soare R.I. (2009)\\
Turing oracle machines, online computing, and three displacements
in computability theory\\
{\it Annals of Pure and Applied Logic} {\bf 160} (3), September 2009,
368--399.

\item
Walsh T.R.S.\ (1982)\\
The busy beaver on a one-way infinite tape\\
{\it SIGACT News} {\bf 14} (1), Winter 1982, 38--43.

\item
Wang K. and Xu S. (1995)\\
New relation between the shift function and the busy beaver function\\
{\it Chinese Journal of Advanced Software Research\/} {\bf 2} (2),
1995, 192--197.

\item
Weimann B.\ (1973)\\
Untersuchungen \"{u}ber Rado's Sigma-Funktion und eingeschr\"{a}nkte Halteprobleme bei Turingmaschinen (German)\\
{[}Studies about Rado's Sigma function and limited halting problems of Turing machines]\\
Dissertation, Rheinische Friedrich-Wilhelms-Universit\"{a}t Bonn, November 1973. 

\item
Weimann B., Casper K.\ and Fenzl W.\ (1973)\\
Untersuchungen \"{u}ber haltende Programme f\"{u}r Turing-Maschinen mit
2 Zeichen und bis zu 5 Befehlen (German)\\
{[}Studies about halting programs for Turing machines with 2 symbols and up to
5 states]
in: {\it Deussen P.\ (eds) GI. Gesellschaft f\"{u}r Informatik e.\ V.\ 
  2.\ Jahrestagung (Karlsruhe, October 1972)},
Lecture Notes In Economics and Mathematical Systems, Vol.\ 78, Springer, 1973,
72--81.

\item 
Yang R., Ding L. and Xu S. (1997)\\
Some better results estimating the shift function in terms of
busy beaver function\\
{\it SIGACT News \/} {\bf 28} (1), March 1997, 43--48.

\item
Yedidia A.\ and Aaronson S.\ (2016)\\
A relatively small Turing machine whose behavior is independent of set theory\\
{\it Complex Systems} {\bf 25} (4), 2016, 297--327\\
(available at \verb+https://www.complex-systems.com/pdf/25-4-4.pdf+).

\end{enumerate}

\bigskip

\begin{Large}
{\bf Websites}
\end{Large}

\begin{enumerate}
\item
Scott Aaronson: \verb+https://www.scottaaronson.com/writings/bignumbers.html+\\
(how to write big numbers with few symbols).

\item
Joan Baez: \verb+https://johncarlosbaez.wordpress.com/2016/05/21/the-busy-beaver-game+\\
(a blog).

\item
Allen H. Brady: \verb+https://www.cse.unr.edu/~al/BusyBeaver.html+\\
(Busy Beaver Competition for machines with 3 states and 3 symbols).


\item
Georgi Gochev (Skelet): \verb+http://skelet.ludost.net/bb+\\
(Busy Beaver Competition for machines with 5 states and 2 symbols, and for
reversal machines).

\item
Google group about busy beavers: \verb+https://groups.google.com/g/busy-beaver-discuss+\\
(created and managed by Shawn Ligocki).

\item
Googology-wiki: \verb+https://googology.fandom.com/wiki/Busy_beaver_function+\\
(a good introduction).

\item
James Harland: \verb+https://sites.rmit.edu.au/jah/+\\
(definitions of some inverse functions for the busy beaver functions).

\item
Tim Hutton: \verb+https://github.com/GollyGang/ruletablerepository/wiki/+\\
\verb+TwoDimensionalTuringMachines+\\
(two-dimensional and higher-dimensional Turing machines)


\item
Heiner Marxen: \verb+http://turbotm.de/~heiner/BB+\\
(the top web reference on busy beavers).

\item
Pascal Michel: \verb+https://webusers.imj-prg.fr/~pascal.michel/+\\
(from which the present paper is written).



\item
Wikipedia: \verb+https://en.wikipedia.org/wiki/Busy_beaver+\\
(a good introduction).

\item
``Wythagoras'': \verb+https://googology.fandom.com/wiki/User:Wythagoras/Rado's_sigma_function+\\
(many lower bounds of $\mathnormal{\Sigma}(n,m)$ for high values of $n$ and $m$).

\item
Hector Zenil: \verb+https://demonstrations.wolfram.com/BusyBeaver+\\
\verb+https://hectorzenil.net+\\
(demonstration on the Wolfram Demonstrations Project, and website).

\end{enumerate}

\end{document}